\documentclass[aop,preprint]{imsart}

\RequirePackage[OT1]{fontenc}
\RequirePackage{amsthm,amsmath,natbib}
\RequirePackage[colorlinks]{hyperref}
\RequirePackage{hypernat}

\RequirePackage[OT1]{fontenc}
\RequirePackage{amsthm,amsmath,natbib}
\RequirePackage[colorlinks]{hyperref}
\RequirePackage{hypernat}

\numberwithin{equation}{section}

\usepackage{amssymb}
\usepackage{comment}
\usepackage{appendix}
\usepackage{graphicx}
\usepackage{subfigure}
\usepackage{enumerate}

\input xy
\xyoption{all}

\startlocaldefs

\newtheorem{theorem}{Theorem}[section]
\newtheorem{lemma}[theorem]{Lemma}
\newtheorem{proposition}[theorem]{Proposition}

\newcommand{\B}{\mathbf{B}}
\newcommand{\C}{\mathbf{C}}
\newcommand{\D}{\mathbf{D}}
\newcommand{\E}{\mathbf{E}}

\newcommand{\h}{\mathbf{H}}

\newcommand{\N}{\mathbf{N}}
\newcommand{\Z}{\mathbf{Z}}
\newcommand{\p}{\mathbf{P}}
\newcommand{\Q}{\mathbf{Q}}
\newcommand{\R}{\mathbf{R}}

\newcommand{\U}{\mathbf{U}}

\newcommand{\CC}{\mathcal {C}}

\newcommand{\CE}{\mathcal {E}}
\newcommand{\CF}{\mathcal {F}}

\newcommand{\CL}{\mathcal {L}}
\newcommand{\CM}{\mathcal {M}}
\newcommand{\CN}{\mathcal {N}}

\newcommand{\CV}{\mathcal {V}}
\newcommand{\CX}{\mathcal {X}}

\newcommand{\CZ}{\mathcal {Z}}

\newcommand{\CG}{\mathcal {G}}
\newcommand{\CH}{\mathcal {H}}

\newcommand{\dist}{{\rm dist}}

\newcommand{\diam}{{\rm diam}}
\newcommand{\var}{{\rm Var}}
\newcommand{\im}{{\rm Im}}

\newcommand{\cov}{{\rm Cov}}

\newcommand{\supp}{{\rm supp }}

\newcommand{\one}{\mathbf{1}}

\newcommand{\SEGGS}{\texttt{SEGGS}}

\newcommand{\ol}{\overline}

\newcommand{\wt}{\widetilde}
\newcommand{\wh}{\widehat}

\endlocaldefs

\begin{document}
\begin{frontmatter}

\title{Fluctuations for the Ginzburg-Landau $\nabla \phi$ Interface Model on  a Bounded Domain}
\runtitle{Fluctuations for the Ginzburg-Landau $\nabla \phi$ interface Model}

\begin{aug}

\author{\fnms{Jason} \snm{Miller}\thanksref{t2}\ead[label=e2]{jmiller@math.stanford.edu}}

\thankstext{t2}{Research supported in part by NSF grants DMS-0406042 and DMS-0806211.}

\runauthor{Jason Miller}

\affiliation{Stanford University}

\address{Stanford University\\
Department of Mathematics\\
Stanford, CA 94305\\
\printead{e2}\\
\phantom{E-mail:\ jmiller@math.stanford.edu}}

\end{aug}
\date{\today}

%% What is left to do:

%% Edit the introduction
%% Check everything one more time

\begin{abstract}
We study  the massless field on $D_n = D \cap \tfrac{1}{n} \Z^2$, where $D \subseteq \R^2$ is a bounded domain with smooth boundary, with Hamiltonian $\CH(h) = \sum_{x \sim y} \CV(h(x) - h(y))$.  The interaction $\CV$ is assumed to be symmetric and uniformly convex.  This is a general model for a $(2+1)$-dimensional effective interface where $h$ represents the height.  We take our boundary conditions to be a continuous perturbation of a macroscopic tilt: $h(x) = n x \cdot u + f(x)$ for $x \in \partial D_n$, $u \in \R^2$, and $f \colon \R^2 \to \R$ continuous.  We prove that the fluctuations of linear functionals of $h(x)$ about the tilt converge in the limit to a Gaussian free field on $D$, the standard Gaussian with respect to the weighted Dirichlet inner product $(f,g)_\nabla^\beta = \int_D \sum_i \beta_i \partial_i f_i \partial_i g_i$ for some explicit $\beta = \beta(u)$.  In a subsequent article, we will employ the tools developed here to resolve a conjecture of Sheffield that the zero contour lines of $h$ are asymptotically described by $SLE(4)$, a conformally invariant random curve.
\end{abstract}



\end{frontmatter}

\maketitle

\section{Introduction}
The object of our study is the massless field on $D_n = D \cap \tfrac{1}{n} \Z^2$, with Hamiltonian $\CH(h) = \sum_{b \in D_n^*} \CV(\nabla h(b))$.  Here, $D \subseteq \R^2$ is a bounded domain with smooth boundary.  The sum is over all edges in the induced subgraph of $\tfrac{1}{n} \Z^2$ with vertices in $D$ and $\nabla h(b) = h(y) - h(x)$ denotes the discrete gradient of $h$ across the oriented bond $b=(x,y)$.  We take our boundary conditions to be a continuous perturbation of a macroscopic tilt: $h(x) = nu \cdot x + f(x)$ when $x \in \partial D_n$ for $u \in \R^2$ and $f \colon \R^2 \to \R$ is a continuous function.  We consider a general interaction $\CV \in C^2(\R)$ which is assumed only to satisfy:
\begin{enumerate}
\item $\CV(x) = \CV(-x)$ (symmetry),
\item $0 < a_\CV \leq \CV''(x) \leq A_\CV < \infty$ (uniform convexity), and
\item $\CV''$ is $L$-Lipschitz.
\end{enumerate}
The purpose of the first condition is merely to simplify the notation since the symmetrization of a non-symmetric potential does not change its behavior.  Note that we can assume without loss of generality that $\CV(0) = 0$.  This is the so-called \emph{Ginzburg-Landau $\nabla \phi$ effective interface (GL) model}, also known as the \emph{anharmonic crystal}.  The variables $h(x)$ represent the heights of a random surface which serves as a model for an interface separating two pure phases.

The macroscopic behavior of the GL model has been the subject of much recent study.  An important step in this development is the construction and classification of \emph{Gibbs states}, which are infinite volume versions of the model.  In two-dimensions, it turns out that the height variable $h(x)$ diverges as the size of the domain tends to infinity, so in order to construct a Gibbs state one must first pass to the gradient field $\nabla h$.  The existence and uniqueness of \emph{gradient Gibbs states} is proved by Funaki and Spohn in \cite{FS97}, where they also study macroscopic dynamics.  Deuschel, Giacomin, and Ioffe in \cite{DGI00} establish a large deviations principle for the surface shape with zero boundary conditions but in the presence of a chemical potential and Funaki and Sakagawa in \cite{FS04} extend this result to the case of non-zero boundary conditions using the contraction principle.  The behavior of the maximum is studied by Deuschel and Giacomin in \cite{DG00} and by Deuschel and Nishikawa in \cite{DN07} in the case of Langevin dynamics.  The central limit theorem for linear functionals of the infinite gradient Gibbs states was proved first by Naddaf and Spencer for the static model with zero tilt in \cite{NS97} and Giacomin, Olla, and Spohn handle the time-varying case with general tilt evolving under Langevin dynamics in \cite{GOS01}.

The special situation in which the interaction is quadratic, i.e. $\CV(x) = \tfrac{1}{2} x^2$, corresponds to the so-called \emph{discrete Gaussian free field} (DGFF) or \emph{harmonic crystal}.  Many of the results in \cite{DG00, DN07, FS04, FS97, GOS01, NS97} have been established separately for the quadratic case and often the results in this setting are more refined.  The reason for the latter is that its Gaussian structure greatly simplifies its analysis; we will discuss this point in more detail later.  For the DGFF, large deviations principles for the surface shape as well as a central limit theorem for the height variable was proved by Ben Arous and Deuschel in \cite{BAD96}, the behavior of the maximum studied by Bolthausen, Deuschel, and Giacomin in \cite{BDG01} and by Daviaud in \cite{DAV06}.  In a particularly impressive and difficult work, Schramm and Sheffield in \cite{SS06} show that the macroscopic level sets are described by a family of conformally invariant random curves which are variants of $SLE(4)$.

Beyond having a Gaussian distribution, the main feature that makes the analysis of the DGFF tractable is that its mean and covariance are completely described in terms of the harmonic measure and Green's function associated with a simple random walk.  These objects are very well understood in the planar case \cite{LAW91}, which often allows for very precise estimates.  The mean and covariance of the more general GL model also admit a representation in terms of a random walk (Helffer-Sj\"ostrand representation or HS random walk \cite{HS94, DGI00}).  The general situation, however, is much more difficult because in addition to being non-Guassian, the corresponding representations involve a \emph{random walk in a dynamic random environment} whose behavior depends non-trivially on the boundary data.  The hypothesis that $\CV$ is uniformly convex is helpful in assuring that the jump rates of the HS random walk are uniformly bounded from above and below.  This implies that its Green's function is comparable to that of a simple random walk, which in turn allows for rough (up to multiplicative constants) variance estimates and, more generally, centered moments in terms of the corresponding moments for the DGFF (Brascamp-Lieb inequalities).  Furthermore, that the jump rates are bounded means that the Nash-Aronson and Nash continuity estimates apply, which give some rough control of the off-diagonal covariance structure.

\subsection{Main Results}

The main result of this article is the following central limit theorem for linear functionals of the height.  Let $D \subseteq \R^2$ be a bounded domain with smooth boundary.  For $\kappa \geq 0$, let $H^\kappa(D)$ be the Sobolev space of degree $\kappa$ and $H^{-\kappa}(D)$ its Banach space dual.

\begin{theorem}[Central Limit Theorem]
\label{thm::clt}
Let $f \colon \R^2 \to \R$ be a continuous function.  Suppose that $D_n = D \cap \tfrac{1}{n} \Z^2$ and that $h^n$ is distributed according to the GL model on $D_n$ with $h^n(x) = \varphi_n(x) + f(x)$ for all $x \in \partial D_n$ where $\varphi_n(x) = n u \cdot x$.  Let $a_u(b) = \E[ \CV''(\eta(b))]$ where $\eta$ has the law of the Funaki-Spohn state with tilt $u$.  Define the linear functional
\[ \xi_\nabla^{n,D}(g) = \sum_{b \in D_n^*} a_u(b) \nabla g(b) \nabla (h^n - \varphi_n)(b) \text{ for } g \in H^{\kappa}(D).\]
For every $\kappa > 4$, the law of $\xi_\nabla^{n,D}$ on $H^{-\kappa}(D)$ converges weakly with respect to the weak topology of $H^{-\kappa}(D)$ to a Gaussian free field on $D$, the standard Gaussian with respect to the weighted Dirichlet inner product $(g_1,g_2)_\nabla^\beta = \int_D \sum_i \beta_i \partial_i g_1 \partial_i g_2$, with boundary condition $f$, where $(\beta_1,\beta_2)$ is proportional to $(\E[\CV''( \eta(0,e_1))], \E[\CV''( \eta(0,e_2))])$.  In particular, $\beta_1 = \beta_2 > 0$ when $u = 0$.
\end{theorem}
\noindent We will recall both the notion of a Sobolev space and the GFF in subsection \ref{subsec::gff}.  We will also review the notion of a Funaki-Spohn state (as well as a new construction for zero-tilt) in subsection \ref{subsec::shift_invariant}.

As we mentioned earlier, central limit theorems have already been established by Naddaf and Spencer in \cite{NS97} and Giacomin, Olla, and Spohn in \cite{GOS01} for linear functionals of infinite gradient Gibbs states of the GL model.  Both articles are based on the beautiful observation that the CLT can be reduced to a homogenization problem using the HS representation.  The proof in \cite{GOS01} has more of a probabilistic flavor while the approach in \cite{NS97} is to use PDE techniques.  The reason that these results are restricted to infinite gradient Gibbs states is that the homogenization techniques they employ fail to carry over to the finite case.  In particular, the main step in \cite{GOS01} is a proof that the macroscopic covariance structure of the gradient Gibbs state for the GL model is the same as that in the GFF by showing that the HS random walk converges in the limit to a Brownian motion.  The key tool here is the so-called Kipnis-Varadhan method \cite{KV86}, which is to represent the random walk as an additive functional of the environment from the perspective of the walker.  When the environment is an infinite, stationary, ergodic Markov process then it remains so when viewed from the walker, thus the convergence to Brownian motion is a consequence of Corollary 1.5 of \cite{KV86}.  If the environment is finite and, in particular, \emph{not ergodic with respect to shifts}, this approach can no longer be used since \emph{the environment viewed from the particle is not ergodic}.

The covariance matrix of the limiting Gaussian in \cite{GOS01} is given in terms of a complicated variational formula.  It is therefore not explicit, except in the case of zero-tilt where it is possible to argue that it is proportional to the identity using rotational invariance.  The proof of Theorem \ref{thm::clt} gives the covariance matrix explicitly, up to a multiplicative constant, which is another new result for gradient Gibbs states.

The careful reader may note that the convergence in Corollary 2.2. of \cite{GOS01} is in $H^{-\kappa}(D)$ for $\kappa > 3$ while we require $\kappa > 4$.  The reason for the distinction is that we assume only continuity of the boundary condition $f$.  This forces us to perform an extra integration by parts, which in turn puts an extra derivative on the test function.  In the more restrictive setting of $C^1$ boundary conditions, our proof also gives convergence in $H^{-\kappa}(D)$ for $\kappa > 3$.

The main step in our proof is motivated by the Markovian structure enjoyed by the quadratic case: the law of a DGFF on $D_n$ with boundary condition $f_n$ is equal in law to that of a \emph{zero-boundary} DGFF on $D_n$ plus the discrete harmonic extension of $f_n$ to $D_n$.  This property is a higher dimensional analog of the fact that a random walk $X_t$ on $\Z$ conditioned to satisfy $X_{t_1} = x_1$ and $X_{t_2} = x_2$ for $t_1 < t_2$ has the law of $Y_t + H_t$ where $Y_t$ is a random walk on $I = [t_1,t_2]$ conditioned to vanish at $t_1,t_2$ and $H$ is the discrete harmonic function on $I$ with boundary values  $H(t_1) = X_{t_1}$ and $H_{t_2} = X_{t_2}$.  Our next theorem is a quantitative estimate of the degree to which this property approximately holds for the GL model.  Although we state it as our second theorem, it is the key step in the proof of Theorem \ref{thm::clt} and much of the article is dedicated to its proof.  In order to give a precise statement of this result, we first need to setup some notation.

Suppose that $D \subseteq \Z^2$ is a bounded subset of diameter $R > 0$.  Fix $\ol{\Lambda} > 0$ and let $\B_{\ol{\Lambda}}^u(D)$ be the set of functions $\phi \colon \partial D \to \R$ satisfying $\max_{x \in \partial D} | \phi(x) - u \cdot x| \leq \ol{\Lambda} (\log R)^{\ol{\Lambda}}$.  For $r > 0$, let $D(r) = \{ x \in D : \dist(x,\partial D) \geq r\}$.  With $\phi \in \B_{\ol{\Lambda}}^u(D)$,  let $\p_D^\phi$ denote the law of the GL model on $D$ with boundary condition $\phi$.  In other words, $\p_D^\phi$ is the measure on functions $h \colon D \to \R$ with density
\[ \frac{1}{\CZ} \exp\left( - \sum_{b \in D^*} \CV(\nabla (h \vee \phi)(b)) \right)\]
with respect to Lebesgue measure on $\R^{|D|}$.  Here, $h \vee \phi$ is used to denote the function
\begin{equation}
\label{intro::eqn::vee} h \vee \phi (x) = \begin{cases} h(x) \text{ if } x \in D,\\
						 \phi(x) \text{ if } x \in \partial D. \end{cases}
\end{equation}

We will write $O_{\ol{\Lambda}}(f(x))$ to denote the set of functions $g$ for which there exists a constant $c_{\ol{\Lambda}}$ depending on $\ol{\Lambda}$ but independent of $R$ so that $|g(x)| \leq c_{\ol{\Lambda}}|f(x)|$.  For $\beta = (\beta_1,\beta_2)$, let
\[ (\Delta^\beta f)(x) = 
	\beta_1 ( f(x+e_1) + f(x-e_1) - 2 f(x)) +
	\beta_2 ( f(x+e_2) + f(x-e_2) - 2 f(x))\]
where $e_1 = (1,0)$ and $e_2 = (0,1)$.
Note that $\Delta^\beta$ is the usual discrete Laplacian for $\beta = (1,1)$.

\begin{theorem}
\label{harm::thm::coupling}
Suppose that $\psi,\wt{\psi} \in \B_{\ol{\Lambda}}^u(D)$ and $\beta = \beta(u)$ as in the statement of Theorem \ref{thm::clt}.  There exists $C,\epsilon,\delta > 0$ depending only on $\CV$ such that if $r \geq CR^{1-\epsilon}$ then the following holds.  There exists a coupling $(h^\psi, h^{\wt{\psi}})$ of $\p_D^\psi, \p_D^{\wt{\psi}}$ such that if $\wh{h} \colon D(r) \to \R$ solves the elliptic problem $\Delta^\beta \wh{h} = 0$ with $\wh{h}|_{\partial D(r)} = \ol{h} = h^\psi - h^{\wt{\psi}}$ then 
\[ \p[ \ol{h} \neq \wh{h} \text{ in } D(r)] = O_{\ol{\Lambda}}(R^{-\delta}).\]
When $u = 0$, we can take $\beta = (1,1)$ so that $\Delta^\beta$ is the usual Laplacian.
\end{theorem}

One of the main challenges in the analysis of the GL model is the lack of useful comparison inequalities for its mean.  The difficulty is that the only explicit formula is given in terms of the annealed first exit distribution of the HS walk \cite{DGI00}.  It is not possible to extract any sort of asymptotic contiguity of this measure with respect to the harmonic measure of simple random walk using only that the HS walk jumps with bounded rates, which is all that is required to prove comparability of the corresponding Green's functions hence also of centered moments with DGFF.  Indeed, examples have been worked out in the continuum setting of diffusions in which the two measures are absolutely singular and that the support of the former has a fractal structure.  The situation is further complicated in the setting of the HS walk since in addition to being dynamic, its jump rates also depend on the boundary conditions, hence it seems difficult to rule out pathological behavior whenever the walk gets close to the boundary and the boundary conditions are rough.

Applying Theorem \ref{harm::thm::coupling} to the special case $\wt{\psi}(x) = u \cdot x$ gives the following estimate of the mean, which we believe to be sufficiently important that we state it as a separate theorem.

\begin{theorem}
\label{harm::thm::mean_harmonic}
Suppose that $\psi \in \B_{\ol{\Lambda}}^u(D)$.  There exists $C, \epsilon,\delta > 0$ such that if $r \geq CR^{1-\epsilon}$ and $\beta = \beta(u)$ as in the statement of Theorem \ref{thm::clt} then the following holds.  If $\wh{h} \colon D(r) \to \R$ is the $\Delta^\beta$-harmonic extension of $\E^\psi h$ from $\partial D(r)$ to $D(r)$ then 
\[ \max_{x \in D(r)} |\E^\psi h(x) - \wh{h}(x)| = O_{\ol{\Lambda}}(R^{-\delta}).\]
When $u = 0$, we can take $\beta = (1,1)$ so that $\wh{h}$ is harmonic with respect to the usual discrete Laplacian.
\end{theorem}

It is worth pointing out that both of these theorems place no restrictions on the regularity of the boundary conditions $\psi, \wt{\psi}$ nor the regularity of $\partial D$.

\begin{comment}
Recalling the relationship between the Green's function of the HS random walk and the covariance structure of the underlying field, Theorem \ref{thm::clt} implies that the Green's function of the HS random walk is the same as that of a Brownian motion.  Using this observation, we are able to prove the homogenization of the HS random walk.

\begin{theorem}
\label{thm::hs_homogenization}
Supposing we are in the setting of Theorem \ref{thm::clt} and that $X_t^n$ denotes the HS random walk associated with the dynamics of $h^n$.  Then $X_t^n$ stopped on its first exit from $D_n$ converges weakly to a standard Brownian motion on $D$ stopped on its first exit from $D$.
\end{theorem}
\end{comment}

\subsection{Sequel}  This article the first in a series of two and will be a prerequisite for the second.  In the sequel, we will make use of many of the estimates developed here in order to resolve a conjecture made by Sheffield (Problem 10.1.3 in \cite{SHE05}) that the macroscopic level lines of the GL model converge in the limit to $SLE(4)$; the case of quadratic potentials is proved by Schramm and Sheffield in \cite{SS06}.  The two papers together are meant to be fairly self-contained.

\subsection{Outline}  The remainder of the article is structured as follows.  The second section is a short discussion of discrete and continuum Gaussian free fields.  We chose to include the former part of this section since the special Markovian structure of the DGFF is the inspiration for Theorem \ref{harm::thm::coupling} and also to serve as an illustration of the complications associated with non-quadratic interaction.  In the latter part, we provide a brief description of the GFF, the standard Gaussian law on $H_0^1(D)$.  A much more thorough introduction can be found in \cite{SHE06}.  In Section \ref{sec::gl}, we will give a formal introduction to the GL model, its Langevin dynamics as well as the HS representation, and the Brascamp-Lieb inequalities.  In Section \ref{sec::dyn}, we will explain how the Langevin dynamics can be used to construct couplings of the GL model and prove an energy inequality for the discrete Dirichlet energy of such a coupling.  This section is concluded with an equivalence of ensembles result: the Funaki-Spohn shift-ergodic gradient Gibbs state can be realized as an infinite volume limit of models on finite domains.  In Section \ref{sec::harm}, we will prove Theorems \ref{harm::thm::coupling} and \ref{harm::thm::mean_harmonic} using an entropy estimate which is based on technical estimates from Section \ref{sec::correlation_decay}.  Finally, Section \ref{sec::clt} is relatively short and deduces the CLT from Theorem \ref{harm::thm::coupling}.  We conclude the article with two appendices containing useful estimates on discrete harmonic functions and symmetric random walks.

\section{Gaussian Free Fields}

In this section we will introduce the discrete and continuum Gaussian free fields (DGFF and GFF).  The reason that we include a discussion of the latter separate from the general case of the GL model is to emphasize its special Markovian structure, which is the motivation behind the ideas used in Section \ref{sec::harm}.

\subsection{Discrete Gaussian Free Field}
\label{subsec::dgff_construction}
Suppose that $G = (V \cup \partial V,E)$ is a finite, undirected, connected graph with distinguished subset $\partial V \neq \emptyset$ and edge weights $\omega > 0$.  The zero-boundary discrete Gaussian free field (DGFF) is the measure on functions $h \colon V \cup \partial V \to \R$ vanishing on $\partial V$ with density
\[ \frac{1}{\CZ_G} \exp\left(-\frac{1}{2} \sum_{b \in E} \omega(b) (\nabla (h \vee 0)(b))^2\right)\]
with respect to Lebesgue measure.  Here, $h \vee 0$ has the same meaning as in \eqref{intro::eqn::vee} and $\CZ_G$ is a normalizing constant so that the above has unit mass.  Equivalently, the DGFF is the standard Gaussian associated with the Hilbert space $H_0^1(V)$ of real-valued functions $h$ on $V$ vanishing on $\partial V$ with weighted Dirichlet inner product
\[ (f,g)_\nabla^\omega = \sum_{b \in E} \omega(b) \nabla f(b) \nabla g(b).\]
This means that the DGFF $h$ can be thought of as a family of Gaussian random variables $(h,f)_\nabla^\omega$ indexed by elements $f \in H_0^1(V)$ with mean zero and covariance
\begin{equation}
 \cov((h,f)_\nabla^\omega,(h,g)_\nabla^\omega) = (f,g)_\nabla^\omega,\ f,g \in H_0^1(V) \label{dgff::covariance}.
 \end{equation}
Although perhaps non-standard since our Hilbert space is finite dimensional, this representation is convenient since it allows for a simple derivation of the mean and covariance of $h$.  Let $\Delta^\omega \colon V \to \R$ denote the discrete Laplacian on $V$, i.e.
\[ \Delta^\omega f(x) = \sum_{b \ni x} \omega(b) \nabla f(b)\]
and let $G^\omega(x,y) = (\Delta^\omega)^{-1} \one_{\{x\}}(y)$ be the discrete Green's function on $V$.  Summation by parts gives that
\[ (f,g)_\nabla^\omega = -\sum_{x \in V} f(x) \Delta^\omega g(x) = -\sum_{x \in V} \Delta^\omega f(x) g(x) \text{ for } f,g \in H_0^1(V).\]
Thus 
\[h(x) = (h,\one_{\{x\}}(\cdot))_{L^2} = -(h,G^\omega(x,\cdot))_\nabla,\]
hence
\[ \cov(h(x),h(y)) = (G^\omega(x,\cdot),G^\omega(y,\cdot))_\nabla = G^\omega(x,y).\]

Suppose that $W \subseteq V$.  Then $H_0^1(V)$ admits the orthogonal decomposition $H_0^1(V) = \CM_I \oplus \CM_B \oplus \CM_O$ where $\CM_I,\CM_B,\CM_O$ are the subspaces of $H_0^1(V)$ consisting of those functions that vanish on $V \setminus W$, are $\Delta^\omega$-harmonic off of $\partial W$, and vanish on $W$, respectively.   It follows that we can write $h = h_I + h_B + h_O$ with $h_I \in \CM_I, h_B \in \CM_B, h_O \in \CM_O$ where $h_I,h_B,h_O$ are independent.  This implies that the DGFF possesses the following Markov property: the law of $h|_W$ conditional on $h|_{V \setminus W}$ is that of a zero boundary DGFF on $W$ plus the $\Delta^\omega$-harmonic extension of $h$ from $\partial W$ to $W$.  In particular, the conditional mean of $h|_W$ given $h|_{V \setminus W}$ is the $\Delta^\omega$-harmonic extension of $h|_{\partial W}$ to $W$.

More generally, if $\phi \colon \partial V \to \R$, the DGFF with boundary condition $\phi$ is the measure on functions $h \colon V \to \R$ with $h|_{\partial V} = \phi$ with density
\[ \frac{1}{\CZ_G} \exp\left( -\frac{1}{2} \sum_{b \in E} \omega(b)(\nabla (h \vee \phi)(b))^2 \right).\]
That is, $h$ has the law of a zero boundary DGFF on $V$ plus the $\Delta^\omega$-harmonic extension of $\phi$ from $\partial V$ to $V$.

\subsection{The Continuum Gaussian Free Field}
\label{subsec::gff}

The GFF is a $2$-time dimensional analog of the Brownian motion.  Just as the Brownian motion can be realized as the scaling limit of many random curve ensembles, the GFF arises as the scaling limit of a number of random surface ensembles \cite{BAD96, GOS01, KEN01, NS97, RV08}, as well as the model under consideration in this article.  In this subsection, we will describe the basic properties of the GFF necessary for our analysis.  Let $D$ be a bounded domain in $\R^2$ with smooth boundary and let $C_0^\infty(D)$ denote the set of $C^\infty$ functions compactly supported in $D$.  We begin with a short discussion of Sobolev spaces; the reader is referred to Chapter 5 of \cite{EVAN02} or Chapter 4 of \cite{TAY96} for a more thorough introduction.  With $\N_0 = \{0,1,\ldots\}$ the non-negative integers, when $f \in C_0^\infty(D)$ and $\alpha = (\alpha_1,\alpha_2) \in \N_0^2$ we let $D^\alpha f = \partial_1^{\alpha_1} \partial_2^{\alpha_2} f$.  For $k \in \N_0$ we define the $H^k(D)$-norm
\begin{equation}
\label{eqn::sobolev} 
\|f\|_{H^k(D)}^2 =  \sum_{|\alpha| \leq k} \int_D |D^\alpha f(x)|^2 dx
\end{equation}
where $|\alpha| = \alpha_1 + \alpha_2$.  The Sobolev space $H_0^k(D)$ is the Banach space closure of $C_0^\infty(D)$ under $\Vert \cdot \Vert_{H^k(D)}$.  If $s \geq 0$ is not necessarily an integer then $H_0^s(D)$ can be constructed via the complex interpolation of $H_0^0(D) = L^2(D)$ and $H_0^k(D)$ where $k \geq s$ is any positive integer (see Chapter 4 section 2 of \cite{TAY96} for more on this construction and also Chapter 4 of \cite{KAT04} for more on interpolation).  A consequence of this is that if $T \colon C_0^\infty(D) \to C_0^\infty(D)$ is a linear map continuous with respect to the $L^2(D)$ and $H^k(D)$ topologies then it is also continuous with respect to $H^s(D)$ for all $0 \leq s \leq k$.  For $s \geq 0$ we define $H^{-s}(D)$ to be the Banach space dual of $H_0^s(D)$ where the dual pairing of $f \in H^{-s}(D)$ and $g \in H_0^s(D)$ is given formally by the usual $L^2(D)$ inner product
\[ (f,g) = (f,g)_{L^2(D)} = \int_D f(x) g(x) dx.\]
More generally, for any $s \in \R$ the $H^s(D)$-topology can be constructed explicitly via the inner product
\begin{equation}
\label{gff::eqn::sobolev_inner_product}
 (f,g)_{s} = \int (1-\ol{\Delta})^{s/2} f \cdot (1-\ol{\Delta})^{s/2} g;
\end{equation}
see the introduction of Chapter 4 of \cite{TAY96}.
We are using the notation $\ol{\Delta}$ for the Laplacian on $\R^2$ to keep the notation consistent since elsewhere in the article $\Delta$ refers to the discrete Laplacian.
Here, 
\[ (1-\ol{\Delta})^p f = \CF^{-1} [(1+\xi_1^2 + \xi_2^2)^p (\CF f)] \text{ for } p \in \R\]
where 
\[  \CF f (\xi) = \int e^{-i \xi \cdot x} f(x) dx\]
is the Fourier transform of $f$.  We will be most interested in the space $H_0^1(D)$.  Fix a positive definite $2 \times 2$ real matrix $A$.  An application of the Poincare inequality (Chapter 4, Proposition 5.2) gives that the norm induced by the weighted Dirichlet inner product
\[ (f,g)_\nabla^A \equiv \int_D \sum_{i,j} a_{ij}  \partial_i  f \partial_j g \text{  for } f,g \in C_0^\infty(D)\]
is equivalent to $\| \cdot \|_{H^1(D)}$.  This choice of inner product is particularly convenient because it is invariant under precomposition by conformal transformations when $A$ is a multiple of the identity.

The $A$-GFF $h$ on $D$ can be expressed formally as a random linear combination of an $(\cdot,\cdot,)_\nabla^A$-orthonormal basis $(f_n)$ of $H_0^1(D)$
\[ h = \sum_n \alpha_n f_n\]
where $(\alpha_n)$ is an iid sequence of standard Gaussians.  Although the sum defining $h$ does not converge in $H_0^1(D)$,  for each $\epsilon > 0$ it does converge almost surely in $H^{-\epsilon}(D)$ (\cite[Proposition 2.7]{SHE06} and the discussion thereafter).  If $f,g \in C_0^\infty(D)$ then an integration by parts gives $(f,g)_\nabla^A = -( f,\ol{\Delta}^A g)$.  Here, $\ol{\Delta}^A = \ol{\nabla} A \ol{\nabla}$.  Using this, we define
\[ (h,f)_\nabla^A = -(h, \ol{\Delta}^A f) \text{  for } f \in C_0^\infty(D).\]
Observe that $(h,f)_\nabla^A$ is a Gaussian random variable with mean zero and variance $(f,f)_\nabla^A$.  Hence by polarization $h$ induces a map $C_0^\infty(D) \to \CG$, $\CG$ a Gaussian Hilbert space, that preserves the Dirichlet inner product.  This map extends uniquely to $H_0^1(D)$ which allows us to make sense of $(h,f)_\nabla^A$ for all $f \in H_0^1(D)$.  We are careful to point out, however, that while $(h, \cdot)_\nabla^A$ is well-defined off of a set of measure zero as a linear functional on $C_0^\infty(D)$ this is not the case for general $f \in H_0^1(D)$.

Suppose that $W \subseteq D$ is a smooth, open set.  Then there is a natural inclusion of $H_0^1(W)$ into $H_0^1(D)$ given by the extension by value zero.  If $f \in C_0^\infty(W)$ and $g \in C_0^\infty(D)$ then as $(f,g)_\nabla^A = -(f,\ol{\Delta}^A g)$ it is easy to see that $H_0^1(D)$ admits the $(\cdot,\cdot)_\nabla^A$-orthogonal decomposition $\CM \oplus \CN$ where $\CM = H_0^1(W)$ and $\CN$ is the set of functions in $H_0^1(D)$ that are $\ol{\Delta}^A$-harmonic on $W$.  Thus we can write
\[ h = h_W + h_{W^c} = \sum_n \alpha_n f_n + \sum_n \beta_n g_n\]
where $(\alpha_n),(\beta_n)$ are independent iid sequences of standard Gaussians and $(f_n)$, $(g_n)$ are orthonormal bases of $\CM$ and $\CN$, respectively.  Observe that $h_W$ has the law of the GFF on $W$, $h_{W^c}$ the  $\ol{\Delta}^A$-harmonic extension of $h|_{\partial W}$ to $W$, and $h_W$ and $h_{W^c}$ are independent.  We arrive at the following proposition:

\begin{proposition}[Markov Property]
\label{gff::prop::markov}
The conditional law of $h|_W$ given $h |_{ D \setminus W}$ is that of the $A$-GFF on $W$ plus the $\ol{\Delta}^A$-harmonic extension of the restriction of $h$ on $\partial W$ to $W$.
\end{proposition}

This proposition will be critical in the proof of Theorem \ref{thm::clt}.  It also allows us to make sense of the $A$-GFF with non-zero boundary conditions: if $f \colon \partial D \to \R$ is a continuous function and $F$ is its $\ol{\Delta}^A$-harmonic extension from $\partial D$ to $D$ then the law of the $A$-GFF on $D$ with boundary condition $f$ is given by the law of $F + h$ where $h$ is a zero boundary $A$-GFF on $D$.

\section{The Ginzburg-Landau Model}
\label{sec::gl}

The Ginzburg-Landau $\nabla \phi$-interface (GL) model is a general effective interface model first studied by Funaki and Spohn in \cite{FS97} and Naddaf and Spencer in \cite{NS97}.  Suppose that $G = (V \cup \partial V,E)$ is a finite, undirected, connected graph with a distinguished set of vertices $\partial V$.  Let $\CV \in C^2(\R)$ satisfy:
\begin{enumerate}
\item $\CV(x) = \CV(-x)$ (symmetry),
\item $0 < a_\CV \leq \CV''(x) \leq A_\CV < \infty$ (uniform convexity), and
\item $\CV''$ is $L$-Lipschitz.
\end{enumerate}
The law of the GL model on $V$ with potential function $\CV$ and boundary condition $\psi \colon \partial V \to \R$ is the measure on functions $h \colon V \to \R$ with $h|_{\partial V} = \psi$ described by the density
\[ \frac{1}{\CZ_\CV} \exp\left( - \sum_{b \in E} \CV(\nabla (h \vee \psi)(b)) \right) \]
with respect to Lebesgue measure and $h \vee \psi$ is as in \eqref{intro::eqn::vee}.

\subsection{Langevin Dynamics}

\begin{figure}
     \centering
          \includegraphics[width=.7\textwidth]{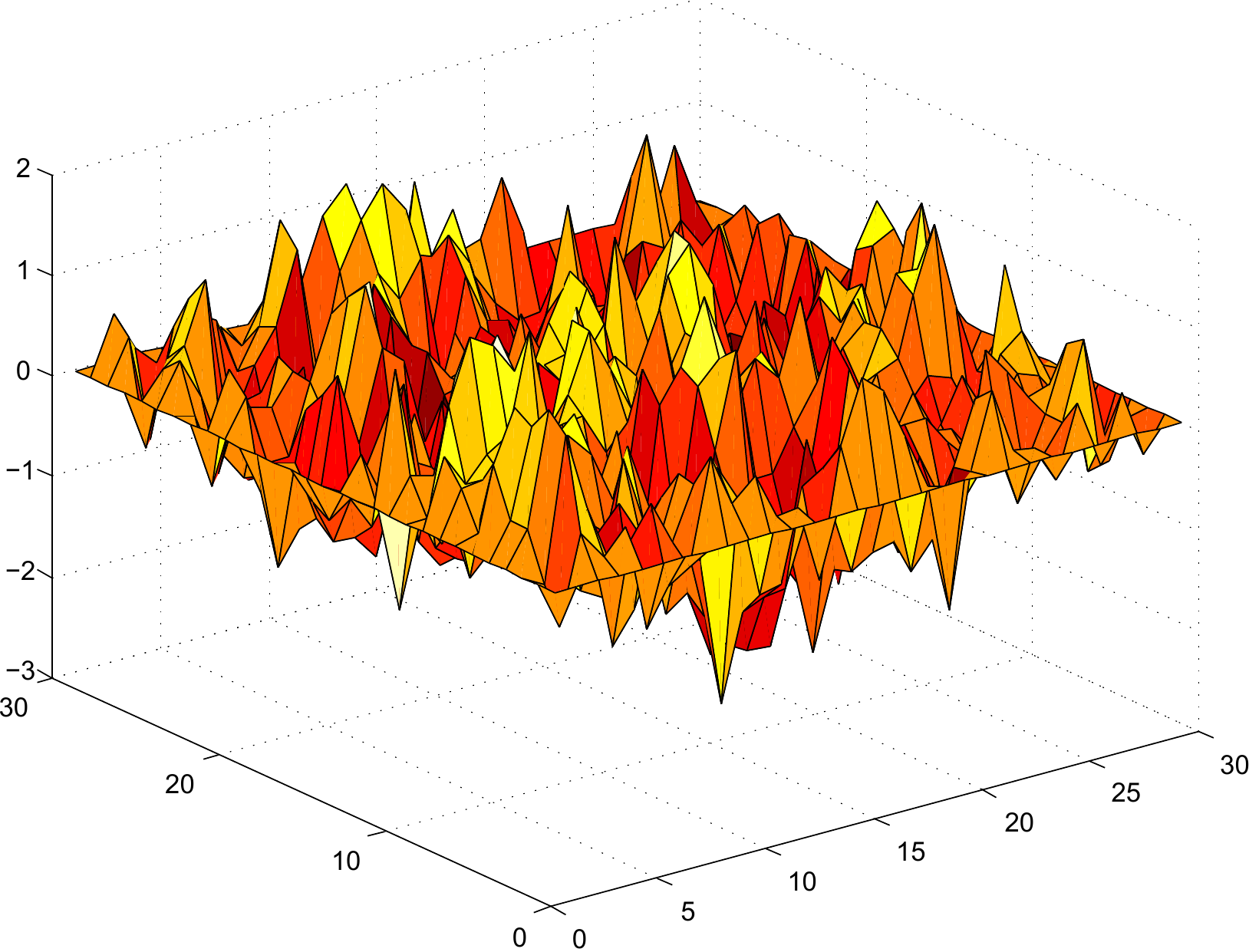}
          \caption{A typical realization of the GL model with zero boundary conditions on $\{0,\ldots,30\}^2$ and potential function $\CV(x) = x^2 + \cos(x)$, sampled using a discretization of the SDS \eqref{gl::eqn::dynam}}
\end{figure}

Consider the stochastic differential system (SDS)
\begin{equation}
\label{gl::eqn::dynam}
 dh_t^\psi(x) = \sum_{b \ni x} \CV'(\nabla (h_t^\psi \vee \psi)(b))dt + \sqrt{2}dW_t(x) \text{ for } x \in D
\end{equation}
where $W_t(x)$, $x \in V$, is a family of independent standard Brownian motions.  The generator for \eqref{gl::eqn::dynam} is given by
\begin{align*}
   \CL \varphi(h) 
&= \sum_{x \in V} \left(\partial_{h(x)}^2 \varphi(h) + \sum_{b \ni x} \CV'(\nabla h \vee \psi(b)) \partial_{h(x)}\varphi(h)\right)\\
&= \sum_{x \in V} e^{\CH_\CV^\psi(h)} \frac{\partial}{\partial h(x)} \left( e^{-\CH_\CV^\psi(h)} \frac{\partial}{\partial h(x)} \varphi(h) \right),
\end{align*}
where 
\[ \CH_\CV^\psi(h) = \sum_{b \in E} \CV( \nabla h \vee \psi(b))\]
is the Hamiltonian for the GL model.
Thus it is easy to see that $\CL$ is self-adjoint in the space $L^2(e^{-\CH(h)})$, hence the dynamics \eqref{gl::eqn::dynam} are reversible with respect to the law of the GL model.  These are the \emph{Langevin dynamics}.

\subsection{The Helffer-Sj\"ostrand Representation}
\label{subsec::hs_representation}

We showed in subsection \ref{subsec::dgff_construction} that if $\CV(x) = \tfrac{1}{2} x^2$ then the mean height is harmonic and that the covariance of heights is described by the discrete Green's function.  Both of these quantities admit simple probabilistic representations:  if $X_t$ is a continuous-time random walk (CTRW) on $G$ that jumps with uniform rate $1$ equally to its neighbors and $\tau$ is the time it first hits $\partial V$, then
\[ \E h(x) = \E_x h(X_\tau) \text{ and } \cov(h(x),h(y)) = \E_x \int_0^\tau \one_{\{X_s = y\}} ds\]
where the subscript $x$ indicates $X_0 = x$.
The idea of the Helffer-Sj\"ostrand (HS) representation, originally developed in \cite{HS94} and reworked probabilistically in \cite{DGI00, GOS01}, is to give an expression for the corresponding quantities for the GL model in terms of the first exit distribution and occupation time of another CTRW.  In contrast to the the quadratic case, the CTRW is rather complicated for non-quadratic $\CV$ as its jump rates are not only \emph{random}, but additionally are \emph{time varying} and \emph{depend} on the boundary data.  Nevertheless, the HS representation is a rather useful analytical tool due to comparison inequalities (Brascamp-Lieb and Nash-Aronson).

Specifically, let $h_t^\psi$ solve \eqref{gl::eqn::dynam} with boundary condition $\psi$.  Conditional on the realization of the trajectory of the time-varying gradient field $(\nabla h_t^\psi(b) : b \in D^*)$, we let $X_t^\psi$ be the Markov process on $G$ with time-varying jump rates $\CV''(\nabla h_t^\psi(b))$.  Let $\tau = \inf\{ t \geq 0 : X_t^\psi \in \partial V\}$.  Let $\p_x$ denote the joint law of $(h_t^\psi, X_t^\psi)$ given $X_0^\psi = x$ and $\E_x$ the expectation under $\p_x$.

\begin{lemma}
\label{gl::lem::hs_mean_cov}
The mean and covariances of $h^\psi$ admit the representation
\begin{align}
   \cov(h^\psi(x),h^\psi(y)) &= \E_x \int_0^\tau \one_{\{X_s^\psi = y\}} ds\\
   \E h^\psi(x) &= \int_0^1 \E_x \psi(X_\tau^{r\psi}) dr. 
\end{align}
\end{lemma}

We refer the reader to \cite[Section 2]{DGI00} for a proof and also a much more detailed discussion on the HS representation.

\subsection{Brascamp-Lieb Inequalities}
\label{subsec::bl}
For $\nu \in \R^{|V|}$, we let
\[ \langle \nu, h^\psi \rangle = \sum_{x \in V} \nu(x) h^\psi(x).\]

The following inequalities, first proved in \cite{BL76} and redeveloped probabilistically in \cite{DGI00}, bound from above the centered moments of $h^\psi$ with those of $h^*$, where $h^*$ is a zero-boundary DGFF on $G$.  Recall that $a_\CV,A_\CV$ are positive, finite constants so that $a_\CV \leq \CV'' \leq A_\CV$.

\begin{lemma}[Brascamp-Lieb inequalities]
\label{bl::lem::bl_inequalities}
 There exists a constant $C > 0$ depending only on $a_\CV,A_\CV$ such that the following inequalities hold:
\begin{align}
 & \var(\langle \nu, h^\psi \rangle) \leq C  \var( \langle \nu, h^* \rangle ) \label{gl::eqn::bl_var}, \\
 & \E \exp(\langle \nu, h^\psi \rangle - \E \langle \nu,h^\psi \rangle) \leq \E \exp( C\langle \nu, h^* \rangle) \label{gl::eqn::bl_exp}
\end{align}
for all $\nu \in \R^{|V|}$.
\end{lemma}
We again refer the reader to \cite[Section 2]{DGI00} for a proof.  The Brascamp-Lieb inequalities allow for the following bound on the moments of the maximum which we will make use of many times throughout the rest of the article.

\begin{lemma}[Moments of the Maximum]
\label{gl::lem::mom_max}
Suppose that $F \subseteq \Z^2$ is bounded and connected with $R = \diam(F)$.  Let $\zeta \in \B_{\ol{\Lambda}}^u(F)$, $h^\zeta \sim \p_F^\zeta$, and $M = \max_{x \in F} |h^\zeta(x) - u \cdot x|$.  For every $\epsilon > 0$ and $p \geq 1$ we have that
\[ (\E M^p)^{1/p} = O_{\ol{\Lambda},p}(R^\epsilon).\] 
\end{lemma}
\begin{proof}
We may assume without loss of generality that $\ol{\Lambda} \geq 1$.  Combining the Brascamp-Lieb and Chebyshev inequalities, we have the tail bound
\[ \p[M \geq t]  = O_{\ol{\Lambda}}(\exp( O_{\ol{\Lambda}}((\log R)^{\ol{\Lambda}}) - t)).\]
Furthermore, we have 
\[ \E M^p \leq \sum_{x \in F} \E |h^\zeta(x)|^p \leq O_{\ol{\Lambda}}( R^2 (\log R)^{p\ol{\Lambda}}).\]
Consequently,
\begin{align*}
 \E M^p &\leq R^{p\epsilon} + \E M^p \one_{\{ M \geq R^{\epsilon}\}}
               \leq R^{p\epsilon} + (\E M^{2p})^{1/2} (\p[ M \geq R^{\epsilon}])^{1/2} \\
              &\leq R^{p\epsilon} + O_{\ol{\Lambda}}( R (\log R)^{p\ol{\Lambda}} \exp(c_p- \tfrac{1}{2}R^{\epsilon})).
 \end{align*}
 Therefore $(\E M^p)^{1/p} = O_{\ol{\Lambda},p}(R^\epsilon)$,
as desired.
\end{proof}

\section{Dynamics}
\label{sec::dyn}
We now specialize to the case where $G$ is a bounded, connected subgraph of $\Z^2$.  We will write $D$ for its vertices, $\partial D = \{ x \in \Z^2 : \dist(x,D) = 1\}$ for its boundary, and $D^* = \{ b = (x_b,y_b) \in (\Z^2)^* : x_b,y_b \in D\}$ for its edges, where $(\Z^2)^*$ denotes the set of edges of $\Z^2$.  Finally, let $\partial D^*$  be the set of edges that are either contained in $\partial D$ or intersect both $\partial D$ and $D$.  The Langevin dynamics are extremely useful for constructing couplings of instances of the GL model with either different boundary conditions, defined on different (though overlapping) domains, or both.  Suppose that $h^\psi, h^{\wt{\psi}}$ are solutions of \eqref{gl::eqn::dynam} driven by the same Brownian motions with boundary conditions $\psi,\wt{\psi}$, respectively.  Let $\ol{\psi} = \psi - \wt{\psi}$ and $\ol{h} = h^\psi - h^{\wt{\psi}}$.  Observe that
\begin{equation}
\label{gl::eqn::diff_dynam}
d \ol{h}_t(x) = \sum_{b \ni x} [\CV'(\nabla(h_t^\psi \vee \psi)(b)) - \CV'(\nabla (h_t^{\wt{\psi}} \vee \wt{\psi})(b))] dt
\end{equation}
Let
\begin{equation}
 \label{gl::eqn::stat_rate_def} c_t(b) = \int_0^1 \CV''(\nabla (h_t^{\wt{\psi}} + s \ol{h}_t)(b)) ds \text{ and } 
 \CL_t f(x) = \sum_{b \ni x} c_t(b) \nabla f(b).
 \end{equation}
Then we can rewrite \eqref{gl::eqn::diff_dynam} more concisely as
\begin{equation}
\label{gl::eqn::diff_dynam_ell}
d \ol{h}_t(x) = \CL_t \ol{h}_t(x) dt.
\end{equation}
The following is \cite[Lemma 2.3]{FS97}:
\begin{lemma}[Energy Inequality]
\label{gl::lem::ee}
For every $T > 0$ we have
\begin{align}
\label{gl::eqn::ee}
&\sum_{x \in D} |\ol{h}_T(x)|^2 + \int_0^T \sum_{b \in D^*} |\nabla \ol{h}_t(b)|^2 dt \notag\\
\leq& C \left(\sum_{x \in D} |\ol{h}_0(x)|^2 + \int_0^T \sum_{b \in \partial D^*} |\ol{\psi}(x_b)||\nabla \ol{h}_t(b)|dt\right)
\end{align}
for $C > 0$ depending only on $a_\CV,A_\CV$.
\end{lemma}
More generally, if $f_t$ solves $\partial_t f = \CL_t f_t$ then $f_t$ also satisfies \eqref{gl::eqn::ee}.

\subsection{Coupling Bounds}
\label{dyn::subsec::coupling_bounds}

The purpose of the next lemma is to show that $\lim_{T \to \infty} (h_{T}^\psi, h_{T}^{\wt{\psi}})$ gives the unique invariant measure of the Markov process $(h_t^\psi,h_t^{\wt{\psi}})$, i.e. where $h_t^\psi,h_t^{\wt{\psi}}$ both solve \eqref{gl::eqn::dynam} with the same driving Brownian motions.

\begin{lemma}$\ $
\label{gl::lem::ergodic}
\begin{enumerate}
 \item \label{gl::lem::unique_invariant} The SDS \eqref{gl::eqn::dynam} is ergodic.
 \item \label{gl::lem::unique_invariant_multiple} More generally, any finite collection $h^1,\ldots,h^n$ satisfying the SDS \eqref{gl::eqn::dynam} and driven by the same family of Brownian motions is ergodic.
\end{enumerate}
\end{lemma}
\begin{proof}
Part \eqref{gl::lem::unique_invariant} follows immediately from Lemma \ref{gl::lem::ee}.  Indeed, the Poincar\'e inequality implies that there exists $c_D > 0$ such that for all functions $f \colon D \to \R$ with $f|_{\partial D} \equiv 0$ we have
\[ \sum_{x \in D} (f(x))^2 \leq c_D \sum_{b \in D^*} (\nabla f(b))^2.\]
Thus Lemma \ref{gl::lem::ee} implies that if $(h_t^\psi,\wt{h}_t^\psi)$ both solve \eqref{gl::eqn::dynam} with the same Brownian motions and boundary data (though with possibility different initial distributions) then with $\ol{h}_t = h_t^\psi - \wt{h}_t^\psi$ we have
\[ \sum_{x \in D} (\ol{h}_T(x))^2 + \frac{1}{c_D} \int_0^T \sum_{x \in D} (\ol{h}_t(x))^2 dt \leq C \sum_{x \in D^*} (\ol{h}_0(b))^2.\]
In particular, the integral is bounded as $T \to \infty$ which implies
\[ \lim_{T \to \infty} \int_T^\infty \sum_{x \in D^*} (\ol{h}_t(x))^2 dt = 0.\]
The energy inequality also gives
\[ \frac{1}{C}\sum_{x \in D} (\ol{h}_T(x))^2 \leq \sum_{x \in D} (\ol{h}_t(x))^2\]
for all $0 < t < T$.  Therefore $\ol{h}_t \stackrel{d}{\to} 0$ as $t \to \infty$.  This proves \eqref{gl::lem::unique_invariant}.

We will now prove part \eqref{gl::lem::unique_invariant_multiple}.  In the interest of keeping the notation simple, we will prove the result in the special case $n=2$.  Suppose that $(h_t^1,h_t^2)$ solve the SDS \eqref{gl::eqn::dynam} driven by the same Brownian motions with boundary conditions $\psi^1,\psi^2$ but with arbitrary initial conditions $(h_0^1, h_0^2)$.  We know that $h_t^i$ converges in distribution $\p_D^{\psi^i}$ by part \eqref{gl::lem::unique_invariant}.  Consequently, the pair $(h_t^1,h_t^2)$ is tight.

We will now prove the existence of the limit $\lim_{t \to \infty} (h_t^1,h_t^2)$.  Suppose that $(T_k)$ and $(S_k)$ are arbitrary increasing sequences diverging to infinity.  Fix $k$ and assume that $T_k \leq S_k$.  Let $(\wt{h}_t^1, \wt{h}_t^2)$ solve \eqref{gl::eqn::dynam} where $\wt{h}_0^i = h_{S_k-T_k}^i$.  Let $(\breve{h}_t^1,\breve{h}_t^2)$ be another pair of solutions to \eqref{gl::eqn::dynam} with the same boundary and initial conditions of $(h^1,h^2)$ but driven by the same Brownian motions as $(\wt{h}^1,\wt{h}^2)$.  Then $(\wt{h}_{T_k}^1,\wt{h}_{T_k}^2) \stackrel{d}{=} (h_{S_k}^1, h_{S_k}^2)$ and $(\breve{h}_{T_k}^1, \breve{h}_{T_k}^2) \stackrel{d}{=} (h_{T_k}^1, h_{T_k}^2)$.  By the energy inequality,
\begin{align*}
 \frac{1}{T_k} \int_{0}^{T_k} \sum_{b \in D^*} |\nabla (\breve{h}_t^i - \wt{h}_t^i)(b)|^2 dt
 &\leq \frac{C}{T_k} \sum_{x \in D} |\breve{h}_0^i(x) - \wt{h}_0^i(x)|^2\\
 &\leq \frac{2C}{T_k} \sum_{x \in D} \big(|h_0^i(x)|^2 + |h_{S_k-T_k}^i(x)|^2 \big).
 \end{align*}
As $T_k \to \infty$, the first term in the summation on the right hand side clearly converges to zero almost surely.  The second term in the summation converges to zero in distribution since $h_t^i$ is tight.  Consequently, for every $\delta > 0$ there exists $k_0$ sufficiently large so that for all $k \geq k_0$ with $S_k \geq T_k$ we have
\[ \p\left[   \frac{1}{T_k} \int_{0}^{T_k} \sum_{b \in D^*} |\nabla (\breve{h}_t^i - \wt{h}_t^i)(b)|^2 dt > \delta \right] < \delta.\]

Let $\epsilon > 0$ be arbitrary and fix $\delta = \epsilon / c_D$ so that if $f \colon D \to \R$ is an arbitrary function vanishing on $\partial D$ with $\sum_{b \in D^*} |\nabla f(b)|^2 \leq \delta$ then $\sum_{x \in D} |f(x)|^2 \leq \epsilon$.  Assume that $k$ is sufficiently large so that with probability $1-\epsilon$ we have 
\[\frac{1}{T_k} \int_{0}^{T_k} \sum_{b \in D^*} |\nabla (\breve{h}_t^i - \wt{h}_t^i)(b)|^2 dt \leq \delta.\]
Then there exists (random) $t_0 \in [0,T_k/2]$ such that $\sum_{b \in D^*} |\nabla (\breve{h}_{t_0}^i - \wt{h}_{t_0}^i)(b)|^2 \leq \delta$ hence $\sum_{x \in D} | \breve{h}_{t_0}^i(x) - \wt{h}_{t_0}^i(x)|^2 \leq \epsilon$ with probability $1-\epsilon$.  Applying the energy inequality once again yields
\begin{align*}
  \sum_{x \in D} | \breve{h}_{T_k}^i(x) - \wt{h}_{T_k}^i(x)|^2
 &\leq  C\sum_{x \in D} |\breve{h}_{t_0}^i(x) - \wt{h}_{t_0}^i(x)|^2
    \leq C\epsilon
 \end{align*}
 with probability $1-\epsilon$.  Of course, we can do the same if $T_k \geq S_k$.  Therefore we conclude that the subsequential limits of $(h_t^1,h_t^2)$ are unique, hence $\mu \stackrel{d}{=} \lim_{t \to \infty} (h_t^1, h_t^2)$ exists.  The same argument also implies that for any $s > 0$ we have $\lim_{t \to \infty} (h_{s+t}^1,h_{s+t}^2)$ exists and has the same distribution as $\mu$.  Therefore $\mu$ is a stationary measure.
 
 To finish proving the lemma, we just need to establish uniqueness.  Suppose that each of the pairs $(h_t^1, h_t^2), (\wt{h}_t^1, \wt{h}_t^2)$ solve the SDS \eqref{gl::eqn::dynam}, $h_t^i, \wt{h}_t^i$ all driven by the same Brownian motions.  Suppose further that both pairs have stationary initial conditions.  Then we can use the energy inequality exactly in the same manner as in the proof of part \eqref{gl::lem::unique_invariant} to deduce that $|h_t^i - \wt{h}_t^i| \stackrel{d}{\to} 0$ as $t \to \infty$.  Since $(h_t^1,h_t^2)$, $(\wt{h}_t^1, \wt{h}_t^2)$ are stationary, it therefore follows that $(h_0^1,h_0^2) \stackrel{d}{=} (\wt{h}_0^1, \wt{h}_0^2)$.
\end{proof}

Suppose that $(h_t^\psi,h_t^{\wt{\psi}})$ is the stationary coupling of instances of the model with boundary conditions $\psi,\wt{\psi} \in \B_{\ol{\Lambda}}^u(D)$, respectively.  Letting $\ol{h} = h^\psi - h^{\wt{\psi}}$, the Caccioppoli inequality \eqref{symm_rw::eqn::caccioppoli} implies that
\begin{equation}
\label{gl::lem::cap_bound}
 \int_{r^2}^{2r^2} \sum_{b \in B^*(x_0,r)} |\nabla \ol{h}_t(b)|^2 dt \leq \frac{C}{r^2} \int_0^{2r^2} \sum_{x \in B(x_0,r)} |\ol{h}_t(x)|^2 dt
 \end{equation}
for $r > 0$ and $x_0 \in D$ with $B(x_0,2r) \subseteq D$.  The maximum principle implies that $\ol{h}$ attains its maximum on $\partial D$, hence $\ol{h} = O_{\ol{\Lambda}}((\log R)^{\ol{\Lambda}})$.  Consequently, taking expectations of both sides of \eqref{gl::lem::cap_bound} and using the stationarity of the dynamics yields
\[ \sum_{b \in B^*(x_0,r)} \E |\nabla \ol{h}_t(b)|^2 \leq O_{\ol{\Lambda}}( (\log R)^{2 \ol{\Lambda}}).\]
Recall that $D(r) = \{ x \in D : \dist(x,\partial D) \geq r\}$ and let $\epsilon > 0$ be arbitrary.  Then $D(R^{1-\epsilon})$ can be covered by $O(R^2 / R^{2-2\epsilon}) = O(R^{2\epsilon})$ balls of radius $R^{1-\epsilon}$, all of which are contained in $D$.  Therefore
\[ \sum_{b \in D^*(R^{1-\epsilon})} \E |\nabla \ol{h}_t(b)|^2 \leq O_{\ol{\Lambda}}( R^{3\epsilon}).\]
Let $D(R_1,R_2) = \{x \in D : R_1 \leq \dist(x,\partial D) < R_2\}$, $r_k = kR^{1-5\epsilon}$, and $D_k = D(R^{1-\epsilon}+r_k, R^{1-\epsilon} + r_{k+1})$.  Note that we can write
\[ \sum_{b \in D^*(R^{1-\epsilon},2R^{1-\epsilon})} \E | \nabla \ol{h}_t(b)|^2 = \sum_{k=0}^{R^{4\epsilon}-1} \sum_{b \in D_k^*} \E | \nabla \ol{h}_t(b)|^2  = O_{\ol{\Lambda}}(R^{3\epsilon}).\]
This implies there exists $0 \leq k \leq R^{4\epsilon}-1$ such that 
\[ \sum_{b \in D_k^*} \E |\nabla \ol{h}_t(b)|^2 \leq O_{\ol{\Lambda}}(R^{-\epsilon})\]
We have proven:

\begin{lemma}
\label{gl::lem::grad_error}
Suppose that $(h^\psi,h^{\wt{\psi}})$ is a stationary coupling of two solutions of the SDS \eqref{gl::eqn::dynam} driven by the same Brownian motions with $\psi,\wt{\psi} \in \B_{\ol{\Lambda}}^u(D)$.  For every $\epsilon > 0$ there exists $R^{1-\epsilon} \leq R_1 \leq 2R^{1-\epsilon}$ such that with $R_2 = R_1 + R^{1-5\epsilon}$ we have that 
\begin{align*}
  \sum_{b \in D^*(R_1, R_2)} \E |\nabla \ol{h}_t(b)|^2 &= O_{\ol{\Lambda}}(R^{-\epsilon}) \text{ and } 
  \sum_{b \in D^*(R_1)} \E |\nabla \ol{h}_t(b)|^2 &= O_{\ol{\Lambda}}(R^{3\epsilon}).
\end{align*}
\end{lemma}

This lemma will be particularly useful for us in Section \ref{sec::harm} in order to construct an intermediate coupling of $\p_D^\psi, \p_D^{\wt{\psi}}$ exhibiting pointwise regularity near $\partial D$ with high probability.

\subsection{Gradient Gibbs States and Equivalence of Ensembles}
\label{subsec::shift_invariant}

By the reverse Brascamp-Lieb inequality \cite[Lemma 2.8]{DGI00}, it follows that if $D_n$ is any sequence of domains tending locally to the infinite lattice $\Z^2$ and, for each $n$, $h^n$ is an instance of the GL model on $D_n$ then $\var(h^n(x)) \to \infty$ as $n \to \infty$.  This holds regardless of the choice of boundary conditions, which suggests that it is not possible to take an infinite volume limit of the height field $h^n(x)$.  However, the Brascamp-Lieb inequality (Lemma \ref{bl::lem::bl_inequalities}) gives that $\var(\nabla h^n(b))$ remains uniformly bounded as $n \to \infty$, indicating that it should be possible to take an infinite volume limit of the \emph{gradient field}.

Working with gradient rather than height fields, though unnecessary for $d \geq 3$, is convenient since it allows for a unified treatment of Gibbs states for all dimensions.  Let $\CX$ be the set of functions $\eta \colon (\Z^d)^* \to \R$.  Let $\CF = \sigma(\eta(b) : b \in (\Z^d)^*)$ be the $\sigma$-algebra on $\CX$ generated by the evaluation maps and, for each $D^* \subseteq (\Z^d)^*$, let $\CF_{D^*} = \sigma(\eta(b) : b \in D^*)$ be the $\sigma$-algebra generated by the evaluation maps in $D^*$.  Suppose that $D \subseteq \Z^d$ is bounded, $\varphi \in \CX$, and $\nabla \phi = \varphi$.  If $h$ is distributed according to $\p_D^\phi$, then the gradient field $\nabla h$ induces a measure $\p_{D^*}^\varphi$ on functions $D^* \to \R$.  We call $\p_{D^*}^\varphi$ the law of the GL model on $D^*$ with \emph{Neumann boundary} conditions $\varphi$.  Let $\mu$ be a measure on $\CX$ and suppose that $\eta$ has the law $\mu$.  We say that $\mu$ is a \emph{gradient Gibbs state} associated with the potential $\CV$ if for every finite $D^* \subseteq (\Z^2)^*$, 
\[ \mu(\cdot | \CF_{(D^*)^c}) =\p_{D^*}^{\eta|\partial D^*} \text{ almost surely.}\]

Fix a vector $u \in \R^d$.  A gradient Gibbs state $\mu$ is said to have \emph{tilt} $u$ if $\E^\mu \eta(x+ b_i) = u \cdot e_i$ where $b_i = (0,e_i)$, $e_i$ a generator of $\Z^d$, and $x \in \Z^d$ is arbitrary.  We say that $\mu$ is \emph{shift invariant} if $\mu \circ \tau_x^{-1} = \mu$ for every $x \in \Z^d$ where $\tau_x \colon \Z^d \to \Z^d$ is translation by $x$.  Finally, a shift-invariant $\mu$ is said to be shift-ergodic if whenever $f$ is a shift-invariant $\CF$-measurable function, then $f$ is $\mu$-almost surely constant.

Funaki and Spohn in \cite{FS97} proved that the shift-ergodic Gibbs states are parameterized according to their tilt $u$; from now on we will refer to the law of such as $\SEGGS_u$.  The natural construction is to take an infinite volume limit of gradient measures $\p_{D^*}^\varphi$ as $D^*$ tends locally to $\Z^d$ with boundary conditions $\varphi(b) = u \cdot (y_b - x_b)$ \cite[Remark 4.3]{STF_Funaki}.  The difficulties with this approach are that $\p_{D^*}^\varphi$ is itself not shift-invariant and it is not clear that the mean gradient field approximately has tilt $u$.  This issue is cleverly circumvented in \cite{FS97} by instead considering the finite volume measure
\[ d\mu_n(\eta) = \frac{1}{\CZ_n} \exp\left( -\sum_{b \in (\Z_n^d)^*} \CV(\eta(b) - (y_b - x_b) \cdot u) \right) d\nu_n(\eta)\]
on gradient fields on the torus, where $\nu_n$ is the uniform measure on the set of functions $\eta \colon (\Z_n^d)^* \to \R$ which can be expressed as the gradient of a function $h \colon \Z_n^d \to \R$.  By construction, $\mu_n$ is shift invariant, has tilt $u$ and both of these properties are preserved in the limit as $n \to \infty$.

We will now explain how to use the method of dynamic coupling to prove that the gradient field of $\p_{D^*}^\varphi$, $\varphi(b) = u \cdot (y_b - x_b)$ as before, converges to $\SEGGS_u$.  Our proof will also yield an alternative construction of the Funaki-Spohn state in the special case $u = 0$.  We will include a statement of this result here as well as a short sketch of the proof since this will be important for the proof of Theorem \ref{thm::clt}.  We note in passing that Theorem \ref{harm::thm::coupling} gives a much better coupling which will be critical for the proof of Theorem \ref{thm::clt}, but this result uses the convergence of the finite volume gradient fields hence we cannot simply apply it here.

Fix a tilt $u$.  Suppose that $D_n$ is any sequence of bounded domains in $\Z^d$ converging locally to $\Z^d$.  For each $n$, let $\eta^n \sim \p_{D_n^*}^{\varphi}$.  Suppose that $\eta \sim \SEGGS_u$.  Fix $x_n \in \partial D_n$ and let $h^n, h^{n,S}$ be the height fields associated with the gradient fields $\eta^n, \eta$, respectively, both set to vanish at $x_n$.  
By the Brascamp-Lieb inequality \eqref{gl::eqn::bl_var},
\begin{equation}
\label{eqn::var_bound}
\var( h^{n,S}(x) - h^{n,S}(y)) \leq  C \log(1+|x-y|).
\end{equation}
Here, $h^{n,S}(x) - h^{n,S}(y)$ is interpreted as $\sum_{i=1}^n \eta(b_i)$ where $b_1,\ldots,b_n$ is any sequence of bonds connecting $x$ to $y$.  Of course the same is also true with $h^n$ in place of $h^{n,S}$.  Let $R_n = \diam(D_n)$.  As $h^n(x) = \E h^{n,S}(x)$ for $x \in \partial D_n$, the Brascamp-Lieb \eqref{gl::eqn::bl_exp} and Chebychev inequalities thus imply that
\begin{equation}
\label{eqn::gs_diff}
\p[ \max_{x \in \partial D_n} |h^{n,S}(x) - h^n(x)| \geq (\log R_n)^2] = O(R_n^{-100}).
\end{equation}
Assume that $(h_t^{n,S}, h_t^n)$ is the stationary coupling of $h^n$ and $h^{n,S}$ conditional on $h^{n,S} |_{\partial D_n}$.  Then as $\ol{h}_t^n = h_t^{n,S} - h_t^n$ satisfies the parabolic equation \eqref{gl::eqn::diff_dynam_ell} and $\ol{h}^n$ is static on $\partial D_n$, the maximum principle implies
\[ \max_{x \in D_n} |\ol{h}_t^n(x)| \leq \max_{x \in \partial D_n} |\ol{h}_0^n(x)|.\]
Combining this with \eqref{eqn::gs_diff} implies
\[ \p[ \max_{x \in D_n} |\ol{h}_t^n(x)|\geq (\log R_n)^2] = O(R_n^{-100}).\]
The Nash continuity estimate (Lemma \ref{symm_rw::lem::nash_continuity_bounded}) applied to $\ol{h}_t^n$ thus implies that
\begin{equation}
\label{gl::eqn::quantitative_bound}
 \E \big[ \max_{b \in D_n^*(R_n^\zeta)} |\nabla \ol{h}_t^n(b)| \big] = O_{\ol{\Lambda}}( R_n^{\epsilon-\zeta \xi_{\rm NC}})
 \end{equation}
for $\epsilon, \zeta > 0$ fixed.  This proves the desired convergence.

Existence in the special case of $u = 0$ can be proved in a very similar manner.  The reason is that, in this case, $h^n \sim \p_{D_n}^0$ hence $\E h^n(x) = 0$.  Thus it is clear that the subsequential limits of gradient fields, which exist by the Brascamp-Lieb inequalities, have zero tilt.  The Brascamp-Lieb and Chebychev inequalities imply that the maximum of $h^n$ is with high probability $O(\log R_n)$.  Thus if $\wt{h}^n$ denotes the law of the GL model on the domain $\wt{D}_n$ given by shifting $D_n$ by one unit, then using the argument of the previous paragraph we can couple $h^n$ and $\wt{h}^n$ such that $\nabla (h^n - \wt{h}^n)$ is with high probability $O(R_n^{\epsilon - \zeta \xi_{\rm NC}})$ at distance at least $R_n^{\zeta}$, $\zeta > 0$, from both $\partial D_n$ and $\partial \wt{D}_n$.  Therefore the subsequential limits of $\p_{D_n}^0$ are shift-invariant which proves the existence of a zero-tilt shift-invariant Gibbs state.  Uniqueness (and also existence of limits) follows by taking two such states $\eta, \eta'$, then applying the argument of the previous paragraph.

We have obtained:
\begin{theorem}[Equivalence of Ensembles] $\ $
\label{gl::thm::infinite_coupling}
\item If $(D_n)$ is any sequence of bounded domains in $\Z^d$ tending locally to $\Z^d$ and, for each $n$, $h^n$ is an instance of the GL model on $D_n$ with boundary conditions $\varphi_n \in \B_{\ol{\Lambda}}^u(D)$, then $\eta^n = \nabla h^n$ converges weakly to $\SEGGS_u$ and we have
\[ \E \big[ \max_{b \in D_n^*(R_n^\zeta)} |\eta^n(b) - \eta(b)| \big] = O_{\ol{\Lambda}}(R^{\epsilon -\zeta \xi_{\rm NC}}).\]
\end{theorem}

\subsection{Invariance under Reflections}
\label{subsec::reflect_invariance}

The following proposition will be especially useful for us later in the next section when applied to $f = \CV''$.  We let $\varphi^h,\varphi^v \colon \R^2 \to \R^2$ be the maps which reflect about the horizontal and vertical axes, respectively, and $\varphi_x^h(b) = \varphi^h(b-x)$, $\varphi_x^v(b) = \varphi^v(b-x)$ for $x \in \Z^2$.

\begin{proposition}
\label{gl::prop::reflection_invariance}
Fix a tilt $u$, $x \in \Z^2$, and let $f \colon \R \to \R$ be an even function.  Let $b^h = \varphi_x^h(b)$ for $b \in (\Z^2)^*$.
If $\eta_u \sim \SEGGS_u$, then $( f(\eta_u(b)) : b \in (\Z^2)^*) \stackrel{d}{=} ( f(\eta_{u}(b^h)) : b \in (\Z^2)^*)$.  The same is also true when the horizontal reflection is replaced with vertical reflection.
\end{proposition}
\begin{proof}
Let $w = \varphi^h(u)$ and $\eta_w \sim \SEGGS_w$.  Let $s_b = 1$ if $b$ is orientated vertically and $-1$ otherwise.
Since $f$ is even, we know that $f(\eta_u(b)) = f( s_b \eta_u(b))$.  Since $(s_b \eta_u(b)) \stackrel{d}{=} (\eta_w(b))$, we thus have that 
\[ (f(\eta_u(b))) \stackrel{d}{=} (f(\eta_w(b))) \stackrel{d}{=} f(\eta_u(b^h)).\]
The reason for the last inequality is that $( \eta_u(b^h) : b \in (\Z^2)^*)$ is still a shift-ergodic Gibbs state but with tilt $w$.
\end{proof}

\section{Correlation Decay}
\label{sec::correlation_decay}

Suppose that $F \subseteq \Z^2$ is a bounded, connected domain with $R = \diam(F)$.  Let $\zeta, \wt{\zeta} \in \B_{\ol{\Lambda}}^u(F)$ and assume that $(h_t^\zeta,h_t^{\wt{\zeta}})$ is the stationary coupling of $\p_F^\zeta, \p_F^{\wt{\zeta}}$.  Throughout, we let $\beta = \beta(u)$ as in the statement of Theorem \ref{thm::clt}.  The main result of this section is that $\CV''(\nabla h^\zeta(b))$ and $\nabla \ol{h}(b)$, $\ol{h} = h^\zeta - h^{\wt{\zeta}}$, are uncorrelated when averaged against a $\Delta^\beta$-harmonic function. 

\begin{theorem}
\label{cd::thm::cd}
Suppose that $x_0 \in F$ with $\dist(x_0,\partial F) \geq R^{\alpha+\epsilon}$ for $\alpha,\epsilon > 0$, $E = B(x_0, R^\alpha)$, and let $g \colon F \to \R$ be $\Delta^\beta$-harmonic.  We have that
\begin{align}
\label{cd::eqn::estimate}
 \E \sum_{b \in E^*} \CV''(\nabla h^\zeta(b)) \nabla \ol{h}(b) \nabla g(b) =& \sum_{b \in E^*} a_u(b) \E[\nabla \ol{h}(b)] \nabla g(b) +\\
  &O_{\ol{\Lambda}}(R^{\epsilon+\alpha(1-\rho_{\rm CD}) } \| \nabla g \|_\infty) \notag
 \end{align}
for $\rho = \rho_{\rm CD}(\CV) > 0$ and $a_u(b) = \E[\CV''(\eta(b))]$ with $\eta \sim \SEGGS_u$.
\end{theorem}

Note that $a_u$ depends on $b$ only through its orientation (either vertical or horizontal) by the shift-invariance of $\eta$.  Our typical choice of $g$ will have $\| \nabla g \|_\infty = O(R^{-\alpha})$, in which case the exponent in the error term is actually negative.  The idea of the proof is to show that replacing $\ol{h}(x)$ by its average $\ol{h}^\rho(x)$ over the ball $B(x,R^{\rho})$ introduces a small amount of error.  The advantage of this replacement is that the time-derivative of $\ol{h}_t^\rho$ possesses more regularity than that of $\ol{h}_t$.  This allows us to replace the left hand side of \eqref{cd::eqn::estimate} with
\[ \E \sum_{b \in E^*} \CV''(\nabla h_T^\zeta(b)) \nabla \ol{h}_0(b) \nabla g(b).\]
The proof is then completed by coupling $h_T^\zeta$ to $\eta \sim \SEGGS_u$ conditional on $\ol{h}_0$, which can be accomplished at the cost of negligible error by the argument used to prove Theorem \ref{gl::thm::infinite_coupling}.

\subsection{Change of Environment}

Let $\eta_t$ follow the $\SEGGS_u$ dynamics independent from $(h_t^\zeta,h_t^{\wt{\zeta}})$.  That is, $\eta_t$ solves the infinite dimensional SDS
\[ d\eta_t((x,y)) = \left( \sum_{b \ni y} \CV'(\eta_t(b)) - \sum_{b \ni x} \CV'(\eta_t(b)) \right) dt + \sqrt{2}(dW_t(y) - dW_t(x))\]
for $b \in (\Z^2)^*$ where $W_t(x), x \in \Z^2,$ is a family of iid standard Brownian motions; see Section 9 of \cite{STF_Funaki} for a discussion of the existence and uniqueness of solutions to this equation.
Fix $T > 0$, let $\breve{c}_t(b) = \CV''( \eta_t(b)),$
and let $\breve{p}$ be the transition kernel of the random walk jumping with rates $\breve{c}_{T-t}(b)$, $t \in [0,T]$, stopped on its first exit from $F$.

\begin{proposition}
\label{cd::prop::environment_change}
Suppose that we have the setup as Theorem \ref{cd::thm::cd}.  Let $\gamma_1, \gamma_2 \in (0,\alpha]$ and $\delta_i = 4 \alpha - 4\gamma_i - \gamma_i \rho_{\rm EC}$.  Let $S_2 = R^{2\gamma_2}$.  There exists $\tfrac{3}{4} R^{2\gamma_1} \leq S_1 \leq R^{2\gamma_1}$ such that the following holds.  Let $\breve{p}^1, \breve{p}^2 \sim \breve{p}$ associated with independent environments $\eta^1,\eta^2 \sim \SEGGS_u$ which are in turn independent of $\ol{h}_0$.  Let $S = S_1 + S_2$ and let
\[ \breve{h}_t(x) = \sum_{y,z} \breve{p}^1(S-t,S_1;x,y) \breve{p}^2(0,S_2;y,z) \ol{h}_0(z),\ \ S - S_2 \leq t \leq S.\]
There exists a coupling of $(\eta^1,\eta^2,\breve{h}_t)$ and $(\nabla h, \ol{h})$ such that
\begin{align*}
     \E \sum_{b \in E^*} (\nabla \breve{h}_{S}(b) - \nabla \ol{h}_{S}(b))^2 
=  O_{\ol{\Lambda}}(R^{\epsilon+2\alpha +2\gamma_2- 4 \gamma_1 + \delta_2}) + O_{\ol{\Lambda}}(R^{\epsilon + \delta_1}),\\
 \E \bigg[ \max_{b \in E^*} \sup_{S/2 \leq t \leq S} | \nabla h_{t}^\zeta(b) - \eta_{t}^1(b)|  \bigg] = O_{\ol{\Lambda}}(R^{\epsilon-\gamma \xi_{\rm NC}}),\label{harm::eqn::ec_change_grad}
\end{align*}
for $\rho_{\rm EC} > 0$ depending only on $\CV$.
\end{proposition}

The following heat kernel estimates are crucial ingredients in the proof of the proposition.

\begin{lemma}
\label{harm::lem::hc_time_space_deriv}
Suppose that $x_0 \in F$ with $\dist(x_0,\partial F) \geq R^{\gamma+\epsilon}$ for $\gamma,\epsilon > 0$ and let $E = B(x_0,R^{\gamma})$.  For $T = R^{2\gamma}$ and $\tfrac{3}{4} T \leq t_1 < t_2 \leq T$, we have
\begin{align*}
    \sum_{x \in E} \sum_{y \in F} | q(u,t_1;x,y) -  q(u,t_2;x,y)|^2  &= O (|t_1-t_2|^2 R^{\epsilon - 2 \gamma \xi_{\rm NC}})
\end{align*}
for $0 \leq u \leq T/4$ and
\begin{align*}
    \sum_{b \in E^*} \sum_{y \in F} \int_0^{T/4} |\nabla q(u,t_1;b,y) - \nabla q(u,t_2;b,y)|^2 du  &= O (|t_1-t_2|^2 R^{\epsilon-2\gamma \xi_{\rm NC}})
\end{align*}
where $q$ is the transition kernel of a random walk on $F$ evolving with rates $a_\CV \leq d_t(b) \leq A_\CV$.
\end{lemma}
\begin{proof}
Using that $\partial_t q(u,t;x,y) = \sum_{b \ni y} d_t(b) \nabla q(u,t;x,b)$, we have 
\begin{align*}
      &  \sum_{b \in E^*} \sum_{y \in F}  \int_{0}^{T/4} |\nabla q(u,t_1;b,y) - \nabla q(u,t_2;b,y)|^2 du \\
\leq& |t_1 - t_2| \sum_{b \in E^*} \sum_{y \in F} \int_{t_1}^{t_2} \int_0^{T/4} |\partial_t \nabla q(u,t;b,y)|^2 du dt\\
\leq& C |t_1 - t_2| \sum_{b \in E^*}\sum_{b' \in F^*}  \int_{t_1}^{t_2} \int_0^{T/4} |\nabla \nabla q(u,t;b,b')|^2 du dt,
\end{align*}
for $C > 0$ depending only on $\CV$.  Applying the Cacciopoli inequality \eqref{symm_rw::eqn::caccioppoli} to the first time and spatial coordinates, we see that this is bounded from above by
\begin{equation}
\label{harm::eqn::hc_mixed_deriv_est}
 O(1) \frac{|t_1- t_2|}{R^{2\gamma}} \sum_{x \in E} \sum_{b' \in F^*} \int_{t_1}^{t_2} \int_{0}^{T/2}  |\nabla  q(u,t;x,b')|^2 du dt.
 \end{equation}
The Nash-Aronson estimates (Lemma \ref{symm_rw::lem::nash_aronson}) imply that the contribution to the sum given by those $b' \in F^*$ with $\dist(b', E) \geq R^{\gamma+\epsilon/2}$ is negligible in comparison to the upper bound we seek to prove.  For $b' \in F^*$ with $\dist(b',E) \leq R^{\gamma+\epsilon/2}$, the Nash continuity and Nash-Aronson estimates (Lemmas \ref{symm_rw::lem::nash_continuity_bounded}, \ref{symm_rw::lem::nash_aronson}) imply
\[ |\nabla  q(u,t;x,b')| = O(R^{-\gamma \xi_{\rm NC}-2\gamma})\]
for $0 \leq  u \leq T/2$ and $t_1 \leq t \leq t_2$.
Inserting this into \eqref{harm::eqn::hc_mixed_deriv_est}, we arrive at the bound
\begin{align*}
    |t_1-t_2|^2 O(R^{4\gamma+\epsilon}) \cdot O(R^{- 2\gamma \xi_{\rm NC}-4\gamma})
= O(|t_1-t_2|^2 R^{\epsilon-2\gamma \xi_{\rm NC}}).
\end{align*}
This proves the second part of the lemma.  The first is exactly the same except the application of the Cacciopoli inequality is unnecessary.
\end{proof}

\begin{lemma}
Suppose that $x_0 \in F$ with $\dist(x_0,\partial F) \geq R^{\gamma+\epsilon}$ for $\gamma,\epsilon > 0$ and let $E = B(x_0,R^\gamma)$, $E_0 = B(x_0,R^{\gamma+\epsilon/2})$, $E_1 = B(x_0,R^{\gamma+\epsilon})$.  Let $q,q'$ be the transition kernels for two random walks on $F$ jumping with rates $a_\CV \leq d_t(b), d_t'(b) \leq A_\CV$, respectively.  With $T = R^{2\gamma}$, assume that $d_t \equiv d_t'$ for all $0 \leq t \leq \tfrac{1}{2} T$ and let
\[ \ol{d}_\infty = \sup_{0 \leq t \leq T} \max_{b \in E_1^*} |d_t(b) - d_t'(b)|.\]
Uniformly in $x \in E_0$ we have that
\[ \sum_{y \in E_1^*} | \ol{q}(0,T;x,y)|^2 +  \sum_{b \in E_1^*} \int_0^{T} |\nabla \ol{q}(0,t;x,b)|^2 dt = 
   O(\ol{d}_\infty^2 R^{-2\gamma})\]
where $\ol{q} = q-q'$.
\end{lemma}
\begin{proof}
By definition,
\[ \partial_t q(0,t;x,y) = [\CL_t q(0,t;x,\cdot)](y) \text{ and } \partial_t q'(0,t;x,y) = [\CL_t' q'(0,t;x,\cdot)](y)\]
where
\[ \CL_t f(y) = \sum_{b \ni y} d_{t}(b) \nabla f(b) \text{ and }  \CL_t' f(y) = \sum_{b \ni y} d_{t}'(b) \nabla f(b).\]
Consequently,
\begin{align}
 \partial_t \ol{q}^2
&= 2\ol{q}(\CL_t \ol{q} +\ol{\CL}_t q')
 \label{harm::eqn::diff_square}
\end{align}
where $\ol{\CL}_t = \CL_t - \CL_t'$.  Using the same proof as the energy inequality (Lemma \ref{gl::lem::ee}), by integrating both sides of \eqref{harm::eqn::diff_square} from $0$ to $T$ then using summation by parts, 
\begin{align*}
 &\sum_{y \in E_1} |\ol{q}(0,T;x,y)|^2 + 2\sum_{b \in E_1^*} \int_0^T d_t(b) |\nabla \ol{q}(0,t;x,b)|^2 dt\\
   \leq&  2\sum_{b \in \partial E_1^*} \int_0^T \bigg[ d_t(b) |\ol{q}(0,t;x,x_b) \nabla \ol{q}(0,t;x,b)| + |\ol{d}_t(b) \ol{q}(0,t;x,x_b)\nabla q'(0,t;x,b)|\bigg] dt+\\ 
        & 2\sum_{b \in E_1^*} \int_0^T |\ol{d}_{t}(b) \nabla \ol{q}(0,t;x,b) \nabla q'(0,t;x,b)| dt.
\end{align*}
By the Nash-Aronson estimates (Lemma \ref{symm_rw::lem::nash_aronson}),
\[ \sum_{b \in \partial E_1^*} \int_0^T |\ol{q}(0,t;x,x_b) \nabla \ol{q}(0,t;x,b)| dt = O(R^{3\gamma+\epsilon}\exp(-c'R^{\epsilon})) = O(\exp(-cR^{\epsilon}))\]
for some $c, c' > 0$ depending only on $a_\CV,A_\CV$.  The other boundary term is similarly of order $O(\exp(-c R^{\epsilon}))$.
Consequently,
\begin{align}
 &\sum_{y \in E_1} |\ol{q}(0,T;x,y)|^2 +   \sum_{b \in E_1^*} \int_{T/2}^T |\nabla \ol{q}(0,t;x,b)|^2 dt
   \leq O(\exp(-cR^{\epsilon})) + \notag\\
 & 
           C\ol{d}_\infty \sum_{b \in E_1^*} \int_{T/2}^T |\nabla \ol{q}(0,t;x,b) \nabla q'(0,t;x,b)| dt   \label{harm::eqn::ec_hc_upper_bound}.
\end{align}
The reason that the lower bound of integration is $T/2$ rather than $0$ is $d_t \equiv d_t'$ for all $t \leq T/2$. 
 If 
\[\sum_{b \in E_1^*} \int_{T/2}^T |\nabla \ol{q}(0,t;x,b)|^2 dt = O(\exp(-cR^{\epsilon/2}))\]
then we are obviously done.  If not, we apply Cauchy-Schwarz to the sum on the right hand side of \eqref{harm::eqn::ec_hc_upper_bound} and rearrange to arrive at
\begin{align*}
 & \sum_{b \in E_1^*} \int_{T/2}^T |\nabla \ol{q}(0,t;x,b)|^2 dt
   \leq C \ol{d}_\infty^2 \sum_{b \in E_1^*} \int_{T/2}^T |\nabla q'(0,t;x,b)|^2 dt,
\end{align*}
increasing $C$ if necessary.
By the Nash-Aronson estimates,
\begin{align*}
    \sum_{y \in E_1^*} |q'(0,T;x,y)|^2 &= O(R^{-2\gamma}) \text{ and}\\
    \sum_{b \in \partial E_1^*} \int_{T/2}^T |q'(0,t;x,x_b)| |\nabla q'(0,t;x,b)|dt &= O(\exp(-c R^{\epsilon})).
\end{align*}
Therefore, by the energy inequality (Lemma \ref{gl::lem::ee}),
\begin{equation}
\label{harm::eqn::hc_int_bound1}
\sum_{b \in E_1^*} \int_{T/2}^T |\nabla q'(0,t;x,b)|^2 dt = O(R^{-2\gamma}),
\end{equation}
which in turn implies
\begin{equation}
\label{harm::eqn::hc_int_bound2}
\sum_{b \in E_1^*} \int_{T/2}^T |\nabla \ol{q}(0,t;x,b)|^2 dt = O(\ol{d}_\infty^2 R^{-2\gamma}).
\end{equation}
Inserting \eqref{harm::eqn::hc_int_bound1}, \eqref{harm::eqn::hc_int_bound2} into \eqref{harm::eqn::ec_hc_upper_bound} proves the lemma.
\end{proof}

\begin{lemma}
Suppose that we have the same setup as Theorem \ref{cd::thm::cd} and let $\breve{p}, \eta$ be as in the introduction of this subsection.  Let $S = R^{2\gamma}$ for $0 \leq \gamma \leq \alpha$ and
\[ \breve{h}_t(y) = \sum_{z} \breve{p}(S-t,S;y,z) \ol{h}_0(z) \text{ for } 0 \leq t \leq S.\]
There exists a coupling of $(\eta,\breve{h})$ and $(\nabla h^\zeta, \ol{h})$ such that
\begin{align}
 \E \sum_{b \in E^*} \int_{3S/4}^{S} | \nabla \breve{h}_t(b) - \nabla \ol{h}_t(b)|^2 dt
   = O_{\ol{\Lambda}}(R^{\epsilon+4\alpha-2\gamma - \gamma \rho_{\rm EC}}), \label{harm::eqn::ec_change_int}\\
\E \sum_{x \in E} | \breve{h}_S(x) - \ol{h}_S(x)|^2 = O_{\ol{\Lambda}}(R^{\epsilon+4\alpha-2\gamma -\gamma \rho_{\rm EC}}), \label{harm::eqn::ec_change_sum}\\
 \E \bigg[ \max_{b \in E^*} \sup_{S/2 \leq t \leq S} | \nabla h_{t}^\zeta(b) - \eta_{t}(b)| \bigg] = O_{\ol{\Lambda}}(R^{\epsilon-\gamma \xi_{\rm NC}}),\label{harm::eqn::ec_change_grad_single}
\end{align}
where $\rho_{\rm EC} > 0$ depends only on $\CV$ and $\eta$ is independent of $(\nabla h_0^\zeta, \ol{h}_0)$.
\end{lemma}

This constant $\rho_{\rm EC}$ from Proposition \ref{cd::prop::environment_change} is the same as that appearing in the statement of this lemma.

\begin{proof}
Let $T = R^{2\alpha}$ and let $\eta_0 \sim \SEGGS_u$ independent of $(h_0^\zeta,h_0^{\wt{\zeta}})$.  Assume further that the evolution of the Brownian motions driving $\eta_t$ in $F$ are independent from those of $(h_t^\zeta, h_t^{\wt{\zeta}})$ until time $T- S$, $S \equiv R^{2\gamma}$, after which they are the same.  Let $\breve{c}_t(b) = \CV''(\eta_{t}(b))$, $c_t(b)$ be as in \eqref{gl::eqn::stat_rate_def} with $h^\zeta, h^{\wt{\zeta}}$ in place of $h^\psi, h^{\wt{\psi}}$, and let $p$ be the transition kernel of a random walk in $F$ stopped on its first exit jumping with rates $c_{T-t}(b)$.  Note that
\[ \ol{h}_t(x) = \sum_{y \in F} p(T-t,T;x,y) \ol{h}_0(y) \text{ for } 0 \leq t \leq T.\]
Set $S' = S - R^{2\sigma}$, $\sigma \in (0,\gamma)$ to be determined later, and define environments 
\[ \wt{c}_t(b) = \begin{cases}
			\breve{c}_t(b) \text{ for } T-S \leq t \leq T\\
			c_t(b) \text{ for } 0 \leq t < T-S
                    \end{cases}, \ \ \ 
  \wt{c}_t'(b) = \begin{cases}
			\breve{c}_t(b) \text{ for } T-S' \leq t \leq T\\
			c_t(b) \text{ for } 0 \leq t < T-S'.
                    \end{cases}\]
Let $\wt{p}, \wt{p}'$ be the transition kernels for the random walks in $F$ stopped on their first exit jumping with rates $\wt{c}_{T-t}(b), \wt{c}_{T-t}'(b)$, respectively.  Finally, let $\ol{p} = p - \wt{p}$ and $\ol{p}' = p - \wt{p}'$.
For $S \leq t \leq T$, we have $p(S,t;x,y) = \wt{p}(S,t;x,y) = \wt{p}'(S,t;x,y)$, hence
\begin{align*}
 &  \sum_{b \in E^*} \sum_{y \in F} (\nabla \ol{p}(u,T;b,y))^2
= \sum_{b \in E^*} \sum_{y \in F} \left( \sum_{z \in F} \nabla \ol{p}(u,S;b,z) p(S, T;z,y)\right)^2.
\end{align*}
By Jensen's inequality, this is bounded from above by
\begin{align*}
  & 4\sum_{b \in E^*} \sum_{z \in F} \bigg( (\nabla \ol{p}(u,S;b,z) - \nabla \ol{p}(u,S';b,z))^2\\
    &+ 
     (\nabla \ol{p}'(u,S';b,z) - \nabla \ol{p}'(u,S;b,z))^2 
     +(\nabla \ol{p}'(u,T;b,z))^2 \bigg) \equiv I_1 + I_2 + I_3.
\end{align*}
Fix a base point $a_0 \in \partial F$, set $h_t^S(a_0) = 0$, and let $h_t^S$ solve $\nabla h_t^S(b) = \eta_t(b)$.  Applying the Nash continuity and Nash-Aronson estimate (Lemma \ref{symm_rw::lem::nash_continuity_bounded}) to $h_t^S - h_t^\zeta$ and $h_t^{\zeta} - h_t^{\wt{\zeta}}$, similar to the proof of Theorem \ref{gl::thm::infinite_coupling}, yields
\begin{equation}
\label{cd::eqn::env_change} \ol{d}_\infty \equiv \sup_{0 \leq t \leq S'} \sup_{b \in E^*} |c_{T-t}(b) - \breve{c}_{T-t}(b)| = O(M_0 R^{-\sigma \xi_{\rm NC}})
\end{equation}
where $M_t = \| \ol{h}_t \|_\infty + \| h_t^\zeta - h_t^S\|_\infty$, as $\CV''$ is Lipschitz.  Here, we are taking the maximum over $x \in F$.  Let $q(s,t;x,y) = p(T-t,T-s;y,x)$ and $q'(s,t;x,y) = \wt{p}'(T-t,T-s;y,x)$.  Then $q,q'$ are the transition kernels for random walks jumping with rates $c_t(b), \wt{c}_t'(b)$, respectively.  The previous lemma thus yields
\[ \sum_{b \in E^*} \sum_{y \in F} \int_{0}^T (\nabla \ol{p}'(u,T;b,y))^2 du = O(M_0^2 R^{\epsilon/2-2\sigma \xi_{\rm NC}})\]
since the contribution to the sum given by those $y \in F$ with $\dist(y, E) \geq R^{\alpha+\epsilon/4}$ is negligible.
Since we can cover $E$ by $O(R^{2(\alpha-\gamma)})$ disks of radius $R^\gamma$, applying Lemma \ref{harm::lem::hc_time_space_deriv} to the terms corresponding to $I_1,I_2$ gives us the bound
\begin{align}
    &\sum_{b \in E^*} \sum_{y \in F} \int_{0}^{S/4} (\nabla \ol{p}(u,T;b,y))^2 du \notag\\ 
 =& O(R^{\epsilon/2+4\sigma + 2(\alpha-\gamma) - 2\gamma \xi_{\rm NC}}) + O(M_0^2 R^{\epsilon/2-2\sigma \xi_{\rm NC}}). \label{harm::eqn::hc_diff_bound}
 \end{align}
Observe that $M$ is the only random quantity on the right hand side.  Let
\[ \wt{h}_t(x) = \sum_{y \in F} \wt{p}(T-t,T;x,y) \ol{h}_0(y) \text{ for } 0 \leq t \leq T.\]
Note that
\begin{align*}
     \wt{h}_t(x) &= \sum_{y,z \in F} \wt{p}(T-t,S;x,y) p(S,T;y,z) \ol{h}_0(z)\\
                        &= \sum_{y \in F} \wt{p}(T-t,S;x,y) \ol{h}_{T-S}(z) \text{ for } T-S \leq t \leq T.
\end{align*}
Hence as $t \mapsto \wt{p}(T-t,S;x,y)$, $T-S \leq t \leq T$, is independent from $\ol{h}_{T-S}$,
it follows that $\wt{h}_{t + (T-S)} \stackrel{d}{=} \breve{h}_t$, $0 \leq t \leq S$, with $\breve{h}$ as in the statement of the proposition.  We can write
\begin{align*}
  &\E \sum_{b \in E^*} \int_{T-S/4}^T (\nabla \ol{h}_t(b) - \nabla \wt{h}_t(b))^2 dt\\
= &  \E \sum_{b \in E^*} \int_{T-S/4}^T \left( \sum_{y \in F} \nabla \ol{p}(T-t,T;b,y) \ol{h}_{T-t}(y) \right)^2 dt.
\end{align*}
The terms in the summation over $y \in F$ which are of distance at least $R^{\alpha+\epsilon/2}$ from $E$ make a negligible contribution to the summation by the Nash-Aronson estimate.  Consequently, by making a change of variables and applying the Cauchy-Schwarz inequality, we see that it suffices to control
\begin{align}
\E \left[ R^{2\alpha+\epsilon} \sum_{b \in E^*} \sum_{y \in F} \int_0^{S/4} (\nabla \ol{p}(u,T;b,y))^2  M_u^2 du \right]. \label{harm::eqn::suff_bound}
\end{align}
Choosing $\sigma = \gamma \xi_{\rm NC}/4$ and applying \eqref{harm::eqn::hc_diff_bound} yields that the expression in \eqref{harm::eqn::suff_bound} is bounded by
\begin{align*}
  & \sup_{0 \leq u \leq S/4} \bigg( \E O(M_u^2 R^{3\epsilon/2+ 4\alpha-2\gamma - \gamma \xi_{\rm NC}}) + \E O(M_0^2 M_{u}^2 R^{3\epsilon/2+2\alpha - 2\gamma \xi_{\rm NC}^2/4}) \bigg).
\end{align*}
Now, Lemma \ref{gl::lem::mom_max} and the Cauchy-Schwarz inequality imply $\E (M_u^2+M_0^2 M_{u}^2) = O_{\ol{\Lambda}}(R^{\epsilon/2})$ uniformly in $u$.  Thus, we are left with the bound
\begin{align*}
   O_{\ol{\Lambda}}(R^{2\epsilon+4\alpha - 2\gamma - \gamma \xi_{\rm NC}^2/2}).
\end{align*}
Taking $\rho_{\rm EC} = \xi_{\rm NC}^2/2$ gives \eqref{harm::eqn::ec_change_int}.  Equation \eqref{harm::eqn::ec_change_sum} follows from exactly the same argument except using the first part of Lemma \ref{harm::lem::hc_time_space_deriv} rather than the second.  The final part of the proposition is immediate from the construction and the Nash continuity estimate.
\end{proof}

We can now prove Proposition \ref{cd::prop::environment_change}.
\begin{proof}[Proof of Proposition \ref{cd::prop::environment_change}]
We now construct couplings as follows.  First, we couple $(\eta^2, \breve{h}^2), (\nabla h, \ol{h})$ as in Proposition \ref{cd::prop::environment_change} for $\gamma = \gamma_2$.  Equation \eqref{harm::eqn::ec_change_sum} implies that with $S_2 = R^{2\gamma_2}$ we have 
\begin{align}
 \E \sum_{x \in \wt{E}} |\ol{h}_{S_2}(x) - \breve{h}_{S_2}^2(x)|^2 &= O_{\ol{\Lambda}}(R^{\epsilon/2 + \delta_2+2\gamma_2}). \label{harm::eqn::ec_coupling_error1}
\end{align}
where $\wt{E} = B(x_0, 2R^{\alpha+\epsilon/100})$.  Now we apply Proposition \ref{cd::prop::environment_change} a second time except with $\gamma = \gamma_1$ and starting at $S_2$ to yield a coupling $\big((\eta^1,\breve{h}^1), (\nabla h, \ol{h}) \big)$.  Equation \eqref{harm::eqn::ec_change_int} implies the existence of $\tfrac{3}{4} R^{2\gamma_1} \leq S_1 \leq R^{2\gamma_1}$ such that with $S = S_1 + S_2$ we have
\begin{equation}
\label{harm::eqn::ec_coupling_error2}
 \E \sum_{b \in E^*} |\nabla \ol{h}_{S}(b) - \nabla \breve{h}_{S}^1(b)|^2 = 
O_{\ol{\Lambda}}(R^{\epsilon + \delta_1}).
\end{equation}
Let $\breve{p}^i$ be the kernel associated with $\eta^i$.  In this coupling, $\ol{h}_{t}$ for $0 \leq t \leq S_2$ is independent from $\breve{p}^1$ and $\breve{p}^1,\breve{p}^2$ are independent.  Let
\begin{align*}
 \breve{h}_t(x) 
&= \sum_{y,z \in F} \breve{p}^1(S-t,S_1;x,y)\breve{p}^2(0,S_2;y,z) \ol{h}_0(z)\\
&= \sum_{y \in F} \breve{p}^1(S-t,S_1;x,y) \breve{h}_{S_2}^{2}(y)
\end{align*}
for $S_2 \leq t \leq S$.
We have
\begin{align*}
   &  \E \sum_{b \in E^*} |\nabla \breve{h}_S(b) - \nabla \breve{h}_S^1(b)|^2
= \E \sum_{b \in E^*} \left(\sum_{y \in F}  \nabla \breve{p}^1(0,S_1;b,y) (\breve{h}_{S_2}^2(y) - \ol{h}_{S_2}(y)) \right)^2\\
\leq&  \left(\sum_{b \in E^*} \sum_{y \in \wt{E}} \E (\nabla \breve{p}^1(0,S_1;b,y))^2\right)  \left(\sum_{y \in \wt{E}}\E (\breve{h}_{S_2}^2(y) - \ol{h}_{S_2}(y))^2\right) + O(\exp(-R^{10^{-5}\epsilon})),
\end{align*}
where the last inequality came from Cauchy-Schwarz and the Nash-Aronson estimates (Lemma \ref{symm_rw::lem::nash_aronson}).
It follows from equation (1.4) of Theorem 1.1 from  \cite{DD05} and the Nash-Aronson estimates that 
\[ \sum_{y \in \wt{E}} \sum_{b \in E^*} \E (\nabla \breve{p}^1(0,S_1;b,y))^2 = O(R^{2\alpha+\epsilon-4\gamma_1}).\]
Combining everything proves the proposition.
\end{proof}

\subsection{Approximation by the Average}

\begin{proposition}
\label{harm::prop::average_approx}
Suppose that we have the same setup as Theorem \ref{cd::thm::cd}.  There exists $\rho_{\rm A} > 0$ depending only on $\CV$ such that the following holds.  If $\ol{h}^\rho(x)$ is the average of $\ol{h}(x)$ on the ball $B(x,R^\rho)$, then
\begin{align}
      &\E \sum_{b \in E^*} \CV''(\nabla h_{T_2}^\zeta(b)) (\nabla \ol{h}_{T_1}(b) - \nabla \ol{h}_{T_1}^\rho(b)) \nabla g(b) \notag\\
   =& O_{\ol{\Lambda}}(R^{\epsilon+\rho+\alpha(1 -\rho_{\rm A})} \| \nabla g \|_\infty) \label{harm::eqn::average_approx}
\end{align}
for $T_1 \leq T_2$.
\end{proposition}
\begin{proof}
While the proof for $T_1 \neq T_2$ does not introduce any additional technical challenges, in the interest of keeping the notation simple we will only provide the proof for the case $T_1=T_2 =S$, with $S$ from Proposition \ref{cd::prop::environment_change}.  
Let 
\[ \gamma_1 = \alpha \left( \frac{ 1 + \tfrac{3}{16} \rho_{\rm EC}}{1+ \tfrac{1}{4} \rho_{\rm EC}} \right)  \text{ and } \gamma_2 = \alpha.\]
Let $\delta_1,\delta_2$ be as in Proposition \ref{cd::prop::environment_change} with these choices of $\gamma_1,\gamma_2$.  Note that then there exists $\rho_{\rm A} > 0$ depending only on $\CV$ such that 
\begin{align*}
   \delta_1 < -\alpha \rho_{\rm A} \text{ and } 2\alpha - 4\gamma_1 + \delta_2 + 2 \gamma_2 < - \alpha \rho_{\rm A}.
\end{align*}
Consequently, by Proposition \ref{cd::prop::environment_change} it suffices to prove \eqref{harm::eqn::average_approx} with $\breve{h}, \breve{h}^\rho$ in place of $\ol{h}, \ol{h}^\rho$ where $\breve{h}_t^{\rho}(x)$ is the average of $\breve{h}_t(x)$ on $B(x,R^\rho)$.  Moreover, \eqref{harm::eqn::ec_change_grad_single} combined with Lemma \ref{gl::lem::grad_error} imply that it suffices to prove \eqref{harm::eqn::average_approx} with $\CV''(\nabla h_{S}^\zeta(b))$ replaced with $\breve{c}^1(b) = \CV''(\eta_{S_1}^1)$.  For $f \colon F \to \R$ let
\[ \CL  f(x) = \sum_{b \ni x} \breve{c}^1(b) \nabla f(b).\]
Summation by parts implies it suffices to bound
\begin{equation}
\label{harm::eqn::avg_err_approx_sbp}
 \sum_{b \in \partial E^*} \breve{c}^1(b) (\breve{h}_{S}(x_b) - \breve{h}_{S}^\rho(x_b)) \nabla g(b) - \sum_{x \in E} [\breve{h}_{S}(x) - \breve{h}_{S}^\rho(x)] \CL g(x).
 \end{equation}
By the Nash continuity estimate (Lemma \ref{symm_rw::lem::nash_continuity_bounded}), we know that
\[ \E | \breve{h}_{S}(x_b) - \breve{h}_{S}^\rho(x_b)| = O_{\ol{\Lambda}}(R^{\epsilon+(\rho-\alpha)\xi_{\rm NC}}),\]
hence the boundary term in \eqref{harm::eqn::avg_err_approx_sbp} is of order $O_{\ol{\Lambda}}(R^{\epsilon+(\rho-\alpha)\xi_{\rm NC}+\alpha}\| \nabla g \|_\infty) = O_{\ol{\Lambda}}(R^{\epsilon + \alpha(1-\rho_{\rm A})} \| \nabla g \|_\infty)$, shrinking $\rho_{\rm A}$ if necessary.

We now deal with the interior term.  For $y, \theta \in \Z^2$, let $y^\theta = y+\theta$.  We are going to omit the times when referring to $\breve{p}^i$ and just write $\breve{p}^i(x,y)$ for $\breve{p}^i(0,S_i;x,y)$.  Say that a bond $b$ is ``positively oriented'' if it points either in the positive hortzonal or vertical directions.  For each triple $(x,y,b)$, let $y_b = y$ if $b$ is positively oriented and $2x-y$ otherwise.  The latter is the reflection of $y$ about $x$.  With $B_\rho = B(0,R^{\rho})$, we can rewrite the interior term of \eqref{harm::eqn::avg_err_approx_sbp} as
\begin{align*}
  \frac{1}{|B_\rho|} \E \sum_{x \in E} \sum_{b \ni x} \sum_{y,z \in F} \sum_{\theta \in B_\rho} &\breve{c}^1(b)\big[ \breve{p}^1(x,y_b) \breve{p}^2(y_b,z)  -\\
     &\breve{p}^1(x_b^\theta,y_b)\breve{p}^2(y_b,z) \big] \ol{h}_0(z) \nabla g(b).
\end{align*}
Using the independence properties of $\eta^i, \ol{h}_0$ as well as the $\Delta^\beta$-harmonicity of $g$, we see that this is the same as
\begin{align*}
  \frac{1}{|B_\rho|}  \sum_{x ,b,y,z,\theta} &w(x,y,b,\theta) \E \big[\breve{p}^2(y_b,z) - \breve{p}^2(x,z)\big] \E[\ol{h}_0(z)] \nabla g(b)
\end{align*}
where the summation is over $x \in E$, $b \ni x$, $y,z \in F$, $\theta \in B_\rho$ and
\[ w(x,y,b,\theta) = \E\big[\breve{c}^1(b) (\breve{p}^1(x,y_b) - \breve{p}^1(x_b^\theta,y_b)) \big].\]
Here we are crucially using that $\breve{p}^2(x,z)$ does not depend on $y$.
For $b \ni x$ let $b' \ni x$ have the opposite orientation of $b$.  By Proposition \ref{gl::prop::reflection_invariance}, we have
\[ w(x,y,b,\theta) =  w(x,y,b',\theta)\]
Consequently, we can rewrite our sum as
\begin{align*}
  &\frac{1}{|B_\rho|}  \sum_{x ,b,y,z,\theta} w(x,y,b,\theta) \bigg( \E \big[(\breve{p}^2(y_b,z) - \breve{p}^2(x,z)) + (\breve{p}^2(y_{b'},z) - \breve{p}^2(x,z)) \big] \nabla g(b) \\
   &+ \E \big[\breve{p}^2(y_{b'},z) - \breve{p}^2(x,z) \big] (\nabla g(b') - \nabla g(b))\bigg) \E[\ol{h}_0(z)]\\
 \equiv& \frac{1}{|B_\rho|}  \sum_{x ,b,y,z,\theta} w(x,y,b,\theta) \bigg( A(x,y,z,b) \nabla g(b) +  B(x,y,z,b) (\nabla g(b') - \nabla g(b))\bigg) \E[\ol{h}_0(z)]
\end{align*}
where the summation is now only over positively oriented bonds. 

We will deal with the term involving $A$ first.  By the Nash-Aronson estimates (Lemma \ref{symm_rw::lem::nash_aronson}), the sum over $y$ of distance from $x$ more than $R^{\gamma_1+\epsilon}$ from $x$ is negligible.  Similarly, we may ignore those $z$ with $|z-x| \geq R^{\gamma_2 + \epsilon}$.  For $|y-x| \leq R^{\gamma_1+\epsilon}$,  it is a consequence of Theorem 1.1 equation (1.5b) of \cite{DD05} and the Nash-Aronson estimates (Lemma \ref{symm_rw::lem::nash_aronson}) that
\begin{equation}
\label{harm::eqn::double_deriv}
 |A(x,y,z,b)| \leq CR^{2\gamma_1+2\epsilon - 4 \gamma_2}.
\end{equation}
Indeed, this can be seen by rewriting the difference as a sum of $O(R^{2\gamma_1+2\epsilon})$ discrete second derivatives.  Using that $\E \ol{h}_0(z) = O_{\ol{\Lambda}}(R^{\epsilon})$ from Lemma \ref{gl::lem::mom_max}, we thus have
\begin{align*}
 &\frac{1}{|B_\rho|}  \sum_{x ,b,y,z,\theta} w(x,y,b,\theta) A(x,y,z,b) \nabla g(b)\\
=& \frac{1}{|B_\rho|}  \sum_{x,b,y,\theta} w(x,y,b,\theta) O_{\ol{\Lambda}}(R^{5\epsilon + 2 \gamma_1 - 2\gamma_2} \| \nabla g \|_\infty).
\end{align*}
Now, equation (1.4) of Theorem 1.1 of \cite{DD05} implies 
\[ \E \sum_{y} |w(x,y,b,\theta)| = O(R^{\rho-\gamma_1})\]
uniformly in $x$ since $|x-x_b^\theta| \leq R^{\rho}$.  Since the sum over $x$ includes $O(R^{2\gamma_2})$ terms, our total error is $O_{\ol{\Lambda}}(R^{10 \epsilon + \rho + \gamma_1} \| \nabla g \|_\infty)$.

We now turn to the term involving $B(x,y,z,b)$.  Since $g$ is $\Delta^\beta$-harmonic so is $x \mapsto \nabla g((x,x+e_i))$ hence we have that $\nabla g(b') - \nabla g(b) = O(R^{-\gamma_2} \| \nabla g \|_\infty)$.  Again applying equation (1.4) of Theorem 1.1 of \cite{DD05} we thus see
\[ B(x,y,z,b) (\nabla g(b') - \nabla g(b)) \E[\ol{h}_0(z)] = O_{\ol{\Lambda}}(R^{3\epsilon+\gamma_1 - 4\gamma_2} \| \nabla g \|_\infty).\]
This is of an even smaller magnitude than the corresponding term with $A$, so we also get an error of $O_{\ol{\Lambda}}(R^{10 \epsilon + \rho + \gamma_1} \| \nabla g \|_\infty)$.  The proposition now follows from the explicit form of $\gamma_1$.
\end{proof}

\subsection{Change of Time}

\begin{proposition}
\label{cd::prop::change_of_time}
Suppose that we have the same setup as Theorem \ref{cd::thm::cd}.  Let $T = R^{2\gamma}$ for $0 \leq \gamma \leq \alpha$ and fix $g \colon F \to \R$.  We have that
\begin{align}
       &\E \sum_{b \in E^*} \CV''(\nabla h_{T}^\zeta(b)) (\nabla \ol{h}_{T}(b) - \nabla \ol{h}_0(b)) \nabla g(b) \label{harm::eqn::change_of_time}\\
    =& O_{\ol{\Lambda}}(R^{\epsilon+\alpha(1-\rho_{\rm CoT}) + \gamma}\| \nabla g\|_\infty) \notag
\end{align}
where $\rho_{\rm CoT} > 0$ depends only on $\CV$.
\end{proposition}

\begin{proof}
For $\rho > 0$, Proposition \ref{harm::prop::average_approx} implies that replacing $\ol{h}_T, \ol{h}_0$ by $\ol{h}_T^\rho, \ol{h}_0^\rho$, respectively, in \eqref{harm::eqn::change_of_time} introduces an error of 
\begin{equation}
\label{harm::eqn::cotexp1}
 O_{\ol{\Lambda}}(R^{\epsilon+\rho+(1-\rho_{\rm A})\alpha} \| \nabla g\|_\infty).
\end{equation}
We will now prove that
\begin{align}
   &\E \sum_{b \in E^*} \CV''( \nabla h_T^{\zeta}(b)) (\nabla \ol{h}_T^\rho(b) - \nabla \ol{h}_0^\rho(b)) \nabla g(b) \notag\\
 = &O_{\ol{\Lambda}}(T R^{\epsilon+\alpha-\rho} \| \nabla g\|_\infty) \label{harm::eqn::cotexp2}.
\end{align}
Since $\sum_{b \in E^*} (\nabla g(b))^2 = O(R^{2\alpha} \| \nabla g \|_\infty^2)$, applying the Cauchy-Schwarz inequality to the expression on the left hand side implies that it suffices to show 
\[ \E \sum_{b \in E^*} (\nabla \ol{h}_T^\rho(b) - \nabla \ol{h}_0^\rho(b))^2 = O_{\ol{\Lambda}}(T^2 R^{\epsilon-2\rho}).\]  As
\[ |\partial_t \ol{h}_t^\rho(x)| \leq \frac{C_0}{R^{2\rho}} \sum_{b \in \partial B^*(x,R^{\rho})} |\nabla \ol{h}_t(b)|\]
we have that
\begin{align*}
&     \E \sum_{b \in E^*} (\nabla \ol{h}_T^\rho(b) -  \nabla \ol{h}_0^\rho(b))^2
 \leq C_1T \E \int_0^T \sum_{b \in E^*} |\nabla \partial_t \ol{h}_t^\rho(b)|^2\\
 \leq& \frac{C_2T}{R^{2\rho}} \E \int_0^T \sum_{b \in E_1^*} |\nabla \ol{h}_t(b)|^2
    = \frac{C_2T^2}{R^{2\rho}} \E \sum_{b \in E_1^*} |\nabla \ol{h}_0(b)|^2,
\end{align*}
where $E_1 = B(x_0,R^{\alpha} + R^{\rho})$ and the final equality comes by stationarity.  We know that the latter quantity is of order $O_{\ol{\Lambda}}(R^{\epsilon})$ by Lemma \ref{gl::lem::mom_max}.  Equating the exponents in \eqref{harm::eqn::cotexp1}, \eqref{harm::eqn::cotexp2} gives the equation
\[ \rho + (1-\rho_{\rm A}) \alpha = 2\gamma + \alpha - \rho,\]
which leads to the choice
\[ \rho = \gamma + \frac{\rho_{\rm A}}{2} \alpha.\]
Combining everything gives the proposition, where $\rho_{\rm CoT} = \tfrac{1}{2} \rho_{\rm A}$.
\end{proof}

\subsection{Proof of Theorem \ref{cd::thm::cd}}
Assume that $T = R^{2\gamma}$ where $\gamma = \alpha \rho_{\rm CoT}/2$.
By Proposition \ref{cd::prop::change_of_time} it suffices to estimate
\begin{equation}
\label{cd::eqn::cd_reduction}
 \sum_{b \in E^*} \E[\CV''(\nabla h_{T}^\zeta(b))  \nabla \ol{h}_0(b))] \nabla g(b)
\end{equation}
provided we pay an error of $O_{\ol{\Lambda}}(R^{\epsilon+\alpha(1- \rho_{\rm CoT}/2)} \| \nabla g \|_\infty)$.  Now the result follows from an argument similar to the proof of \eqref{cd::eqn::env_change} from the proof of Proposition \ref{cd::prop::environment_change}.  Indeed, we let $\eta_t$ follow the $\SEGGS_u$ dynamics driven by the same Brownian motions as $(h_t^\zeta, h_t^{\wt{\zeta}})$ but with $\eta_0$ independent from $(h_0^\zeta,h_0^{\wt{\zeta}})$, then use the Nash continuity estimate (Lemma \ref{symm_rw::lem::nash_continuity_bounded}) to argue that 
\[ \p[ \max_{b \in E^*} |\nabla h_T^\zeta(b) - \eta_T(b)| \geq R^{\epsilon-\gamma \xi_{\rm NC}}] = O_{\ol{\Lambda}}(R^{-100}).\]
Thus applying the Cauchy-Schwarz inequality, we see that replacing $\CV''(\nabla h_T^\zeta(b))$ with $\CV''(\eta_T(b))$ in \eqref{cd::eqn::cd_reduction} introduces an error of order
\[ \left( \sum_{b \in E^*} \E[ (\nabla \ol{h}_0(b))^2]\right)^{1/2} \cdot O_{\ol{\Lambda}}(R^\alpha \| \nabla g \|_\infty) \cdot O_{\ol{\Lambda}}(R^{\epsilon - \gamma \xi_{\rm NC}}).\]
The result now follows as the first term is of order $O_{\ol{\Lambda}}(R^{\epsilon})$ by Lemma \ref{gl::lem::grad_error}.
\qed

\section{Harmonic Coupling}
\label{sec::harm}

Throughout this section we will make use of the following notation.  For $F \subseteq \Z^2$ and $\phi \colon \partial F \to \R$, we let $\p_F^\phi$ be the law of the GL model on $F$ with boundary condition $\phi$.  We will denote by $h^\phi$ a random variable distributed according to $\p_F^\phi$, where $F$ is understood through the domain of definition of $\phi$.  Finally, for $g \colon F \to \R$ we let $\Q_F^{\phi,g}$ be the law of $(h^\phi - g)$.  We will write $\E^\psi$ for the expectation under $\p_F^\psi$ if we want to emphasize $\psi$ and, similarly, $\E^{\psi,\wt{\psi}}$ for the expectation under a coupling of $\p_F^\psi, \p_F^{\wt{\psi}}$ if we want to emphasize both $\psi$ and $\wt{\psi}$.  We also let $\beta = \beta(u)$ as in the statement of Theorem \ref{thm::clt} for a fixed tilt $u \in \R^2$.

Suppose that $\mu,\nu$ are measures with $\mu$ absolutely continuous with respect to $\nu$.  Recall that the relative entropy of $\mu$ with respect to $\nu$ is the quantity
\[ \h(\mu|\nu) = \E_\mu \left[ \log \frac{d\mu}{d\nu} \right].\]
We begin by fixing $D \subseteq \Z^2$ with $R = \diam(D) < \infty$.  Fix $\ol{\Lambda} > 0$ and let $\psi,\wt{\psi} \in \B_{\ol{\Lambda}}^u(D)$. 
Morally, the idea of our proof is to get an explicit upper bound on the rate of decay of the symmetrized relative entropy
\[ \h(\p^{\wt{\psi}}|\Q^{\psi,\wh{h}}) + \h(\Q^{\psi,\wh{h}}|\p^{\wt{\psi}}),\]
where $\wh{h}$ is the $\Delta^\beta$-harmonic extension of $\psi - \wt{\psi}$ from $\partial D \to D$, as $R \to \infty$, then invoke Pinsker's inequality, the well-known bound that the total variation distance of measures is bounded from above by the square-root of their relative entropy \cite{DZ98}:
\begin{equation}
\label{harm::eqn::pinsker} \| \mu - \nu \|_{TV}^2 \leq \tfrac{1}{2} \h(\mu|\nu) .
\end{equation}
We will show shortly that the symmetrized relative entropy takes the form
\begin{equation}
\label{harm::eqn::entropy_form}
\sum_{b \in D^*} \E^{\psi,\wt{\psi}} c(b) \nabla \wh{h}(b) (\nabla \wh{h}(b) - \nabla \ol{h}(b))
\end{equation}
where $\ol{h} = h^\psi - h^{\wt{\psi}}$ and $c(b)$ is a collection of conductances which are \emph{random} but uniformly bounded from above and below.  In the Gaussian case, $c(b) \equiv c$ is constant, hence one can sum by parts, then use the harmonicity of $\wh{h}$ to get that the entropy vanishes.  The idea of our proof is to use Theorem \ref{cd::thm::cd} repeatedly to show that this approximately holds in expectation:
\begin{align}
\label{harm::eqn::entropy_approx_constant}
   &\sum_{b \in D^*} \E^{\psi,\wt{\psi}} c(b) \nabla \wh{h}(b)(\nabla \wh{h}(b)-\nabla \ol{h}(b))\notag\\
=& \sum_{b \in D^*} a_u(b) \E^{\psi,\wt{\psi}} \nabla \wh{h}(b) (\nabla\wh{h}(b)-\nabla\ol{h}(b)) + O_{\ol{\Lambda}}(R^{-\delta})
\end{align}
for some $\delta > 0$, where $a_u(b) = \E[ \CV''(\eta(b))]$ for $\eta \sim \SEGGS_u$.  Note that $a_u(b)$ depends only on $\CV,u,$ and the orientation of $b$.

Theorem \ref{cd::thm::cd} is only applicable if the distance of $b$ to $\partial D$ is $\Omega(R^\xi)$.  This will force us to deal with a boundary term, the magnitude of which will in turn depend on the regularity of both $\nabla \ol{h}$ and $\wh{h}$ near $\partial D$.  Since we make no hypotheses on $\psi, \wt{\psi}$ other than being pointwise bounded it may very well be that neither $\ol{h}$ nor $\wh{h}$ possess any regularity near $\partial D$.  

We will resolve this issue by invoking the length-area comparison technique of Section \ref{sec::dyn}.  This gives us that, for $\epsilon > 0$ fixed, there exists $R^{1-\epsilon} \leq R_D^1 \leq R_D^2 \equiv R_D^1 + R^{1-5\epsilon} \leq 2R^{1-\epsilon}$ such that
\[ \E^{\psi,\wt{\psi}} \sum_{b \in D^*(R_D^1, R_D^2)} |\nabla \ol{h}(b)|^2 = O_{\ol{\Lambda}}(R^{-\epsilon})\]
where $D(R_1,R_2) = \{ x \in D : R_1 \leq \dist(x,\partial D) < R_2\}$.  Let $g$ be the $\Delta^\beta$-harmonic extension of $\ol{h}$ from $\partial D(R_D^1, R_D^2)$ to $D(R_D^1,R_D^2)$.  Note that $g$ is the minimizer of the variational problem
\[    \wh{g} \mapsto \sum_{b \in D^*(R_D^1, R_D^2)}  a_u(b) (\nabla \wh{g}(b))^2,\ \ \wh{g}(x) = \ol{h}(x) \text{ for } x \in \partial D(R_D^1,R_D^2).\]
Indeed, the first order conditions for optimality are exactly that $\Delta^\beta \wh{g}(x) = 0$ for $x \in D(R_D^1, R_D^2)$.  Consequently, 
\[ \E^{\psi,\wt{\psi}} \sum_{b \in D^*(R_D^1, R_D^2)} |\nabla g(b)|^2 = O_{\ol{\Lambda}}(R^{-\epsilon}).\]
Going back to \eqref{harm::eqn::entropy_form}, by invoking Pinsker's inequality, this implies that we can construct our initial coupling so that $\ol{h}$ is harmonic in $D(R_D^1,R_D^2)$ on an event $\CH$ with $\p[\CH] = 1-O_{\ol{\Lambda}}(R^{-\epsilon/2})$.  On $\CH$, we have that $\nabla \ol{h}(b) = O_{\ol{\Lambda}}(R^{6\epsilon-1})$ uniformly in $b \in \partial D(R_D)$ where $R_D = R_D^1 + \tfrac{1}{2} R^{1-5\epsilon}$.  Thus with high probability $\ol{h}$ has plenty of regularity a bit away from $\partial D$ while $\psi - \wt{\psi}$ need not have any.

Moving to the subdomain $D(R_D)$ from $D$ is also useful since it possesses the $r$-exterior ball property for $r = R_D$.  This means that for every $x \in \partial D(R_D)$ there exists $y \notin D(R_D)$ such that $B(y,R_D) \cap D(R_D) = \emptyset$ and $x \in \partial B(y,R_D)$.  The importance of this property is that it implies pointwise regularity of harmonic functions in $D(R_D)$ near $\partial D(R_D)$, more so than one has for such functions in $D$ near $\partial D$ without further hypotheses.  This is related to the notion of ``stochastic regularity'' \cite{KS98} and that random walk ``exits much more quickly'' from such domains.

\begin{figure}
     \centering
     \subfigure[Domain without the exterior ball property]{
          \includegraphics[width=.30\textwidth]{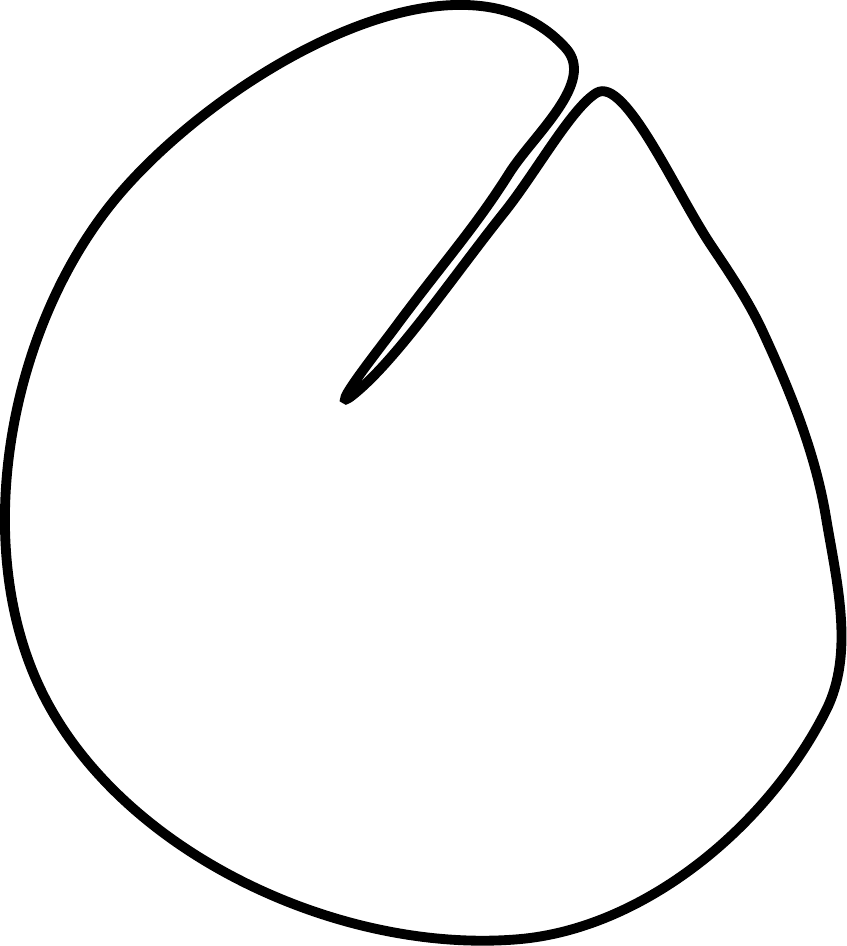}}
          \hspace{0.1in}
     \subfigure[Domain with the $r$-exterior ball property]{
          \includegraphics[width=.30\textwidth]{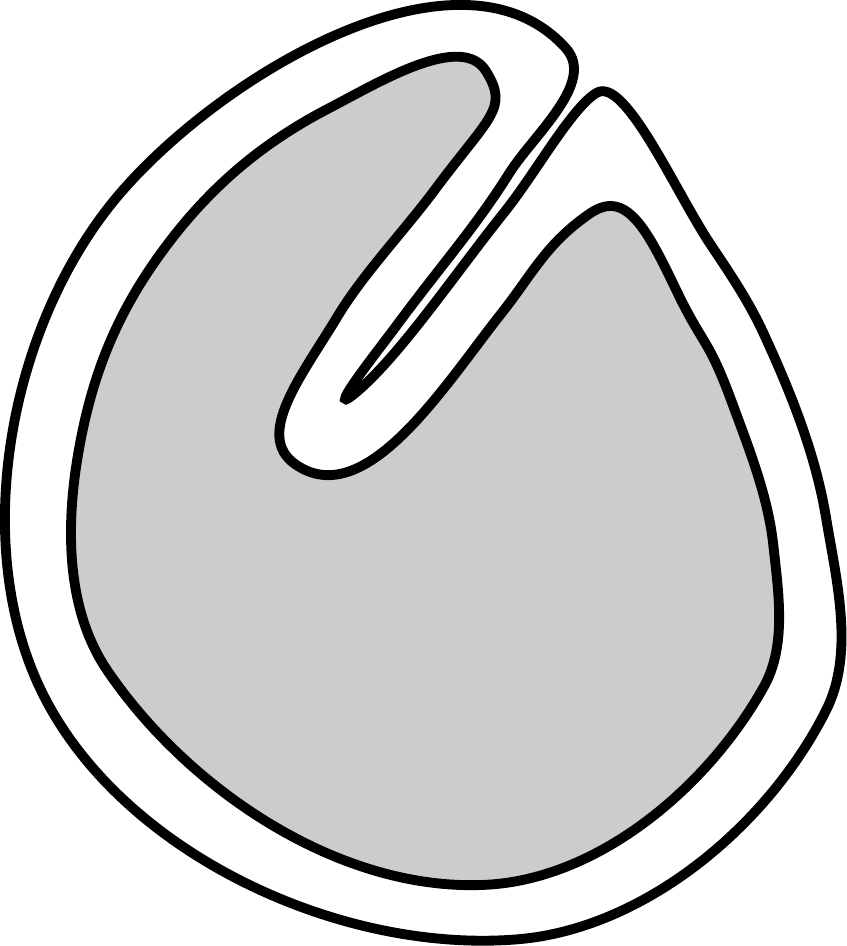}}
          \caption{The domain on the left hand side does not have the $r$-exterior ball property for due to the fjord in its upper right corner.  However, the inner domain $D(r) = \{ x \in D: \dist(x,\partial D) \geq r\}$ shaded in light grey on the right hand side trivially does possess this property.  The distinction is important since the $r$-exterior ball property is related to the regularity near the boundary of discrete harmonic functions.}
\end{figure}

The rest of this section is organized as follows.  In subsection \ref{harm::subsec::entropy}, we will justify \eqref{harm::eqn::entropy_form}.  The purpose of subsection \ref{subsec::entropy_est} is to prove a form of \eqref{harm::eqn::entropy_approx_constant}.  Finally, in subsection \ref{harm::subsec::proof} we will put everything together to prove Theorems \ref{harm::thm::coupling} and \ref{harm::thm::mean_harmonic}.

Before we proceed, we would like to emphasize that in this section we will prove that a certain symmetrized relative entropy decays like a small, negative power of $R = \diam(D)$.  Because of this, many of our estimates will be accurate only up to very small powers of $R$.  In particular, we do not try to derive the ``best possible'' exponents and make many cautious choices in order to avoid carrying around overly complicated exponents.

\subsection{The Symmetrized Relative Entropy}
\label{harm::subsec::entropy}

\begin{lemma}
\label{harm::lem::entropy_form}
Suppose that $F \subseteq \Z^2$ is bounded and $\zeta,\wt{\zeta} \colon \partial F \to \R$ are given boundary conditions.  If $g \colon F \to \R$ is any function such that $g|_{\partial F} = \zeta - \wt{\zeta}$ then
\begin{align*}
& \h(\p_F^{\wt{\zeta}}| \Q_F^{\zeta,g}) + \h(\Q_F^{\zeta,g}|\p_F^{\wt{\zeta}})\\
=& \sum_{b \in F^*} \E^{\zeta,\wt{\zeta}} \bigg[ \CV''(\nabla h^{\zeta}(b)) \nabla g(b) \nabla(g-\ol{h})(b) + O( (|\nabla \ol{h}(b)|^2 + |\nabla g(b)|^2)|\nabla g(b)|) \bigg].
\end{align*}
where $\E^{\zeta,\wt{\zeta}}$ denotes the expectation under any coupling of $\p_F^{\zeta},\p_F^{\wt{\zeta}}$ and $\ol{h} = h^\zeta - h^{\wt{\zeta}}$.
\end{lemma}
\begin{proof}
The densities $p,q=q_g$ of $\p_F^{\wt{\zeta}}$ and $\Q_F^{\zeta,g}$ with respect to Lebesgue measure are given by
\begin{align*}
p(h) &= \frac{1}{\CZ_p} \exp\left(-\sum_{b \in F^*} \CV( \nabla (h \vee \wt{\zeta})(b))\right),\\
q(h) &= \frac{1}{\CZ_q} \exp\left(-\sum_{b \in F^*} \CV( \nabla [(h + g) \vee \zeta](b))\right).
\end{align*}

With $\E^{\zeta,g}$ the expectation under $\Q_F^{\zeta,g}$, we have
\begin{align*}
   \h(\Q_F^{\zeta,g}|\p_F^{\wt{\zeta}}) + \log \frac{\CZ_q}{\CZ_p}
&= \sum_{b \in F^*} \E^{\zeta,g} \bigg[\CV( (\nabla h \vee \wt{\zeta})(b)) - \CV( \nabla ((h + g) \vee \zeta)(b)) \bigg]\\
&= -\sum_{b \in F^*} \left( \E^{\zeta} \int_0^1 \CV'(\nabla [(h + (s-1)g) ](b))ds  \right) \nabla g(b)
\end{align*}
Similarly, with $\E^{\wt{\zeta}}$ the expectation under $\p_F^{\wt{\zeta}}$,
\begin{align*}
   \h(\p_F^{\wt{\zeta}}|\Q_F^{\zeta,g}) + \log \frac{\CZ_p}{\CZ_q}
&= \sum_{b \in F^*} \E^{\wt{\zeta}} \bigg[\CV(\nabla ((h + g) \vee \zeta)(b)) - \CV( (\nabla h \vee \wt{\zeta})(b)) \bigg]\\
&= \sum_{b \in F^*} \left( \E^{\wt{\zeta}} \int_0^1 \CV'(\nabla (h + s g)(b)) ds \right) \nabla g(b)
\end{align*}
Thus
\begin{align*}
 &  \h(\p_F^{\wt{\zeta}}|\Q_F^{\zeta,g}) + \h(\Q_F^{\zeta,g}|\p_F^{\wt{\zeta}})\\
=&\sum_{b \in F^*} \left( \E^{\wt{\zeta}} \int_0^1 \CV'(\nabla (h + s g)(b)) ds - \E^{\zeta} \int_0^1 \CV'(\nabla [(h + (s-1)g) ](b))ds\right) \nabla g(b)\\
=& \sum_{b \in F^*} \left( \E^{\zeta,\wt{\zeta}} \int_0^1 \int_0^1 \CV''(\nabla (h^{\zeta} + (s-1) g)(b) + r \nabla (g-\ol{h})(b)) dr ds\right)\\
&\ \ \ \ \ \ \ \ \ \ \ \ \ \ \ \nabla g(b) \nabla (g-\ol{h})(b).
\end{align*}
As $\CV''$ is Lipschitz,
\begin{align*}
  & \int_0^1 \int_0^1 \CV''(\nabla (h^{\zeta} + (s-1) g)(b) + r \nabla (g-\ol{h})(b)) dr ds\\
=& \CV''(\nabla h^{\zeta}(b)) + O(\nabla \ol{h}(b)) + O(\nabla g(b)).
\end{align*}
The lemma now follows by an application of Cauchy-Schwarz.
\end{proof} 

\subsection{Estimating the Entropy}
\label{subsec::entropy_est}

Suppose that $E \subseteq \Z^2$ with $\diam(E) = R$ and $\zeta,\wt{\zeta} \in \B_{\ol{\Lambda}}^u(E)$.  We are now going to give a general estimate of the symmetrized relative entropy of the previous lemma when $g \colon E \to \R$ is the $\Delta^\beta$-harmonic extension of $\ol{\zeta} = \zeta - \wt{\zeta}$ from $\partial E$ to $E$.  The error is going to be a function of the regularity of $g$, $\ol{\zeta}$, and the number of balls required to cover annuli near $\partial E$.  When all of the boundary data is smooth, the error is actually negligible.
To this end, we let
\[ \| \ol{\zeta}\|_{\nabla}^E = \max_{x,y \in \partial E} \frac{|\ol{\zeta}(x) - \ol{\zeta}(y)|}{|x-y|}.\]
Fix $\epsilon > 0$, let $\gamma_k = \epsilon k$, $M$ be the largest integer so that $\gamma_M < 1$, and $N_k$ be the number of balls of radius $R^{\gamma_k}$ necessary to cover $E(R^{\gamma_k}, R^{\gamma_{k+1}}) = \{x \in E : R^{\gamma_k} \leq \dist(x, \partial E) < R^{\gamma_{k+1}}\}$.  Finally, with $1 \leq \ell \leq M$ fixed, let
\begin{align*}
  \CE_E &= |E^*| \| \nabla g \|_\infty^3 + \CE \| \nabla g \|_\infty,\\
  \CE_B^\ell &= |(E')^*(R^{\gamma_\ell})| \| \nabla g \|_\infty \CE^\ell,\\
  \CE_I^\ell &= \sum_{k=\ell+1}^M N_k R^{\epsilon + \gamma_k(1-\rho_{\rm CD})} \| \nabla g \|_\infty,
\end{align*}
where
\[ \CE^\ell = ( R^{1/2} \| \ol{\zeta}\|_\nabla^E + R^{(\gamma_\ell - 1/2)\rho_{\rm B}}\| \ol{\zeta}\|_\infty ) \text{  and  } \CE = |\partial E^*| \| \ol{\zeta}\|_\infty \CE^1.\]
\begin{proposition}
\label{harm::prop::ent_est}
We have that
\[ \h(\p_E^\zeta | \p_E^{\wt{\zeta}}) +  \h( \p_E^{\wt{\zeta}} | \p_E^\zeta) = O_{\ol{\Lambda}}(\CE_E + \CE_I + \CE_B).\]
\end{proposition}
\begin{proof}
Let $(h^\zeta,h^{\wt{\zeta}})$ be the stationary coupling of $\p_E^\zeta, \p_E^{\wt{\zeta}}$ and $\ol{h} = h^\zeta - h^{\wt{\zeta}}$.  We begin with an \emph{a priori} estimate on the Dirichlet energy of $\ol{h}$.  First fix $b \in \partial E^*$.  Lemma \ref{symm_rw::lem::beurling} implies that $\nabla \ol{h}(b) = O(R^{1/2}  \| \ol{\zeta}\|_\nabla^E + R^{-\rho_{\rm B}/2} \| \ol{\zeta}\|_\infty)$.
Applying Lemma \ref{gl::lem::ee} to the stationary dynamics $(h_t^\zeta,h_t^{\wt{\zeta}})$, we have that
\begin{align*}
   & \sum_{b \in E^*}  |\nabla \ol{h}(b)|^2 
\leq C \sum_{b \in \partial E^*} |\ol{\zeta}(x_b)| |\nabla \ol{h}(b)|.
\end{align*}
Consequently, we have
\begin{align}
   \sum_{b \in E^*}  |\nabla \ol{h}(b)|^2 
\leq C\CE. \label{harm::eqn::dir_en_bound}
\end{align}

We now break  the right hand side in the statement of Lemma \ref{harm::lem::entropy_form} into three terms:

\begin{align}
 \sum_{b \in E^*} &\E^{\zeta,\wt{\zeta}} (|\nabla g(b)|^2 + |\nabla \ol{h}(b)|^2)|\nabla g(b)| \label{harm::eqn::ent_extra},\\
 \sum_{b \in (E')^*(R^{\gamma_\ell})} &\E^{\zeta,\wt{\zeta}} \bigg[ \CV''(\nabla h^{\zeta}(b)) \nabla g(b) \nabla(g-\ol{h})(b) \bigg]  \label{harm::eqn::ent_boundary},\\
 \sum_{b \in E^*(R^{\gamma_\ell})} &\E^{\zeta,\wt{\zeta}} \bigg[ \CV''(\nabla h^{\zeta}(b)) \nabla g(b) \nabla(g-\ol{h})(b) \bigg] \label{harm::eqn::ent_interior}.
\end{align}

\noindent{\it Estimate of \eqref{harm::eqn::ent_extra}}

The first term in the summation is trivially bounded by $|E^*| \| \nabla g \|_\infty^3$.  By \eqref{harm::eqn::dir_en_bound} the second is at most $\|\nabla g \|_\infty \CE.$  This gives an error of $O(\CE_E)$.

\noindent{\it Estimate of \eqref{harm::eqn::ent_boundary}}

The error is easily seen to be $|(E')^*(R^{\gamma_\ell})| \|\nabla g\|_\infty \delta$ where $\delta = \max_{x \in E'(R^{\gamma_\ell})} | g(x) - \ol{h}(x)|$, so we just need to estimate $\delta$.  Fix $x \in E'(R^{\gamma_\ell})$.  We know that we can write $g(x) = \E \ol{\zeta} (X_\tau)$ where $X$ is a random walk initialized at $x$ with bounded rates and $\tau$ its first exit from $E$.  By Lemma \ref{symm_rw::lem::beurling}, the probability that $X_\tau$ exits at distance less than $R^{1/2}$ from $x \in E'(R^{\gamma_\ell})$ is $1-O(R^{ (\gamma_\ell - 1/2) \rho_{\rm B}})$.   Now, $\ol{h}$ admits a similar representation though the random walk has time-varying rates.  Nevertheless, the same statement still holds.  Consequently, $|\ol{h}(x) - g(x)| \leq \CE^\ell$, which leads to the desired bound.  This gives an error of $O(\CE_B^\ell)$.

\noindent{\it Estimate of \eqref{harm::eqn::ent_interior}}

We break the summation over $E^*(R^{\gamma_\ell})$ into the annuli $E_k = E^*(R^{\gamma_k}, R^{\gamma_{k+1}})$.  Each annulus can be covered by $N_k$ balls of radius $R^{\gamma_k}$ by hypothesis.  Let $B$ be such a ball.  Applying Theorem \ref{cd::thm::cd} on $B$ yields an error of $O_{\ol{\Lambda}}(R^{\epsilon + \gamma_k(1-\rho_{\rm CD})} \| \nabla g \|_\infty)$.  Repeating this on all $N_k$ balls for $\ell+1 \leq k \leq M$ implies that \eqref{harm::eqn::ent_interior} is equal to
\begin{align*}
\sum_{b \in E^*(R^{\gamma_\ell})} a_u(b) \E[ \nabla(\ol{h} - g)(b)] \nabla g(b) + 
\sum_{k=\ell}^M N_k O_{\ol{\Lambda}}(R^{\epsilon + \gamma_k(1-\rho_{\rm CD})} \| \nabla g \|_\infty)
\end{align*}
where $a_u(b) = \E[ \CV''(\eta(b))]$, $\eta \sim \SEGGS_u$.  Note that the error term is precisely $\CE_I^\ell$.  By our bound of \eqref{harm::eqn::ent_boundary}, we can rewrite the above expression as
\begin{align*}
\sum_{b \in E^*} a_u(b) \E[ \nabla(\ol{h} - g)(b)] \nabla g(b) + 
O_{\ol{\Lambda}}(\CE_E + \CE_I).
\end{align*}
By summation by parts, we see that this is equal to
\begin{align*}
\sum_{b \in E^*} \E[ (\ol{h} - g)(x)] \Delta^\beta g(x) + 
O_{\ol{\Lambda}}(\CE_E + \CE_I) = O_{\ol{\Lambda}}(\CE_E + \CE_I),
\end{align*}
where $\beta = \beta(u)$ as $g$ is $\Delta^\beta$-harmonic.  This gives an error of $O(\CE_I^\ell)$.
\end{proof}

\subsection{Proof of Theorems \ref{harm::thm::coupling} and \ref{harm::thm::mean_harmonic}}
\label{harm::subsec::proof}

\subsubsection{The Initial Coupling}

\begin{figure}
     \centering
     \subfigure[Stage 1: Langevin Dynamics in $D$]{
          \includegraphics[width=.28\textwidth]{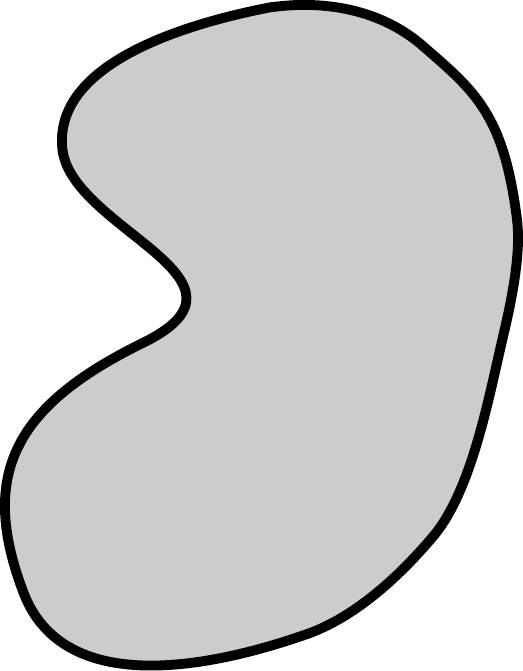}}
          \hspace{0.1in}
     \subfigure[Stage 2: Harmonic in $D(R_D^1, R_D^2)$]{
          \includegraphics[width=.28\textwidth]{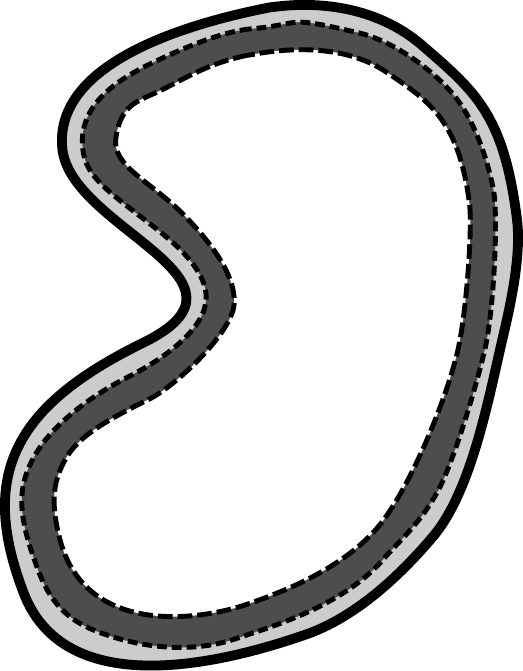}}
          \hspace{0.1in}
      \subfigure[The domain $E = D(R_D)$]{
          \includegraphics[width=.28\textwidth]{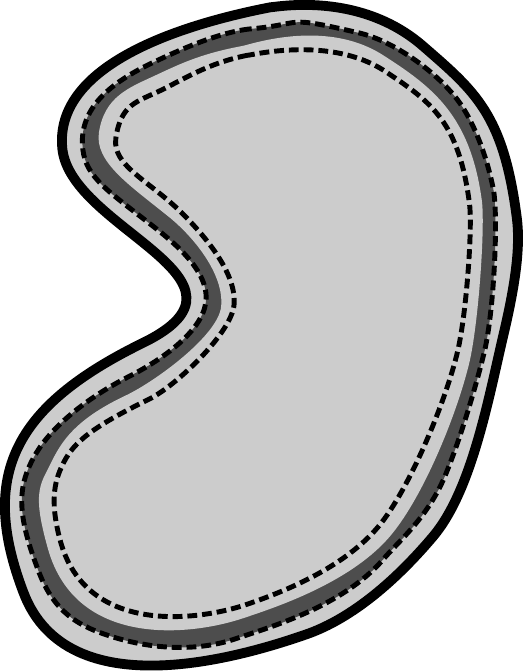}}   
          \caption{The entropy estimate consists of two stages of coupling, indicated by the images above.  The annulus surrounded by a dashed line is $D(R_D^1, R_D^2)$, and inner light grey region in (c) is $E = D(R_D)$.}
\end{figure}

In this subsection we are going to show that we can construct a coupling of $\p_D^\psi, \p_D^{\wt{\psi}}$ so that with high probability the error terms from Proposition \ref{harm::prop::ent_est} are with high probability negligible when applied to $\p_E^\zeta, \p_E^{\wt{\zeta}}$ where $E = D(R_D)$, $R_D = cR^{1-\epsilon_D}$, some $\epsilon_D > 0$, and $(\zeta,\wt{\zeta}) = (h^\psi,h^{\wt{\psi}})|_{(\partial E)^2}.$  We will accomplish this using the following steps:

\begin{enumerate}
\item Take the stationary coupling $(h^\psi,h^{\wt{\psi}})$ of $\p_D^\psi, \p_D^{\wt{\psi}}$.
\item Invoke Lemma \ref{gl::lem::grad_error} to find an annulus $D(R_D^1, R_D^2)$ on which the Dirichlet energy of $\ol{h}$ is controlled.
\item Use Lemma \ref{harm::lem::entropy_form} to show that we can recouple the laws on $D(R_D^1, R_D^2)$ so that with high probability $\ol{h}$ is $\Delta^\beta$-harmonic.  On this event we will have all of the regularity that we need on $\partial D(R_D)$, where $R_D^1 < R_D < R_D^2$.
\end{enumerate}

Fix $\epsilon_D > 0$ so small that $ \epsilon_D  < (10^{-100}\rho_{\rm CD} \wedge \rho_{\rm B} \wedge \xi_{{\rm NC}})^2$, where $\xi_{{\rm NC}}$ is the exponent from the Nash continuity estimate (Lemma \ref{symm_rw::lem::nash_continuity_bounded}), $\rho_{\rm B}$ is the exponent from Lemma \ref{symm_rw::lem::beurling}, and $\rho_{\rm CD}$ is from Theorem \ref{cd::thm::cd}.  We assume $R^{1-\epsilon_D} \leq R_D^1 \leq 2R^{1-\epsilon_D}$ has been chosen such that with $R_D^2 = R_D^1 + R^{1-5\epsilon_D}$ and $R_D = R_D^1 + \tfrac{1}{2} R^{1-5\epsilon_D}$:
\begin{align}
\label{harm::eqn::h_reg_assump_annulus}
  &\E^{\psi,\wt{\psi}} \sum_{b \in D^*(R_D^1, R_D^2)} |\nabla \ol{h}(b)|^2 = O_{\ol{\Lambda}}(R^{-\epsilon_D}),\\
  \label{harm::eqn::h_reg_assump_interior}
  &\E^{\psi,\wt{\psi}} \sum_{b \in D^*(R_D)} |\nabla \ol{h}(b)|^2 = O_{\ol{\Lambda}}(R^{3\epsilon_D}).
\end{align}
That such a choice is possible is ensured by Lemma \ref{gl::lem::grad_error}.  Let $E = D(R_D)$.

In order to apply Proposition \ref{harm::prop::ent_est} we need to make sure that we can arrange for the number of balls required to cover annuli near the boundary is not too large.  Such estimates would come for free if $D$ was a lattice approximation of a smooth domain in $\R^2$.  As we want Theorem \ref{harm::thm::coupling} to hold for general bounded subsets of $\Z^2$, such estimates do not necessarily hold uniformly but only on the average provided we are far enough from $\partial D$.  Let $\gamma_k = k \epsilon_D$ and let $N_k$ be as in Proposition \ref{harm::prop::ent_est}. In particular, by using the averaging technique of Lemma \ref{gl::lem::grad_error} we can arrange for $R_D$ to be such that
\begin{align}
  N_{k} &= O(R^{1+\gamma_2 - \gamma_k}),\ \ 
  |E(R^{\gamma_k})| = O(R^{1+\gamma_{k+2}}),\ \ 
  |\partial E| = O(R^{1+\gamma_1}) \label{harm::eqn::boundary_growth}
 \end{align}
for all $k \leq M$, $M$ the largest integer such that $\gamma_M < 1$.

\begin{lemma}[Harmonic Coupling at the Boundary]
\label{harm::lem::harmonic_coupling}  There exists a coupling $(h^{\psi},h^{\wt{\psi}})$ of $\p_D^\psi, \p_D^{\wt{\psi}}$ such that
\[ \CH = \{ \ol{h} = h^{\psi} - h^{\wt{\psi}} \text{ is $\Delta^\beta$-harmonic in }D(R_D^1, R_D^2)\}\]
occurs with probability $1- O_{\ol{\Lambda}}(R^{-\epsilon_D/2})$.
\end{lemma}
\begin{proof}
Let $\ol{h} = h^{\psi} - h^{\wt{\psi}}$, $F = D(R_D^1, R_D^2)$, and let $g \colon F \to \R$ be $\Delta^\beta$-harmonic in $F$ with boundary values $\ol{h}$ on $\partial F$. By our choice of $R_D$,
\[ \sum_{b \in F^*} \E^{\psi,\wt{\psi}} (\nabla g(b))^2 = O_{\ol{\Lambda}}(R^{-\epsilon_D})\]
as $g$ is harmonic in $F$, has the same boundary values as $\ol{h}$, and $\ol{h}$ satisfies the same estimate.
 Let $\zeta, \wt{\zeta} = (h^{\psi},h^{\wt{\psi}})|_{ \partial F \times \partial F}$.  Conditional on $(\zeta,\wt{\zeta})$, let $\p_F^\zeta, \p_F^{\wt{\zeta}}$ have the laws of the GL model on $F$ with boundary conditions $\zeta,\wt{\zeta}$, respectively, and let $\Q_F^{\zeta,g}$ have the law of $h^{\zeta} - g$ where $h^{\zeta} \sim \p_F^{\zeta}$.  It follows from the Cauchy-Schwarz inequality and Lemma \ref{harm::lem::entropy_form} that
\[ \E^{\psi,\wt{\psi}}[ \h(\p_F^{\wt{\zeta}}|\Q_F^{\zeta,g}) + \h(\Q_F^{\zeta,g}|\p_F^{\wt{\zeta}})] = O_{\ol{\Lambda}}(R^{-\epsilon_D}).\]
The lemma follows by invoking Pinsker's inequality \eqref{harm::eqn::pinsker}.
\end{proof}

\subsubsection{Regularity Estimate}
In the following lemma we will use $(h^{\psi},h^{\wt{\psi}})$ to indicate a random variable with joint law given by the coupling of $\p_D^\psi,\p_D^{\wt{\psi}}$ from Lemma \ref{harm::lem::harmonic_coupling} and $\CH$ the corresponding event.  Let $(\zeta,\wt{\zeta}) = (h^\psi,h^{\wt{\psi}}) |_{ \partial E \times \partial E}$.  Let $g \colon E \to \R$ be the $\Delta^\beta$-harmonic extension of $\ol{\zeta} = \zeta - \wt{\zeta}$ from $\partial E$ to $E$.    Recall the definition of $\|\ol{\zeta}\|_\nabla^E$ just before the statement of Proposition \ref{harm::prop::ent_est}.

\begin{lemma}
\label{harm::lem::harmonic_gradient}
There exists $1 \leq c_D \leq 10$ so that
\begin{align}
  \left(\E^{\psi,\wt{\psi}}(\| \ol{\zeta}\|_{\nabla}^E)^p \one_\CH\right)^{1/p} &= O_{\ol{\Lambda},p}(R^{c_D \epsilon_D-1}), \label{harm::eqn::boundary_reg}\\
  \left( \E^{\psi,\wt{\psi}}[ \max_{b \in E^*} |\nabla g(b)|^p \one_{\CH}] \right)^{1/p} &= O_{\ol{\Lambda},p}(R^{c_D \epsilon_D-1}) \label{harm::eqn::harm_reg}
\end{align}
for every $p \geq 1$.
\end{lemma}
\begin{proof}
By construction, $E$ has the $r$-exterior ball property for $r = R^{1-\epsilon_D}$.  Furthermore, if $x,y \in \partial E$ with $|x-y| \leq \tfrac{1}{8}R^{1-5\epsilon_D}$ then the shortest path connecting $x$ to $y$ in $\Z^2$ is contained in $B(x,\tfrac{1}{4} R^{1-5\epsilon_D})$.  Consequently, 
\begin{equation}
\label{harm::eqn::boundary_reg1}
|\ol{\zeta}(x) - \ol{\zeta}(y)| \leq \breve{g}|x-y|
\end{equation}
where 
\[ \breve{g} = \max\{ |\nabla \ol{h}(b)| : b \in D^*(R_D^1 + \tfrac{1}{4} R^{1-5\epsilon_D}, R_D^1 + \tfrac{3}{4} R^{1-5\epsilon_D})\}.\]
Let $M = \max\{ |\ol{h}(x)| : x \in D\}$.  Since $\ol{h}$ is harmonic in $D(R_D^1,R_D^2)$ on $\CH$, we have that 
\begin{equation}
\label{harm::eqn::g_bound}
 \breve{g} \leq \ol{g} \equiv \frac{C M}{R^{1-5\epsilon_D}}.
\end{equation}
We may assume without loss of generality that $C \geq 100$.
If $|x-y| \geq \tfrac{1}{8}R^{1-5\epsilon_D}$ then 
\[ \ol{g} |x-y| \geq \tfrac{1}{8} CM \geq 2M \geq |\ol{\zeta}(x) - \ol{\zeta}(y)|.\]
Therefore Lemma \ref{dhf::lem::exterior_regularity} implies, by possibly increasing $C > 0$, that
\begin{equation}
\label{harm::eqn::max_gradient}
 \max_{b \in E^*} |\nabla g(b)| \leq \frac{CM}{R^{1-5\epsilon_D}} \left[  \log R +  \frac{R}{R^{1-\epsilon_D}} \right] \leq \frac{C M}{R^{1-6\epsilon_D}}
\end{equation}
on $\CH$.  Trivially,
\[ M \leq \max_{x \in D} |h^\psi(x)| + \max_{x \in D} |h^{\wt{\psi}}(x)|\]
and by Lemma \ref{gl::lem::mom_max} we know that $(\E M^p)^{1/p} = O_{\ol{\Lambda},p}(R^{\epsilon})$.  This clearly gives \eqref{harm::eqn::harm_reg}.  Combining \eqref{harm::eqn::boundary_reg1} with \eqref{harm::eqn::g_bound} gives \eqref{harm::eqn::boundary_reg}.
\end{proof}

\subsubsection{Putting Everything Together}

\begin{proof}[Proof of Theorem \ref{harm::thm::coupling}]
To prove the theorem we just have to estimate the error terms from Proposition \ref{harm::prop::ent_est} on $\CH$.  First of all, by Lemmas \ref{gl::lem::mom_max}, \ref{harm::lem::harmonic_gradient} we observe
\begin{equation}
 (\E[ (\CE)^p \one_\CH])^{1/p} = O_{\ol{\Lambda},p}(R^{1+2\epsilon_D - a_1\rho_{\rm B}})
\end{equation}
for $a_1 = 1/10$.  Consequently,
\begin{align}
   \E[\CE_E \one_{\CH}] =& O_{\ol{\Lambda}}(R^2 \cdot R^{3c_D \epsilon_D - 3} + R^{1+2\epsilon_D - a_1\rho_{\rm B} + c_D \epsilon_D - 1}) \notag\\
     =& O_{\ol{\Lambda}}(R^{c_1\epsilon_D - a_1\rho_{\rm B}}) \label{harm::eqn::exp1}
 \end{align}
 for $c_1 < 100$.  Using Lemmas \ref{gl::lem::mom_max}, \ref{harm::lem::harmonic_gradient} again, we see that
 \begin{align}
   \E[\CE_B^\ell \one_{\CH}] =& O_{\ol{\Lambda}}(R^{1+2\epsilon_D + \gamma_{\ell+2} - 1}) O_{\ol{\Lambda}}(R^{1/2+c_D\epsilon_D - 1} + R^{\epsilon_D+(\gamma_\ell - 1/2) \rho_{\rm B}}) \notag\\
   =& 
   O_{\ol{\Lambda}}(R^{c_1\epsilon_D + 2\gamma_{\ell+2} - a_1\rho_{\rm B}}) \label{harm::eqn::exp2},
\end{align}
the last line coming as $\rho_{\rm B} \leq 1$.
Finally, Lemma \ref{harm::lem::harmonic_gradient} clearly implies
\begin{align}
   \E[ \CE_I^\ell \one_\CH] &= \sum_{k=\ell+1}^M O_{\ol{\Lambda}}(R^{1+\gamma_2-\gamma_k} R^{\epsilon_D + \gamma_k(1-\rho_{\rm CD})} R^{c_2 \epsilon_D - 1}) \notag \\
    &= \sum_{k=\ell+1}^M O_{\ol{\Lambda}}(R^{c_3 \epsilon_D - \gamma_\ell \rho_{\rm CD}}) \label{harm::eqn::exp3}
\end{align}
for $c_3 \leq 1000$.  The exponents from \eqref{harm::eqn::exp1}, \eqref{harm::eqn::exp2}, \eqref{harm::eqn::exp3} are
\[ c_1 \epsilon_D - a_1 \rho_{\rm B},\ \ c_1 \epsilon_D + 2 \gamma_{\ell+2} - a_1 \rho_{\rm B},\ \ c_3 \epsilon_D -  \gamma_{\ell} \rho_{\rm CD}.\]
Choosing $\ell > 10^5$ we see that the second and third exponents are negative and the first is clearly negative.
\end{proof}

We finish this section with the short proof of Theorem \ref{harm::thm::mean_harmonic}.

\begin{proof}[Proof of Theorem \ref{harm::thm::mean_harmonic}]
Assume that we still have the setup of the previous theorem except $\wt{\psi} = 0$.  Then there exists $\delta > 0$ so that we can find a coupling of $\p_D^\psi, \p_D^{\wt{\psi}}$ so that $\ol{h}$ is harmonic in $E$ on an event $\wt{\CH}$ (this is different from $\CH$ in the proof of Theorem \ref{harm::thm::coupling}) with probability $1-O_{\ol{\Lambda}}(R^{-\delta})$. Let $g,\wh{h}$ be the harmonic extensions of $\ol{h}$, $\E h^\psi(x)$ from $\partial E$ to $E$, respectively.  For $x \in E$ we have that
\begin{align*}
   \E h^\psi(x) &= \E \ol{h}(x) = \E g(x) (1-\one_{\wt{\CH}^c}) + \E \ol{h}(x) \one_{\wt{\CH}^c}\\
      &= \wh{h}(x) + \E (\ol{h}(x) - g(x)) \one_{\wt{\CH}^c}
\end{align*}
since $\E h^{\wt{\psi}}(x) = 0$.
Since $\psi \in \B_{\ol{\Lambda}}^u(D)$, Lemma \ref{gl::lem::mom_max} implies both $\E \ol{h}^2(x) = O_{\ol{\Lambda}}(R^{\delta/2})$ and $\E g^2(x) = O_{\ol{\Lambda}}(R^{\delta/2})$.  Consequently, an application of Cauchy-Schwarz yields
\[ \E\big[ |\ol{h}(x)| + |g(x)| \big] \one_{\wt{\CH}^c} = O_{\ol{\Lambda}}(R^{-\delta/4}),\]
from which the theorem follows.
\end{proof}

\section{The Central Limit Theorem}
\label{sec::clt}

We will prove Theorem \ref{thm::clt} in this section.  The primary inputs are Theorem \ref{harm::thm::coupling} and the main result of either \cite{GOS01} or \cite{NS97}.  Throughout, we let $D \subseteq \R^2$ be a connected, bounded, smooth domain and for each $n$ let $D_n = D \cap \tfrac{1}{n} \Z^2$.  We will be dealing with both discrete and continuum derivatives, so to keep the notation consistent with the rest of the article we will still let $\nabla, \Delta$ denote the discrete gradient and Laplacian, respectively, and use $\ol{\nabla}$ and $\ol{\Delta}$ for their continuum counterpart.  We begin with a simple analysis lemma.  Let $\beta = \beta(u)$ be as in the statement of Theorem \ref{thm::clt}.

\begin{lemma}
\label{clt::lem::discrete_de}
For each $\epsilon > 0$ there exists $c > 0$ so that for all $g \in C^\infty(D)$ we have that
\begin{align*}
   \sum_{b \in D_n^*} |\nabla g(b)|^2 &\leq c \sup_{x \in D} \|\ol{\nabla} g(x)\|^2,\ \ 
   \sum_{b \in D_n^*} |\nabla g(b)|^2 \leq c \| g \|_{H^{2+\epsilon}(D)}^2,\\
   \sum_{b \in \partial D_n^*} |\nabla g(b)| &\leq c \| g \|_{H^{2+\epsilon}(D)},\ \ 
   \sum_{x \in D_n} | \Delta^\beta g(x)| \leq c \| g\|_{H^{3+\epsilon}(D)}.
\end{align*}
\end{lemma}
\begin{proof}
The first claim is obvious since $|\nabla g(b)| \leq \tfrac{1}{n} \sup_{x \in D} \|\ol{\nabla} g(x)\|$.
As for the second claim, we note that the Sobolev embedding theorem implies that for each $\epsilon > 0$ there exists $c,c' > 0$ so that 
\[ \sup_{x \in D} \|\ol{\nabla} g(x) \| \leq c' \| \ol{\nabla} g \|_{H^{1+\epsilon}(D)} \leq c \| g\|_{H^{2+\epsilon}(D)}.\]
See, for example, from Proposition 1.3 in Chapter 4 of \cite{TAY96}.  The final two claims are proved similarly.
\end{proof}

Let $\varphi_n(x) = nu \cdot x$.  Fix a continuous function $f \colon \R^2 \to \R$ and let $h^n$ have the law of the GL model on $D_n$ with $h^n(x) = f(x) + \varphi_n(x)$ for $x \in \partial D_n$, $\eta^{n,D} = \nabla h^n$, and for $g \in H^{3+\epsilon}(D)$ define
\begin{align*}
     \xi_\nabla^{n,D}(g) &= \sum_{b \in D_n^*} a_u(b) \nabla g(b) ( \eta^{n,D}(b) - \nabla \varphi_n(b))
\end{align*}
where $a_u(b) = \E[ \CV''(\eta(b))]$ for $\eta \sim \SEGGS_u$.  Note that $a_u(b)$ depends only on the orientation of $b$.

\begin{lemma}
\label{clt::lem::exp_moments}
For each $\epsilon > 0$ there exists $c > 0$ so that for all $g \in H^{3+\epsilon}(D)$ we have
\[ \E \exp( \xi_\nabla^{n,D}(g)) \leq \exp\left[ c \left( \|f\|_\infty^2 + c \| g \|_{H^{3+\epsilon}(D)}^2\right) \right].\]
\end{lemma}
\begin{proof}
Using summation by parts, we have
\begin{align*}
      |\E \xi_\nabla^{n,D}(g)|
 \leq& \sum_{x \in D_n} |\E (h^n(x) - \varphi_n(x)) \Delta^\beta g(x)| + c_1\sum_{b \in \partial D_n^*} |f(x_b) \nabla g(b)|\\
  \leq& c_2 \| f \|_{\infty} \| g \|_{H^{3+\epsilon}(D)}. 
\end{align*}
In the final inequality we used that $|\E(h^n(x) - \varphi_n(x))| \leq \|f \|_\infty$, which is a consequence of Lemma \ref{gl::lem::hs_mean_cov}, in addition to Lemma \ref{clt::lem::discrete_de}.
The exponential Brascamp-Lieb inequality (Lemma \ref{bl::lem::bl_inequalities}) combined with \eqref{dgff::covariance} and the previous lemma implies there exists $c_3,c_4 > 0$ so that
\begin{align*}
           \E \exp( \xi_\nabla^{n,D}(g))
&\leq  \exp\left[ c_3 \left(\|f\|_\infty \| g \|_{H^{3+\epsilon}(D)} + \sum_{b \in D_n^*} (\nabla g(b))^2 \right) \right]\\
 &\leq  \exp\left[ c_4 \left(\|f\|_\infty^2 +\| g \|_{H^{3+\epsilon}(D)}^2 \right) \right].
\end{align*}
\end{proof}

\begin{lemma}
For each $\kappa > 4$, the law of $\xi_\nabla^{n,D}$ induces a tight sequence on $H^{-\kappa}(D)$ equipped with the weak topology.
\end{lemma}
\begin{proof}
Fix $\kappa > 4$.  It suffices to show that for each $\delta > 0$ there exists $M = M(\delta)$ such that
\[ \p[ \| \xi_\nabla^{n,D} \|_{H^{-\kappa}(D)} \geq M] \leq \delta\]
as the Banach-Alaoglu theorem yields that the ball $\{ \| g \|_{H^{-\kappa}(D)} \leq M\}$ is compact in the weak topology of $H^{-\kappa}(D)$.

Let $\epsilon = \tfrac{1}{2}(\kappa-4)$ and $\kappa' = \kappa-1-\epsilon > 3$.  Let $(\wt{g}_k)$ be the eigenvectors of $\ol{\Delta}$ on $D$, normalized to be orthonormal in $L^2(D)$, with negative eigenvalues $(\lambda_k)$ ordered to be non-increasing in $k$.  Let $g_k = (1-\ol{\Delta})^{-\kappa/2} \wt{g}_k$.  By \eqref{gff::eqn::sobolev_inner_product}, $(g_k)$ is an orthonormal basis of $H^{\kappa}(D)$.  As $g_k / (1-\lambda_k)^{(1+\epsilon)/2} = (1-\ol{\Delta})^{-(1+\epsilon)/2} g_k$ we have that
\[ \| g_k \|_{H^{\kappa'}(D)} = \| (1-\ol{\Delta})^{-(1+\epsilon)/2} g_k \|_{H^{\kappa}(D)} = \frac{1}{(1-\lambda_k)^{(1+\epsilon)/2}}.\]
The Weyl formula implies that $k/(-\lambda_k)$ tends to a constant $c = c_D$ depending only on $D$ as $k \to \infty$.  Therefore there exists $c_D' \geq c_D$ so that
\[ \| k^{(1+\epsilon)/2} g_k \|_{H^{\kappa'}(D)} \leq c_D' \text{ for all } k \in \N.\]
Combining the above with Chebychev's inequality yields
\[ \p[ |\xi_\nabla^{n,D}(g_k)| \geq M/k^{(1+\epsilon/2)/2} ] = \p[ |\xi_\nabla^{n,D}(k^{(1+\epsilon)/2} g_k)| \geq M k^{\epsilon/4}] \leq \exp(c - M k^{\epsilon/4} )\]
where $c = c(\epsilon,D,f)$.
Consequently, letting $A_M^n = \cap_k \{ |\xi_\nabla^{n,D}(g_k)| \leq M/k^{(1+\epsilon/2)/2}\}$, a union bound yields $\p[ A_M^n] \to 1$ as $M \to \infty$.
Note that if $g \in H^{\kappa}(D)$, $g = \sum_k \alpha_k g_k$ for $(\alpha_k) \in \ell^2$, we have 
\[ \xi_\nabla^{n,D}(g) = \sum_k \alpha_k \xi_\nabla^{n,D}(g_k)\]
since $H^{\kappa}(D)$-convergence implies uniform convergence as $\kappa > 4$ and $\xi_\nabla^{n,D}$ is obviously continuous in the uniform topology on functions $D_n \to \R$.
Let 
\[ N_\epsilon = \sum_{k} \frac{1}{k^{1+\epsilon/2}}.\] 
On $A_M^n$, note that
\[ |\xi_{\nabla}^{n,D}(g)| \leq \left(\sum_k \alpha_k^2\right)^{1/2} \left( \sum_k (\xi_\nabla^{n,D}(g_k))^2\right)^{1/2} \leq \| g\|_{H^{\kappa}(D)} M \sqrt{N_\epsilon}.\]
Consequently,
\begin{align*}
   \p\bigg[ \sup_{ \|g \|_{H^{\kappa}(D) }\leq 1} |\xi_\nabla^{n,D}(g)| \geq M \sqrt{N_\epsilon} \bigg]
&\leq \p[(A_M^n)^c].
\end{align*}
This proves for every $\delta > 0$ there exists $\wt{M} = \wt{M}(\delta)$ sufficiently large so that 
\[ \p[ \| \xi_\nabla^{n,D}\|_{H^{-\kappa}(D)} \geq \wt{M}] \leq \delta\]
for every $n \in \N$.
\end{proof}

Let $\eta \sim \SEGGS_u$, but thought of as a random gradient field on $(\tfrac{1}{n} \Z^2)^*$.  Fix a base point $x^* \in \partial D$ and let $x_n$ be a point in $\partial D_n$ with minimal distance to $x^*$.  Set $h^{n,0}(x_n) = 0$ and let $h^{n,0}$ be the function satisfying $\nabla h^{n,0} = \eta$.
Let
\begin{align*}
     \xi_\nabla^{n}(g) &= \sum_{b} a_u(b) \nabla g(b)(\eta^{n}(b) - \nabla \varphi_n(b)) \text{ for } g \in C_0^\infty(\R^2).
\end{align*}
Corollary 2.2 of \cite{GOS01} implies that for any $g_1,\ldots,g_k \in C_0^\infty(\R^2)$ fixed, the random vector $(\xi_\nabla^n(g_1),\ldots, \xi_\nabla^n(g_k))$ converges in distribution to a zero-mean Gaussian vector $(Z_1,\ldots,Z_k)$ with covariance $\cov(Z_i, Z_j) = (g_i,g_j)_\nabla^A$ for $A = A(u,\CV)$ depending only on the tilt $u$ and $\CV$.  Note that our definition of $\xi_\nabla^n$ differs from $\xi^\epsilon$ in \cite{GOS01} in that we do not have a normalization of $n^{-1}$.  The reason is that $\xi_\nabla^n$ operates on discrete gradients of $C_0^\infty(\R^2)$ functions, which themselves are of order $n^{-1}$.

\begin{lemma}
\label{clt::lem::bc_harm}
There exists $\ol{\Lambda}$ depending only on $\CV$ such that
\[ \p[ h^{n,0} |_{\partial D_n} \in \B_{\ol{\Lambda}}^u(D_n)] = 1- O(n^{-8}).\]
\end{lemma}
\begin{proof}
Let $x \in D_n$ and $\ol{x} = x-x_n$.  Combining the exponential Brascamp-Lieb inequality with \eqref{eqn::var_bound} yields for $x \in \partial D_n$ that 
\[ \E \exp( h^{n,0}(x) - \varphi_n(\ol{x})) = \E \exp( h^{n,0}(x) - h^{n,0}(x_n) - \varphi_n( \ol{x})) \leq \exp( C \log n) = n^C.\]
Assume without loss of generality that $C \geq 1$.  By Chebychev's inequality,
\[ \p[ |h^{n,0}(x) - \varphi_n(\ol{x})| \geq 10 C \log n] \leq n^{-9}.\]
Using a union bound we thus have
\[ \p[ \max_{x \in \partial D_n} |h^{n,0}(x) - \varphi_n(\ol{x})| \geq 10 C \log n] \leq O(n^{-8}).\]
Consequently, taking $\ol{\Lambda} = 10C$ we have that $\p[ h^{n,0}|_{\partial D_n} \in \B_{\ol{\Lambda}}^u(D_n)] = 1-O(n^{-8})$.
\end{proof}

By the same proof as Lemma \ref{clt::lem::exp_moments} we have
\begin{equation}
\label{clt::eqn::seggs_exp_moment}
 \E \exp( \xi_\nabla^n(g)) \leq \exp(c \| g\|_{H^{3+\epsilon}(D)}^2).
\end{equation}

\begin{proof}[Proof of Theorem \ref{thm::clt}]
Fix $\kappa > 4$ and $f \colon \R^2 \to \R$ continuous.  Let $h^0,h^f$ be GFFs on $D$ with respect to $(\cdot,\cdot)_\nabla^A$, $A = A(u,\CV)$ as before, where $h^0$ has zero boundary conditions and the boundary condition of $h^f$ is given by $f|_{\partial D}$.  Let $\xi_\nabla^D$ be a weak-$H^{-\kappa}(D)$ subsequential limit of $(\xi_\nabla^{n,D})$.  We will prove for any $g_1,\ldots,g_k \in C^\infty(D)$ that 
\[ (\xi_\nabla^D(g_1),\ldots,\xi_\nabla^D(g_k)) \stackrel{d}{=} ( (h^f,g_1)_\nabla^A,\ldots,(h^f,g_k)_\nabla^A)\]
since the continuity of $\xi_\nabla^D$ implies that its law is determined by its projections onto $C^\infty(D)$, a dense subset of $H^\kappa(D)$.  We will identify $A$ at the end of the proof.  To establish this, it suffices to show that 
\[ (\xi_\nabla^{n,D}(g_1),\ldots,\xi_\nabla^{n,D}(g_k)) \stackrel{d}{\to} ((h^f,g_1)_\nabla^A,\ldots,(h^f,g_k)_\nabla^A).\]
We will first prove the result for $C_0^\infty(D)$, then using an approximation argument generalize to $C^\infty(D)$.

By Lemma \ref{clt::lem::bc_harm}, with probability $1-O(n^{-8})$ we can apply Theorem \ref{harm::thm::coupling} to $\p_{D_n}^\psi, \p_{D_n}^{\wt{\psi}}$ where $\psi = f + \varphi_n$ and $\wt{\psi} = h^{n,0}|_{\partial D_n}$. This implies the existence of $\epsilon,\delta > 0$ independent of $n$ such that we can couple $h^n, h^{n,0}$ so that with $\wh{h}^n$ the $\Delta^\beta$-harmonic extension of $h^n - h^{n,0}$ from $\partial D_n(n^{-\epsilon})$ to $D_n(n^{-\epsilon})$, we have
\[ \p[ \CH^c ] = O(n^{-\delta}) \text{ where } \CH = \{\ol{h}^n = \wh{h}^n \text{ in } D_n(n^{-\epsilon})\}\]
for all $n$ large enough.  The reason that we see $n^{-\epsilon}$ rather than $n^{1-\epsilon}$ as in the statement of Theorem \ref{harm::thm::coupling} is that $D_n = D \cap \tfrac{1}{n} \Z^2$, so all of our distances need to be scaled by a factor of $n^{-1}$.
Fix $g_1,\ldots,g_k \in C_0^\infty(D)$ and assume that $n$ is sufficiently large so that $\supp(g_1),\ldots,\supp(g_k) \subseteq D_n(n^{-\epsilon})$.  On $\CH$, for each $1 \leq i \leq k$ we have that
\begin{equation}
\label{clt::eqn::sum_by_parts}
 \xi_\nabla^{n,D}(g_i) = \xi_\nabla^{n}(g_i) + \sum_{b \in D_n^*} a_u(b) \nabla g_i(b) \nabla \ol{h}^n(b) = \xi_\nabla^{n}(g_i).
\end{equation}
The second equality follows from summation by parts and the $\Delta^\beta$-harmonicity of $\wh{h}^n$; there is no boundary term since $g_i$ vanishes near $\partial D_n(n^{-\epsilon})$.  

Combining Lemma \ref{clt::lem::exp_moments}, \eqref{clt::eqn::seggs_exp_moment}, and the  Cauchy-Schwarz inequality yields
\begin{align*}
 &  \E |\xi_\nabla^{n,D}(g_i) - \xi_\nabla^{n}(g_i)|
  = \E |\xi_\nabla^{n,D}(g_i) - \xi_\nabla^n(g_i)| \one_{\CH^c}\\
\leq& \big[ O(1) O(n^{-\delta}) \big]^{1/2}
 \leq O(n^{-\delta/2}).
\end{align*}
Therefore $(\xi_\nabla^D(g_1),\ldots,\xi_\nabla^D(g_k))$ is a Gaussian vector with $\cov( \xi_\nabla^D(g_i), \xi_\nabla^D(g_j)) = \cov( (h,g_i)_\nabla^A, (h,g_j)_\nabla^A)$ where $h$ is an $A$-GFF on $\R^2$.  Proposition \ref{gff::prop::markov} implies that $h$ restricted to $D$ admits the decomposition $h = h^0 + \wh{h}$ where $h^0$ is a zero-boundary $A$-GFF on $D$ and $\wh{h}$ is a $\ol{\Delta}^A$-harmonic function.  Integration by parts implies that $(\wh{h},g_i)_\nabla^A \equiv 0$ for all $i$, consequently the covariance structure of $(\xi_\nabla^D(g_1),\ldots,\xi_\nabla^D(g_k))$ is the same as $( (h^0,g_1)_\nabla^A,\ldots,(h^0,g_k)_\nabla^A)$.

We now turn to the general case that $g_1,\ldots,g_k \in C^\infty(D)$ do not necessarily have compact support in $D$.  Note that we can write
\[  g_i = (g_i - \wt{g}_i - \wh{g}_i) + \wt{g}_i + \wh{g}\]
where $\wh{g}_i$ is the $\ol{\Delta}^\beta$-harmonic extension of $g_i$ from $\partial D$ to $D$ and $\wt{g}_i \in C_0^\infty$ satisfies $\|g_i - \wh{g}_i - \wt{g}_i(x) \|_{H^1(D)} \leq \delta_1$.  Note that such an approximation exists since $g_i - \wh{g}_i \in H_0^1(D)$.
Since $\wh{g}_i$ is harmonic with smooth boundary conditions, we have $\Delta^\beta \wh{g}_i = o(1) n^{-2}$ uniformly in $D$.  Thus summing by parts twice and using that $h^n(x) = f(x) + \varphi_n(x)$ on $\partial D_n$, with $f^n$ the $\Delta^\beta$-harmonic extension of $f$ from $\partial D_n$ to $D_n$ we have
\begin{align*}
   \E \xi_\nabla^{n,D}(\wh{g}_i) = \sum_{b \in \partial D_n^*} a_u(b) f(x_b) \nabla \wh{g}_i(b) + o(1) =
   \sum_{b \in D_n^*} a_u(b) \nabla f^n(b) \nabla \wh{g}_i(b) + o(1).
\end{align*}
Thus it is not difficult to see that if $F$ denotes the $\ol{\Delta}^\beta$-harmonic extension of $f$ from $\partial D$ to $D$  then
\[ \lim_{n \to \infty} \E \xi_\nabla^{n,D}(\wh{g}_i) = \int_D \ol{\nabla} F A_u \ol{\nabla} \wh{g}_i = \int_D \ol{\nabla} F A_u \ol{\nabla} g_i\]
where $A_u$ is the diagonal matrix with entries $\beta = (\beta_1,\beta_2)$.

Applying summation by parts,
\[ \xi_{\nabla}^{n,D}(\wh{g}_i) = \sum_{b \in \partial D_n^*} a_u(b) f(x_b) \nabla \wh{g}_i(b) - \sum_{x \in D_n} f(x) \Delta^\beta \wh{g}_i(x).\]
Note that the first summation on the right hand side is deterministic.  Consequently, combining the Brascamp-Lieb inequalities with \eqref{dgff::covariance} implies
\begin{align*}
 &      \var\left( \xi_\nabla^{n,D}(\wh{g}_i) \right)
   = O(1) \sum_{x,y \in D_n} |\Delta^\beta \wh{g}_i(x) \Delta^\beta \wh{g}_i(y)| G_n(x,y)\\
 &= o(1) \sum_{x \in D_n} n^{-2} \sum_{y \in D_n} n^{-2} G_n(x,y)
  = o(1)
\end{align*}
where $G_n$ is the discrete Green's function on $D_n$.  This takes care of $\wh{g}_i$.   We already know that the limiting behavior of $\xi_\nabla^{n,D}(\wt{g}_i)$, which leaves us to deal with $g_i-\wh{g}_i-\wt{g}_i$.  Invoking Lemma \ref{clt::lem::discrete_de} and the Brascamp-Lieb inequality, we have
\begin{align}
\limsup_{n \to \infty} \E ( \xi_\nabla^{n,D} (g_i-\wh{g}_i-\wt{g}_i))^2 &\leq 
 C \limsup_{n \to \infty} \sum_{b \in D_n^*} |\nabla (g_i - \wh{g}_i - \wt{g}_i)(b)|^2 \notag\\
&= C \| g_i - \wh{g}_i - \wt{g}_i \|_{H^1(D)}^2 \leq c \delta_1^2 \label{clt::eqn::approx_bound}. 
\end{align}
In the equality, we are using that $g_i - \wh{g}_i - \wt{g}_i \in C^\infty(D)$.  We also have
\[ \E |(h^0, g_i-\wh{g}_i - \wt{g}_i)_\nabla^A|^2 \leq  c\| g_i - \wh{g}_i - \wt{g}_i \|_{H^1(D)}^2 \leq c\delta_1^2.\]

Assume that $(\xi_\nabla^{n,D}(\wt{g}_1),\ldots,\xi_\nabla^{n,D}(\wt{g}_k))$ and $((h^0,\wt{g}_1)_\nabla^A,\ldots,(h^0,\wt{g}_k)_\nabla^A)$ have been embedded into a common probability space so that $\lim_n \xi_\nabla^{n,D}(\wt{g}_i) = (h^0,\wt{g}_i)_\nabla^A$ almost surely for each $1 \leq i \leq k$.  By \eqref{clt::eqn::approx_bound}, 
\begin{align*}
 &\E \limsup_{n \to \infty} |\xi_\nabla^{n,D}(g_i) - (h^0,g_i)_\nabla^A - \int_D \ol{\nabla} F A_u \ol{\nabla} g_i|\\
=& \E \limsup_{n \to \infty} |\xi_\nabla^{n,D}(g_i-\wh{g}_i-\wt{g}_i) - (h^0,g_i - \wh{g}_i -\wt{g}_i)|\\  \leq& 2\delta_1.
\end{align*}
Since $\delta_1 > 0$ was arbitrary, we therefore have that $(\xi_\nabla(g_1),\ldots,\xi_\nabla(g_k))$ has the same distribution as 
\[\left((h^0,g_i)_\nabla^A + \int \ol{\nabla} F A_u \ol{\nabla} g_i : 1 \leq i \leq k \right).\]

We will now explain why $A = A_u$, which will complete the proof.  We will not spell out all of the details exactly in order to avoid repetition.  Suppose that $U \subseteq D$ is a smooth open subset, $U_n = U \cap \tfrac{1}{n} \Z^2$, and let $\psi_n$ be the function which is $\Delta^\beta$ harmonic in $U_n$ and is equal to $h^n(x) - \varphi_n(x)$ in $D_n \setminus U_n$. Consider the auxiliary functional
\[ \xi_\nabla^{n,U}(g) = \sum_{b \in D^*} a_u(b) \nabla (h^n - \varphi_n - \psi_n)(b) \nabla g(b).\]
Exactly the same argument implies that $\xi_\nabla^{n,U}$ converges to a zero-boundary $A$-GFF on $U$, say $\xi_\nabla^U$.  Since $\psi_n = \xi_\nabla^{n,D} - \xi_\nabla^{n,U}$, as linear functionals, we also know that $\psi_n$ has a limit, say $\psi = \xi_\nabla^D - \xi_\nabla^U$.  It is not difficult to see that $\psi$ is $\ol{\Delta}^\beta$ harmonic and depends on $\xi_\nabla^D$ only through its values on $D \setminus U$.  Since $\xi_\nabla^D$ is an $A$-GFF on $D$, we know that it admits the decomposition $\xi_\nabla^D = \wt{\xi}_\nabla^U + \wt{\psi}$ where $\wt{\psi}$ is $\ol{\Delta}^A$-harmonic and $\wt{\xi}_\nabla^U$ is a zero-boundary $A$-GFF on $U$ independent of $\xi_\nabla^U$.  Therefore we have that
\[ \xi_\nabla^U = \wt{\xi}_\nabla^U + (\wt{\psi}-\psi).\]
This implies that $\wt{\psi} = \psi$ almost surely since a zero-boundary $A$-GFF plus an independent function does not have the law of a zero-boundary $A$-GFF.  This finishes the proof of the theorem.
\end{proof}

\appendix

\section{Discrete Harmonic Functions}

Suppose that $D \subseteq \Z^2$ is bounded and connected.  We say that $D$ satisfies the $r$-exterior ball property if for each $x \in \partial D$ there exists $y \in \Z^2$ such that $x \in \partial B(y,r)$ and $B(y,r) \cap D = \emptyset$.

\begin{figure}
     \centering
          \includegraphics[width=.70\textwidth]{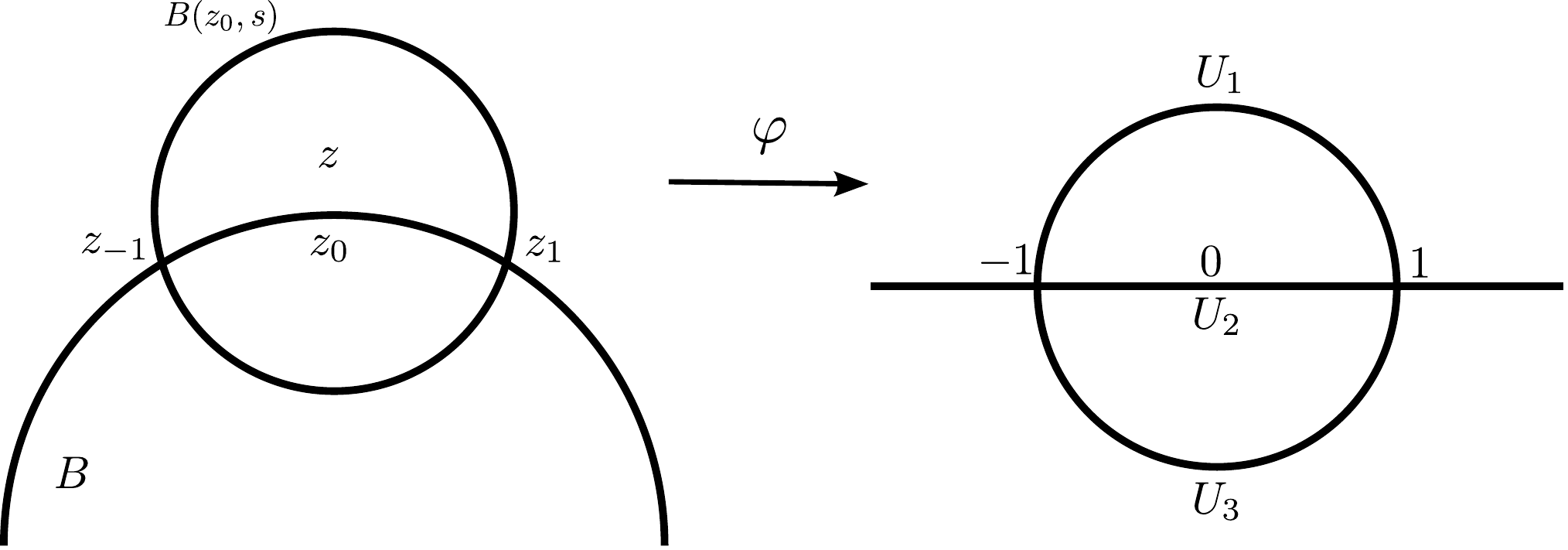}
          \caption{The setup for the proof of Lemma \ref{dhf::lem::beurling}.  Here, we choose a M\"obius transformation $\varphi \colon \C \to \C$ determined by $\varphi(z_{-1}) = -1$, $\varphi(z_0) = 0$, and $\varphi(z_1) = 1$.  In particular, $\varphi(B(z_0,s))$ is mapped to $\D$, the unit disk in $\C$, and $\varphi(B(z_0,s) \setminus B)$ is sent to $\U$, the upper half of the unit disk.}
\end{figure}

\begin{lemma}
\label{dhf::lem::beurling}
Suppose that $D \subseteq \Z^2$ is bounded, connected, and satisfies the $r$-exterior ball property.  Let $\beta_1, \beta_2 > 0$ with $\beta_1+\beta_2 = 1/2$.  Let $X$ be a random walk on $D$ with $X_0 = x$ that jumps up and down with probability $\beta_1$ and left and right with probability $\beta_2$.  Let $\tau_D = \inf\{ t \geq 0 : X_t \notin D\}$, $\tau_s = \inf\{t \geq 0 : |X_t - x| = s\}$, and $d = \dist(x,\partial D)$.  There exists a constant $C = C(\beta_1,\beta_2) > 0$ such that
\[ \p_x[\tau_s \leq \tau_D] \leq C \frac{d}{s \wedge r}.\]
\end{lemma}
\begin{proof}
We are first going to prove a related result for Brownian motion, then explain how to deduce the corresponding result for random walk.  Clearly, we may assume that $d \leq s/4$ and $s \leq r/4$.  Suppose that $z \in \C$, $B$ is a ball of radius $r$, $\dist(z,B) = d$, and $z_0$ is the point in $\partial B$ closest to $z$.  Let $W$ be a Brownian motion initialized from $w \in \C$, $\tau_B^W = \inf\{t \geq 0 : W_t \in B\}$, and $\tau_s^W = \inf\{ t \geq 0 : |W_t - z_0| \geq s\}$.  Finally, let $u(w) = \p_w[ \tau_s^W \leq \tau_B^W]$.  We claim there exists $C > 0$ independent of the setup such that if $w \in B(z_0,s/4)$ then
\[ u(w) \leq \frac{C}{s} \dist(w,B).\]
We know that $u$ solves the Dirichlet problem
\[ \Delta u = 0 \text{ in } B(z_0,s) \setminus B,\ \ u|_{(\partial B(z_0,s)) \setminus B} \equiv 1,\ \  u|_{(\partial B) \setminus B^c(z_0,s)} = 0.\]
Let $z_{-1}, z_{1}$ be the two points in $\partial B \cap \partial B(z_0,s)$ and let $\varphi$ be the M\"obius transformation satisfying $\varphi(z_{-1}) = -1, \varphi(z_0) = 0, \varphi(z_1) = 1$.  Let $\D = \{ z \in \C : |z| \leq 1\}$, $\U = \{ z \in  \D : \im(z) \geq 0\}$, $U_1 = \{z \in \partial \U : \im(z) > 0\}$, and $U_2 = \{ z \in \U : \im(z) = 0\}$.  Then 
\[ \varphi(B(z_0,s) \setminus B) = \U,\ \varphi( (\partial B(z_0,s)) \setminus B) = U_1,\ \varphi( (\partial B \setminus B^c(z_0,s)) = U_2.\]
Consequently, $v = u \circ \varphi^{-1}$ solves
\[ \Delta v = 0 \text{ in } \U,\ \ v|_{U_1} \equiv 1,\ \ v|_{U_2} \equiv 0.\]
Let $U_3 = \{ z \in \partial \D : \im(z) < 0\}$ and let $\wt{v}$ solve
\[ \Delta \wt{v} = 0 \text{ in } \D,\ \ \wt{v}|_{U_1} \equiv 1,\ \ \wt{v}|_{U_3} \equiv -1.\]
By symmetry, $\wt{v}|_{U_2} \equiv 0$, hence $v = \wt{v}$ in $\U$.  Therefore $u$ is the restriction of $\wt{u} = \wt{v} \circ \varphi$.  Since $\wt{u}$ is harmonic in $B(z_0,s)$ with $\| \wt{u}|_{\partial B(z_0,s) }\|_\infty \leq 1$ it follows that $\wt{u}$ is Lipschitz in $B(z_0,s/2)$ with constant $C/s > 0$ and $C$ is independent of the setup.  Fix $w \in B(z_0,s/4)$ and let $w_0$ be the point in $\partial B$ closest to $w$. As $w,w_0 \in B(z_0,s/2)$, we have
\[ |u(w)| = |\wt{u}(w) - \wt{u}(w_0)| \leq \frac{C}{s}|w - w_0| = \frac{C}{s} \dist(w,\partial B).\]

This proves our claim, from which we will now deduce the lemma.  Fix $x \in D$ with $\dist(x,\partial D) = d$.  Since $D$ satisfies the exterior ball property, there exists $y \in \Z^2$ such that with $B = B(y,r)$ we have $B \cap D = \emptyset$ and $\dist(B,x) = d$.  By monotonicity, it suffices to show that $\p_x[ \tau_s \leq \tau_B] \leq Cd / (s \wedge r)$, where $\tau_B = \inf\{ t \geq 0 : X_t \in B\}$.  This is slightly different than the setting for the Brownian motion case because $\tau_s$ is the first time $X$ has distance $s$ from $x$, whereas before we considered the first time $W$ has distance $s$ from $z_0$, which was the point in $B$ closest to $z$.  By the obvious monotonicity of the problem, the desired bound in either case implies the corresponding bound in the other.  The simple random walk estimate obviously follows from the Brownian motion estimate.  The non-simple case follows since the Brownian motion estimate holds even after applying a non-degenerate linear transformation.
\end{proof}

\begin{lemma}
\label{dhf::lem::exterior_regularity}
There exists a constant $C > 0$ such that the following holds.  Suppose that $D \subseteq \Z^2$ is bounded, connected, and satisfies the $r$-exterior ball property.  Suppose further that $g \colon D \to \R$ is a $\Delta^\beta$-harmonic function, $\beta_1, \beta_2 > 0$, such that there exists $\ol{g} \geq 0$ such that $|g(x)-g(y)| \leq k\ol{g}$ if $x,y \in \partial D$ and $|x-y| = k$.  Then
\[ \max_{b \in D^*} |\nabla g(b)| \leq C \ol{g} \left[ \log R +  \frac{R}{r}\right] \]
where $R = \diam(D)$. 
\end{lemma}
\begin{proof}
Fix $b =(x,y) \in D^*$.  Let $X_t$ be a random walk, jumping up and down with rate $\beta_1$ and left and right with rate $\beta_2$, initialized from $x$, $Y_t = X_t + (y-x)$, and $\sigma = \inf\{t \geq 0 : X_t \in \partial D \text{ or } Y_t \in \partial D\}.$  Since $g(X_t) - g(Y_t)$ is a bounded martingale, the optional stopping theorem implies that
\[ |\nabla g(b)| \leq \E |g(X_\sigma) - g(Y_\sigma)|.\]
Therefore it suffices to show that there exists $C > 0$ so that for every $b =(x,y)$ with $x \in \partial D$ and $y \in D$ we have
\[ |\nabla g(b)| \leq C\ol{g} \left[ \log R +  \frac{R}{r}\right].\]
Let $Y_t$ be a simple random walk started from $y$, $\tau_D = \inf\{t \geq 0 : Y_t \in \partial D\}$, and let $p_k = \p[ |Y_{\tau_D} - x| \geq k]$.  By the previous lemma, we know that there exists $C > 0$ so that $p_k \leq C / (k \wedge r)$ when $k \geq 1$.  By summation by parts,
\begin{align*}
    |\nabla g(b)|
&\leq \sum_{k=1}^{R} (p_{k-1} - p_k) k \ol{g}
  \leq \sum_{k=0}^{R-1} p_k \ol{g}
  \leq C \ol{g} \left[ \sum_{k=1}^r \frac{1}{k} + \sum_{k=r+1}^{R-1} \frac{1}{r}\right],  
\end{align*}
which proves the lemma.
\end{proof}

\section{Symmetric Random Walks}

Throughout we suppose that we have a continuous time random walk $X_t$ on $\Z^2$ with time-dependent jump rates $c_t(b)$ satisfying
\[ 0 < a \leq c_t(b) \leq A < \infty\]
uniformly in the edges $b$ and time $t$.  Let $p(s,t;x,y)$ be the transition kernel of $X$.  One of the important tools in the analysis of such walks is the Nash-Aronson estimates, and their time dependent generalization, which give a comparison between $p$ and the transition kernel of a standard random walk.

\begin{lemma}[Nash-Aronson Estimates]
\label{symm_rw::lem::nash_aronson}
Let $D \subseteq \Z^2$ and $\wt{p}(u,t;x,y)$ be the transition kernel of $\wt{X}$, the random walk in $D$ jumping with rates $c_t(b)$ stopped on its first exit from $D$.  There exists $C \geq 1, \delta > 0$ depending only on $a,A$ such that
\begin{align}
 \wt{p}(s,t;x,y) &\leq \frac{C}{1 \vee (t-s)} \exp\left( - \frac{|x-y|}{C(1 \vee (t-s)^{1/2})} \right),\\
 \wt{p}(s,t;x,y) &\geq \frac{\delta}{1 \vee (t-s)}  \text{ for } |x-y| \leq \sqrt{t-s}
\end{align}
provided $s \leq t$ with $|t-s| \leq r^{2-2\epsilon}$ and $x,y \in B(x_0,r)$ with $B(x_0,2r) \subseteq D$
\end{lemma}
\begin{proof}
In the time-independent setting with $D = \Z^2$ these follow from the usual Nash-Aronson estimates, see \cite{SZ97}.  The time-dependent extension, also for $D = \Z^2$, is proved in Propositions B3 and B4 of \cite{GOS01} and also Propositions 4.2 and 4.3 of \cite{DD05}.  

This leaves us to handle the case that $D \neq \Z^2$.  We can couple $X, \wt{X}$ together so that they are equal until $\wt{X}$ exits $D$.  By the case for $D = \Z^2$, the probability that this happens before time $s+r^{2-2\epsilon}$ given $\wt{X}_s = x \in B(x_0,r)$ is of order $\exp(-r^{\epsilon}/C)$ where $C$ is the constant of that lemma.  That is,
\[ \sup_{x \in B(x_0,r)} \sup_{y \in \Z^2} |\wt{p}(s,t;x,y) - p(s,t;x,y)| = O(\exp(- r^\epsilon/C)) \text{ for } t \leq s+r^{2-2\epsilon},\]
which implies the result.
\end{proof}

Let $\CL_t$ be the operator
\[ \CL_t f (x) = \sum_{b \ni x} c_t(b) \nabla f(b).\]

\begin{lemma}[Nash Continuity Estimate]
\label{symm_rw::lem::nash_continuity_bounded}
Let $D \subseteq \Z^2$.  Suppose that $\wt{f}_t \colon [0,\infty) \times D \to \R$ is a solution of the equation
\[ \partial_t \wt{f}_t(x) = \CL_t \wt{f}_t(x).\]
There exists $\xi_{{\rm NC}}, C > 0$ depending only on $a,A$  such that
\begin{align}
      &|\wt{f}_t(x) - \wt{f}_s(y)|\\ 
\leq& C\| \wt{f}_u\|_\infty \left[ \left( \frac{|t-s|^{1/2} \vee |x-y|}{(t \wedge s)^{1/2}} \right)^{\xi_{{\rm NC}}} + O( \exp(-r^{\epsilon}/C)) \right] \notag
\end{align}
for all $u \leq s \leq t$ with $|t-u| \leq r^{2-2\epsilon}$ and $x,y \in B(x_0,r)$ with $B(x_0,2r) \subseteq D$.
\end{lemma}
\begin{proof}
In the same manner as Lemma \ref{symm_rw::lem::nash_aronson}, this follows from the usual Nash continuity estimate for $D = \Z^2$, which is proved in the time independent setting in  \cite{SZ97} and the time dependent version is Proposition B6 of \cite{GOS01}.
\end{proof}

The following is Proposition 4.1 of \cite{DD05}.
\begin{lemma}[Caccioppoli inequality]
\label{symm_rw::lem::caccioppoli}
There exists $C > 0$ depending only on $a,A$ such that the following holds.  If $f_t$ solves $\partial_t f_t = \CL_t f_t$ on $[0,2r^2] \times B(x_0,2r)$ then
\begin{equation}
\label{symm_rw::eqn::caccioppoli}
  \int_{r^2}^{2r^2} \sum_{b \in B^*(x_0,r)} |\nabla f_t(b)|^2 dt \leq \frac{C}{r^2} \int_0^{2r^2} \sum_{y \in B(x_0,r)} f_t^2(y) dt.
\end{equation}
\end{lemma}

We now prove a Beurling-type estimate for $X_t$.

\begin{lemma}
\label{symm_rw::lem::beurling}
There exists constants $\rho_{\rm B}, C > 0$ depending only $a,A$ such that the following holds. Fix $x \in \Z^2, r > 0$ and let $H$ be a connected graph with $H \cap B(x,r)^c \neq \emptyset$ and $\dist(H,x) = d \leq r$.  If $\tau_r = \inf\{t \geq 0: |X_t - x| = r\}$ and $\tau_H = \inf \{ t \geq 0 : X_t \in H\}$ then
\[ \p[ \tau_r \leq \tau_H] \leq C \left( \frac{d}{r} \right)^{\rho_{\rm B}}.\]
\end{lemma}

The idea is to show that the probability that $X$ runs around an annulus with inner and outer radii $r$, $2r$, respectively, is strictly positive independent of $r$.  The following auxiliary estimate will be useful for the proof.

\begin{lemma}
\label{symm_rw::lem::first_hitting}
Fix $\alpha \in (0,1)$ and let $F \subseteq B(x,\alpha R)$.  Let $\tau = \inf\{t > 0 : X_t \notin B(x,R)\}$ and $\tau_F = \inf\{t > 0 : X_t \in F\}.$  There exists $C > 0$ depending only on $a,A,\alpha$ such that
\[ \p_x[\tau_F \leq \tau] \geq C \frac{|F|}{R^2}.\]
\end{lemma}
\begin{proof}
Employing Lemma \ref{symm_rw::lem::nash_aronson} in the final step,
\begin{align*}
    &  \p_x[\tau_F \leq \tau]
\geq \p_x[X_t \in F, \tau > t] = \p_x[X_t \in F] - \p_x[X_t \in F, \tau \leq t]\\
\geq& \sum_{y \in F} \left( p(0,t;x,y) - \E [ \E_{X_\tau}[p(\tau,t;X_\tau,y)] \one_{\{\tau \leq t\}}] \right)\\
\geq& \sum_{y \in F} \left( p(0,t;x,y) - \E \left[ \frac{C}{1 \vee (t-\tau)} \exp\left( - \frac{|X_\tau-y|}{C(1 \vee (t-\tau)^{1/2})} \right)\one_{\{\tau \leq t\}}\right] \right)
\end{align*}
Obviously, $|X_\tau - y| \geq (1-\alpha)R$.  It is easy to see that $x \mapsto C x^{-1} \exp(-(1-\alpha) R C^{-1} x^{-1/2})$ is maximized when $x = (1-\alpha)^2 R^2 / (4 C^2)$.  Consequently, if $t = \gamma R^2$ for $\gamma > 0$ then
\begin{align*}
       \p_x[\tau_F \leq \tau]
\geq  \sum_{y \in F} \left( \frac{\delta}{1 \vee t} - \frac{C'}{(1-\alpha)^2 R^2} \right)
\geq \left(\frac{\delta}{\gamma} - \frac{C'}{(1-\alpha)^2}\right) \frac{|F|}{R^2}
\end{align*}
where $C'$ does not depend on $\alpha,\gamma$ provided $F \subseteq B(x,\sqrt{\gamma} R)$.  This implies that there exists $\alpha_0 > 0$ such that the result holds whenever $\alpha \in (0,\alpha_0)$.  To see the general case, let $x_1,\ldots,x_n$ be an $\alpha_0/4$-net of $B(x,\alpha R)$.  If $x_j \in B(0, \tfrac{\alpha_0}{2})$, then it follows from our a priori estimate that $X_t$ hits $B(x_j, \tfrac{\alpha_0}{4})$ before exiting $B(x, R)$ with strictly positive probability $\rho_0$.  Iterating this argument, it follows that the probability that $X_t$ hits $B(x_k, \tfrac{\alpha_0}{4})$ before exiting $B(x,R)$ with strictly positive probability $\rho_1$.  The result is now follows by combining this fact once again with our a priori estimate.
\end{proof}

\begin{proof}[Proof of Lemma \ref{symm_rw::lem::first_hitting}]
The previous lemma implies the existence of $\rho_1 > 0$ depending only on $\CV$ such that the following holds.  Let $A(x,r_1,r_2) = \{y : r_1 \leq |x-y| \leq r_2\}$ be the annulus with inner and outer radii $r_1,r_2$ satisfying $r_2=2r_1$.  Let $\CC$ be the event that $X$ runs a full circle around $A(x,r_1,r_2)$ after hitting $\partial B(x, \tfrac{3}{2} r_1)$ without hitting $\partial A(x,r_1, r_2)$.  Then
\begin{equation}
\label{symm_rw::eqn::cc_bound}
 \p[\CC] \geq \rho_1 > 0.
\end{equation}
Now set $r_k = 2^k d$.  The largest index $m$ so that $r_m \leq r$ is $\lfloor \log_2 \tfrac{r}{d} \rfloor$.  By \eqref{symm_rw::eqn::cc_bound}, the probability that $X$ makes it to distance $r$ without running a full circle around one of the $A(x,r_i,r_{i+1})$ is at most 
\[ (1-\rho_1)^m \leq \exp(1 - \log_2(r/d) \rho_1) \equiv C \left( \frac{d}{r} \right)^{\rho_{\rm B}} .\]
\end{proof}

\section*{Acknowledgements}  I thank Amir Dembo for endless patience, encouragement, and inspiration without which this work would not have been possible. I also thank Scott Sheffield for suggesting the usage of Theorem \ref{harm::thm::coupling} in order to prove Theorem \ref{thm::clt} as well as Mykhaylo Shkolnikov for comments on an earlier draft of this article.  Finally, I thank Jean-Dominique Deuschel for suggesting the extension of Theorem \ref{harm::thm::coupling} to the case of non-zero tilt.

\bibliographystyle{acmtrans-ims.bst}
\bibliography{gl_level_lines}

\end{document}